\documentclass[11pt]{amsart}

\usepackage{amsfonts,amssymb,amsmath,mathrsfs,graphicx}
\usepackage{float}
\usepackage{epsfig}

\usepackage{CJK,color}

\definecolor{red}{rgb}{1.00,0.00,0.00}
\definecolor{blue}{rgb}{0.00,0.00,0.63}
\definecolor{black}{rgb}{0.00,0.00,0.00}

\newtheorem{theorem}{Theorem}[section]
\newtheorem{lemma}{Lemma}[section]

\newtheorem{proposition}{Proposition}[section]
\newtheorem{remark}{Remark}[section]

\renewcommand{\theequation}{\thesection.\arabic{equation}}

\def\charf {\mbox{{\text 1}\kern-.30em {\text l}}}
\def\G{\Gamma}

\def\mb{\mathbf}

\def\f{\frac}
\def\t{\theta}
\def\d{\delta}
\def\a{\alpha}
\def\b{\beta}
\def\i{\infty}

\def\di{\displaystyle}
\def\r{\rho}
\def\wt{\widetilde}
\def\k{\kappa}
\def\l{\lambda}

\def\o{\omega}

\def\charf {\mbox{{\text 1}\kern-.30em {\text l}}}
\def\mb{\mathbf}
\def\f{\frac}
\def\t{\theta}
\def\o{\omega}
\def\d{\delta}
\def\s{\sigma}
\def\a{\alpha}
\def\b{\beta}
\def\p{\phi}
\def\T{\Theta}
\def\i{\infty}
\def\l{\lambda}

\def\di{\displaystyle}
\def\r{\rho}
\def\wt{\widetilde}
\def\k{\kappa}

\begin{document}
\begin{CJK*}{GBK}{song}
%%%%%%%%%%%%%%%%

\title[Stability of nonlinear waves to bipolar VPB]{Stability of nonlinear wave patterns to the bipolar Vlasov-Poisson-Boltzmann system}

\author[Li]{Hailiang Li}
\address[Hailiang Li]{\newline School of Mathematical Sciences, Capital Normal University, Beijing 100048, China
\newline and Beijing Center of Mathematics
and Information Sciences, Beijing 100048, China}
\email{Hailiang-Li@mail.cnu.edu.cn}

\author[Wang]{Yi Wang}
\address[Yi Wang]{\newline CEMS, HCMS, NCMIS, Academy of Mathematics and Systems Science, Chinese Academy of Sciences, Beijing 100190, China
\newline and School of Mathematical Sciences, University of Chinese Academy of Sciences, Beijing 100049, China
\newline and Beijing Center of Mathematics
and Information Sciences, Beijing 100048, China}
\email{wangyi@amss.ac.cn}

\author[Yang]{Tong Yang}
\address[Tong Yang]{\newline Department of Mathematics, City University of Hong Kong, Hong Kong, China
}
\email{matyang@cityu.edu.hk}

\author[Zhong]{Mingying Zhong}
\address[Mingying Zhong]{\newline Department of Mathematics, Guangxi University, China
}
\email{zhongmingying@sina.com}

\date{\today}

%\subjclass{$\cdots$}


%\iffalse
%\thanks{\textbf{Acknowledgment.} The work of H. L. Li is partially supported by NNSFC grants No. 11171228, 11231006
%and 11225102, and by the Key Project of Beijing Municipal Education Commission no. CIT$\&$TCD20140323.
%The work of Y. Wang is partially supported by the National Natural Sciences Foundation
%of China (No. 11171326 and 11322106). The work of T. Yang is supported by the NSFC-RGC Grant, the General Research Fund
%of Hong Kong, CityU 103412, N-CityU102/12.  The work of M. Y. Zhong is partially supported by the NNSFC grants No. 11301094 and Project supported by Beijing Postdoctoral Research Foundation No. 2014ZZ-96.}
%\fi

\maketitle

\begin{abstract}
The main purpose of the present paper is to investigate the nonlinear stability of viscous shock waves and rarefaction wave for bipolar Vlasov-Poisson-Boltzmann (VPB) system. To this end, motivated by the micro-macro decomposition to the Boltzmann equation in \cite{Liu-Yu, Liu-Yang-Yu}, we first set up a new micro-macro decomposition around the local Maxwellian related to the bipolar VPB system and give a unified framework to study the nonlinear stability of the basic wave patterns to the system.
Then, as the applications of this new decomposition, the time-asymptotic stability of the two typical nonlinear wave patterns, viscous shock waves and rarefaction wave, are proved for the 1D bipolar Vlasov-Poisson-Boltzmann system. More precisely, it is first proved that the linear superposition of two Boltzmann shock profiles in the first and third characteristic fields is nonlinearly stable to the 1D bipolar VPB system up to some suitable shifts without the zero macroscopic mass conditions on the initial perturbations. Then the time-asymptotic stability of rarefaction wave fan to compressible Euler equations is proved to 1D bipolar VPB system.
These two results are concerned with the nonlinear stability of wave patterns for Boltzmann equation coupled with additional (electric) forces, which together with spectral analysis made in \cite{LYZ1} sheds light on understanding the complicated dynamic behaviors around the wave patterns in the transportation of charged particles under the binary collisions, mutual interactions, and the effect of the electrostatic potential forces.
\end{abstract}

\maketitle \centerline{\date}

\tableofcontents

%%%%%%%%%%%%%%%%%%%%%%%%%%%%%%%%%%%%%%%%%%%%%%%%%%%%%%%%%%%%%%%%%%%
%
%                             Sec 1. Introduction
%
%%%%%%%%%%%%%%%%%%%%%%%%%%%%%%%%%%%%%%%%%%%%%%%%%%%%%%%%%%%%%%%%%%%
\section{Introduction}
\setcounter{equation}{0}
It is an interesting and challenging problem to investigate the nonlinear wave phenomena and understand the dynamical behaviors of charged particles transport under the influence of external force such as electrostatic potential, magnetic field, or electromagnetic fields, etc.  To begin with, we first investigate wave phenomena for the bipolar  Vlasov-Poisson-Boltzmann system which is used to simulate the  transport of two dilute charged particles (e.g. ions and electrons) affected by the self-consistent electrostatic potential force~\cite{Ma}.  In spatial three-dimensional space, the bipolar Vlasov-Poisson-Boltzmann system takes the form
\begin{equation}\label{VPB0}
\left\{
\begin{array}{ll}
\di
F_{At} +v\cdot\nabla_x F_A+\nabla_x\Pi\cdot\nabla_v F_A= Q(F_A, F_A)+Q(F_A,F_B), \\[2mm]
\di
F_{Bt} +v\cdot\nabla_x F_B-\nabla_x\Pi\cdot\nabla_v F_B= Q(F_B,F_A)+Q(F_B,F_B), \\
\di \Delta \Pi=\int (F_A-F_B)dv,
\end{array}
\right.
\end{equation}
where $v=(v_1,v_2,v_3)\in {\mb R}^3$, $x=(x_1,x_2,x_3)\in {\mb R}^3$, $t\in \mb{R}^+$ and $F_A=F_A(t,x,v), F_B=F_B(t,x,v)$ are the
density distribution function of two-species particles (e.g. ions and electrons) at time-space $(t,x)$ with velocity $v$, respectively, and $\Pi=\Pi(x,t)$ is the electric field potential. For the hard sphere model, the collision operator
$Q(f,g)$ takes the following  bilinear form
$$
 Q(f,g)(v) \equiv \f{1}{2}
\int_{{\mb R}^3}\!\!\int_{{\mb S}_+^2} \Big(f(v')g(v_*')-f(v) g(v_*) \Big)
|(v-v_*)\cdot\Omega)|
 \; d v_* d\Omega,
$$
where the unit vector $\Omega \in {\mb S}^2_+=\{\Omega\in {\mb S}^2:\ (v-v_*)\cdot \Omega\geq 0\}$, $(v,v_*)$ and $(v',v_*')$ are the two particle velocities before and after the binary elastic collision respectively, which together with the conservation laws of momentum and energy, satisfy the following relations
$$
 v'= v -[(v-v_*)\cdot \Omega] \; \Omega, \qquad
 v_*'= v_* + [(v-v_*)\cdot \Omega] \; \Omega, \quad \Omega \in {\mb S}^2.
$$

In the case that the effect of electrons is neglected, the bipolar Vlasov-Poisson-Boltzmann (abbreviated as VPB for simplicity below) system~\eqref{VPB0} can be reduced to the unipolar VPB equations~\cite{Ma}
\begin{equation}\label{VPB0.z}
\left\{
\begin{array}{ll}
\di
F_{At} +v\cdot\nabla_x F_A+\nabla_x\Pi\cdot\nabla_v F_A= Q(F_A, F_A), \\[2mm]
\di \Delta \Pi=\int F_Adv - \bar\rho(x),
\end{array}
\right.
\end{equation}
with $\bar\rho(x)>0$ a given function representing the background doping profile.

Both the bipolar VPB system~\eqref{VPB0} and unipolar VPB system can be viewed as the Boltzmann equations under the affect of the electrostatic potential force determined through the self-consistent Poisson equation related to the macroscopic density of charged particles.
However, this bipolar or unipolar VPB system is by no means the simple extension of the Boltzmann equation. Indeed, there are complicated asymptotical behaviors, which are different completely from those of the Boltzmann equation,  had been observed and justified rigorously for the VPB systems~\eqref{VPB0}--\eqref{VPB0.z} in \cite{DS,LYZ,LYZ1} due to the combined effects such as the binary elastic collision between the particles of same species, the electrostatic potential force, and/or the mutual interactions among the charged particles of two different species. To be more precise,  it was shown in  \cite{LYZ} that for the unipolar VPB system for motion of one species, the influence of the electric field affects the spectrum structure of the linearized VPB system and causes the slower but optimal (compared with the Boltzmann equation) time-asymptotical convergence rate of global solution to the equilibrium state, and there is no wave pattern propagation (such as the shock profile and rarefaction wave compared with the Boltzmann equation)  due to effect of the electric field~\cite{LYZ1}, one can refer to \cite{DS,LYZ,LYZ1,Zhong} and references therein for more details.
On the other hand, however, a completely different dynamical phenomena/behaviors of global solution are observed for bipolar VPB system in \cite{LYZ1}. Therein, it was shown that the linearized VPB system around the global Maxwellian consists of
a decoupled system: one is the linear Boltzmann equation for the distribution function  $F_1=\f{F_A+F_B}{2}$ of the average of the total charged particles which admits the wave modes at lower frequency,  the other is the equation of unipolar VPB type for the neutral function $F_2=\f{F_A-F_B}{2}$ of the particles with different charge which admits spectral gap at lower frequency and causes the strong neutrality in the sense that the neutral function $F_2$  and the electric field related the neutral  function $F_2$ decay exponentially in time. In addition, the multi-dimensional pointwise diffusive properties similar to the Boltzmann equation is also shown in \cite{LYZ1}. A natural problem follows then, can one observe the nonlinear wave pattern propagation and justify the combined influence of the electrostatic potential force and/or the mutual interactions among the charged particles for the bipolar VPB system?

\bigskip

The main purpose of the present paper is to investigate the nonlinear wave phenomena and understand the dynamical behaviors of charged particles transport under the influence of the electrostatic potential force.  It is well-known that the Boltzmann equation is asymptotically equivalent to the compressible Euler equations as illustrated by the famous Hilbert expansions. The system of compressible Euler equations is a typical example of hyperbolic conservation laws system. There are three basic wave patterns to the hyperbolic conservation laws: two nonlinear waves, i.e., the shock wave and rarefaction wave in the genuinely nonlinear field,  and one linearly degenerate wave called contact discontinuity. Therefore, Boltzmann equation has the rich wave phenomena as for the macroscopic fluid mechanics, and there have been made important progress on the nonlinear stability of these basic wave patterns of the Boltzmann equation, refer for instance to~\cite{Caflish-Nicolaenko, Liu-Yu, Liu-Yu-1, Liu-Yang-Yu-Zhao, Huang-Xin-Yang, Huang-Yang, Yu, Wang-Wang} and references therein. Yet, we should mention here that the pioneering study on the stability and positivity of viscous shock wave was first made by Liu-Yu \cite{Liu-Yu} in energy space with the zero total macroscopic mass condition based on the micro-macro decomposition proposed by Liu-Yu~\cite{Liu-Yu}. Furthermore, Yu \cite{Yu} made an important breakthrough to establish the stability of single viscous shock profile without the zero mass condition by the elegant point-wise method based on the Green function around the shock profile. Then, the stability of rarefaction wave is proved by  Liu-Yang-Yu-Zhao \cite{Liu-Yang-Yu-Zhao} and the stability of viscous contact wave, which is the viscous version of contact discontinuity, by Huang-Yang \cite{Huang-Yang} with the zero mass condition and Huang-Xin-Yang \cite{Huang-Xin-Yang} without the zero mass condition.  Recently, Wang-Wang \cite{Wang-Wang} proved the stability of superposition of two viscous shock profiles to the Boltzmann equation without the zero mass condition by the weighted characteristic energy method.

% {\red Netherverness, the hydrodynamic limit of Boltzmann equation to the generic Riemann solution to Euler system is investigated successfully, refer to \cite{Yu1, Huang-Wang-Yang, Xin-Zeng, Huang-Wang-Yang-2, Huang-Wang-Yang-3, hwwy} and references therein.}\\
%
%
\bigskip

Therefore, due to the appearance of wave modes and the spectral gap of the linearized bipolar VPB system as shown in \cite{LYZ1}, it is natural and interesting to investigate and understand the nonlinear wave phenomena of the bipolar VPB system under the influence of electric field and mutual interactions between charged particles.
To this end, we first consider the nonlinear stability of viscous shock waves and rarefaction wave for the bipolar VPB system. Yet, compared with the stability analysis made for the Boltzmann equation, it is not straightforward to study the stability of viscous shock waves and rarefaction wave under the influence of the electric field effect and the mutual interactions among the charged particles. Moreover, there is no generic framework made concerned with the stability of basic wave patterns to the bipolar VPB system as far as we know.  To overcome these difficulties, for $F_1=\f{F_A+F_B}{2}$ satisfying the Boltzmann-type equation with the additional electric fields, we employ the micro-macro type decomposition as in \cite{Liu-Yu, Liu-Yang-Yu} for Boltzmann equation, while in VPB type equation for $F_2=\f{F_A-F_B}{2}$, we introduce a new micro-macro type decomposition around the local Maxwellian with respect to $F_1$. More importantly, we can derive a new diffusion equation with the damping from the macroscopic part of $F_2$, which crucially implies that the electric fields are strongly dissipative and guarantees the stability of wave patterns.

Note that this new decomposition for $F_2$ is quite universal and will play an important role in the stability analysis towards wave patterns to the bipolar VPB system \eqref{VPB0}. Then as the applications of this new decomposition, the stability of viscous shock waves and rarefaction wave are proved for the 1D bipolar VPB system as the first step. Note that for the stability of superposition of two viscous shock waves in the first and third characteristic fields, there is no zero macroscopic mass conditions for the initial perturbations by introducing suitable shifts on the two viscous shock waves, linear diffusion wave in the second characteristic field and the coupled diffusion waves, which is motivated by Liu \cite{Liu-1985}, Szepessy-Xin \cite{Szepessy-Xin}, Huang-Matsumura \cite{Huang-matsumura} and Wang-Wang \cite{Wang-Wang}. Roughly speaking, it is first proved that the linear superposition of two Boltzmann shock profiles is nonlinearly stable time-asymptotically to the 1D bipolar VPB system up to some suitable shifts without imposing the zero macroscopic conditions on the initial perturbations. Moreover, we proved the nonlinear stability of the rarefaction wave solution to the Riemann problem of inviscid Euler system time-asymptotically to the 1D bipolar VPB system. The precise statements of the stability of viscous shock waves and rarefaction waves can be referred to Theorem \ref{thm0} and Theorem \ref{thm}, respectively. Future works will be done for the stability of other wave patterns and their linear superpositions.

\bigskip

There have been important works on the existence and behavior of solutions to the VPB system. The global existence of renormalized solution for large initial data was proved in Mischler \cite{Mi}. The first global
existence result on classical solution in torus when the initial data is near a global Maxwellian was established
in Guo \cite{Guo1}. And the global existence of classical solution in $\mathbf{R}^3$
was given in \cite{YYZ, YZ} . The case
with general stationary background density function
 was studied in \cite{DY}, and the perturbation of vacuum
was investigated in \cite{DYZ, DS}. Recently, Li-Yang-Zhong \cite{LYZ, LYZ1}  analyze the
spectrum of the linearized VPB system (unipolar and bipolar) and
obtain the optimal decay rate of solutions to the nonlinear system near global Maxwellian. See also the works on the stability of global Maxwellian and the optimal time decay rate in \cite{YY, W} and on boundary value problems of stationary VPB system \cite{BGV}. Recently, Duan-Liu \cite{Duan-Liu} proved the stability of rarefaction wave to a unipolar VPB system, which can be viewed as an approximation of bipolar VPB system  \eqref{VPB0} when the electron density is very rarefied and reaches a local equilibrium state with small electron mass compared with the ion. However, there is no any analysis made concerned with the stability of  viscous shock waves to the bipolar VPB system \eqref{VPB0} as far as we know.

\bigskip

It should be also mentioned that deep investigation has been achieved on the asymptotic stability of wave patterns for viscous conservation laws, which are extremely  helpful to understand the wave phenomenon to the kinetic equations. The time-asymptotic stability of viscous shock profile started from Goodman \cite{Goodman} for the uniformly viscous conservation laws and Matsumura-Nishihara \cite{MN-85} for the compressible Navier-Stokes equations independently by the anti-derivative methods. Note that in the both above results the zero mass conditions are imposed on the initial perturbation. Then Liu \cite{Liu-1985} and Szepessy-Xin \cite{Szepessy-Xin} removed the zero mass condition by introducing the linear and nonlinear diffusion waves and the coupled diffusion waves in the transverse characteristic field for the uniformly viscous conservation laws and Liu-Zeng \cite{Liu-Zeng} proved the physical viscosity case. And Zumbrun
\cite{zum} proved the stability of large-amplitude shock waves of compressible Navier-Stokes equations by Evans function approach.
Then the stability of rarefaction waves for the compressible Navier-Stokes are proved by Matsumura-Nishihara \cite{MN-86} and Nishihara-Yang-Zhao \cite{NYZ}. The stability of viscous contact wave for the uniformly viscous conservation laws was proved by Liu-Xin \cite{Liu-Xin} and Xin \cite{Xin2} with the zero mass condition. Then for the compressible Navier-Stokes equations with physical viscosities, the stability of viscous contact wave was proved by Huang-Matsumura-Xin \cite{Huang-Matsumura-Xin} with the zero mass condition, and by Huang-Xin-Yang \cite{Huang-Xin-Yang} without the zero mass condition.
For the composite waves, Huang-Matsumura \cite{Huang-matsumura}
first studied the asymptotic stability of two viscous shock waves under
general initial perturbation without zero mass conditions on initial perturbations for the full Navier-Stokes system and
Huang-Li-Matsumura \cite{Huang-Li-matsumura} justified the stability of a combination wave of a viscous contact wave and rarefaction waves.

\bigskip

The rest part of the paper is arranged as follows. In section 2 we present the classical micro-macro decomposition and introduce new micro-macro decomposition for the bipolar VPB system \eqref{VPB0}. Then, the main result on the stability of  viscous shock waves and rarefaction wave are stated and proved in section 3 and 4, respectively. Finally, Appendix A and B are devoted to a priori estimates for the stability of viscous shock waves and rarefaction wave, respectively.

%%%%%%%%%%%%%%%%%%%%%%%%%%%%%%%%%%%%%%%%%%%%%%%%%%%
%
%
%
%
%%%%%%%%%%%%%%%%%%%%%%%%%%%%%%%%%%%%%%%%%%%%%%%%%%%

\section{Micro-Macro Decompositions}
%\label{decompisition}
\setcounter{equation}{0}
We reformulate the bipolar VPB system \eqref{VPB0} and give a new micro-macro decomposition around the local Maxwellian in order to study the nonlinear stability of basic wave patterns to the system \eqref{VPB0}.  Set
$$
F_1=\f{F_A+F_B}{2},\qquad F_2=\f{F_A-F_B}{2},
$$
then the system~\eqref{VPB0} is changed into
\begin{equation}\label{VPB}
\left\{
\begin{array}{ll}
\di
F_{1t} +v\cdot\nabla_x F_1+\nabla_x\Pi\cdot\nabla_v F_2=2 Q(F_1,F_1), \\
\di
F_{2t} +v\cdot\nabla_x F_2+\nabla_x\Pi\cdot\nabla_v F_1=2 Q(F_1,F_2), \\
\di \Delta \Pi=2\int F_2dv.
\end{array}
\right.
\end{equation}
We present the micro-macro decompositions around the local Maxwellian to the bipolar VPB system \eqref{VPB}.
The equation $\eqref{VPB}_1$ can be viewed as the Boltzmann equation with additional electric potential force, we make use of the micro-macro decomposition as introduced by Liu-Yu \cite{Liu-Yu} and Liu-Yang-Yu \cite{Liu-Yang-Yu}. In fact, for any solution $F_1(t,x,v)$ to Eq.~$\eqref{VPB}_1$, there are five macroscopic (fluid) quantities: the mass density $\r(t,x)$, the momentum $m(t,x)=\r u(t,x)$,
and the total energy $E(t,x)=\r\big(e+\frac12|u|^2\big)(t,x)$ defined by
\begin{equation}
\left\{
\begin{array}{l}
\di\rho(t,x)=\int_{\mb{R}^3}\xi_0(v)F_1(t,x,v)dv,\\
\di\rho u_i(t,x)=\int_{\mb{R}^3}\xi_i(v)F_1(t,x,v)dv,~i=1,2,3,\\
\di\rho(e+\f{|u|^2}{2})(t,x)=\int_{\mb{R}^3}\xi_4(v)F_1(t,x,v)dv,
\end{array}
\right. \label{macro}
\end{equation}
where  $\xi_i(v)$ $(i=0,1,2,3,4)$ are the collision
invariants given by
\begin{equation}
 \xi_0(v) = 1,~~~
 \xi_i(v) = v_i~  (i=1,2,3),~~~
 \xi_4(v) = \f{1}{2} |v|^2,
\label{collision-invar}
\end{equation}
and satisfy
$$
\int_{{\mb R}^3} \xi_i(v) Q(g_1,g_2) d v =0,\quad {\textrm
{for} } \ \  i=0,1,2,3,4.
$$
Define the local Maxwellian $\mb{M}$ associated to the solution $F_1(t,x,v)$  to Eq.~$\eqref{VPB}_1$  in terms of the fluid quantities by
\begin{equation}
\mb{M}:=\mb{M}_{[\rho,u,\t]} (t,x,v) = \f{\rho(t,x)}{\sqrt{ (2 \pi
R \t(t,x))^3}} e^{-\f{|v-u(t,x)|^2}{2R\t(t,x)}},  \label{maxwellian}
\end{equation}
where $\t(t,x)$ is the temperature which is related to the internal
energy $e(t,x)$ by $e=\frac{3}{2}R\t$ with $R>0$ the gas constant,  and $u(t,x)=\big(u_1(t,x),u_2(t,x),u_3(t,x)\big)^t$
is the fluid velocity. Then, the collision operator of $Q(f,f)$ can be linearized to be $\mb{L}_\mb{M}$ with respect to the local Maxwellian $\mb{M}$ by
\begin{equation}
\mb{L}_\mb{M} g=2Q(\mb{M}, g)+ 2Q(g,\mb{M}).  \label{L_M}
\end{equation}
The null space $\mathfrak{N}_1$ of $\mb{L}_\mb{M}$ is spanned by $\xi_i(v)~(i=0,1,2,3,4)$.

Define an inner product $ \langle g_1,g_2\rangle_{\widetilde{\mb{M}}}$ for $g_i\in L(\mb{R}^3_v)$ with respect to the given local Maxwellian $\wt{\mb{M}}$ as:
\begin{equation}
 \langle g_1,g_2\rangle_{\widetilde{\mb{M}}}\equiv \int_{{\mb R}^3}
 \f{1}{\widetilde{\mb{M}}}g_1(v)g_2(v)d v.\label{product}
\end{equation}
For simplicity, if $\widetilde{\mb{M}}$ is the local Maxwellian $\mb{M}$ in \eqref{maxwellian}, we shall use the notation $\langle\cdot,\cdot\rangle$ instead of $\langle\cdot,\cdot\rangle_{\mb{M}}$.
Furthermore, there exists a positive constant $\wt\sigma_1>0$ such
that it holds for any function $g(v)\in \mathfrak{N}_1^\bot$ (cf. \cite{CC,Grad}) that
\begin{equation}
\langle g,\mb{L}_\mb{M}g\rangle
\le
 -\wt\sigma_1\langle \nu(|v|)g,g\rangle ,  \label{H-thm}
\end{equation}
where $\nu(|v|)\sim (1+|v|)$ is the collision frequency for the
hard sphere collision.

  With respect to the inner product $\langle\cdot,\cdot\rangle$,  the following pairwise orthogonal basis span the macroscopic space  $\mathfrak{N}_1$
\begin{equation}
\left\{
\begin{array}{l}
 \chi_0(v) \equiv {\di\f1{\sqrt{\rho}}\mb{M}}, \quad
 \chi_i(v) \equiv {\di\f{v_i-u_i}{\sqrt{R\t\rho}}\mb{M}} \ \ {\textrm {for} }\ \  i=1,2,3, \\[2mm]
 \chi_4(v) \equiv
 {\di\f{1}{\sqrt{6\rho}}(\f{|v-u|^2}{R\t}-3)\mb{M}},\quad \langle\chi_i,\chi_j\rangle=\delta_{ij}, ~i,j=0,1,2,3,4.
 \end{array}
\right.\label{orthogonal-base}
\end{equation}
In terms of above orthogonal  basis, the macroscopic projection $\mb{P}_0$ from $L^2(\mb{R}^3_v)$ to  $\mathfrak{N}_1$ and the microscopic projection $\mb{P}_1$  from $L^2({\mb R}^3_v)$ to $\mathfrak{N}_1^\bot$ can be defined as
\begin{equation*}
 \mb{P}_0g = {\di\sum_{j=0}^4\langle g,\chi_j\rangle\chi_j},\qquad
 \mb{P}_1g= g-\mb{P}_0g.
\end{equation*}
A function $g(v)$ is said to be microscopic or non-fluid, if it holds
$$
\int g(v)\xi_i(v)dv=0,~i=0,1,2,3,4,
$$
where $\xi_i(v)$ are the collision invariants defined in \eqref{collision-invar}.

Based on the above preparation, the solution $F_1(t,x,v)$  to Eq.~$\eqref{VPB}_1$ can be decomposed into
the macroscopic (fluid) part, i.e., the local Maxwellian $\mb{M}=\mb{M}(t,x,v)$
defined in \eqref{maxwellian}, and the microscopic (non-fluid) part, i.e.
$\mb{G}=\mb{G}(t,x,v)$:
$$
F_1(t,x,v)=\mb{M}(t,x,v)+\mb{G}(t,x,v),\quad
\mb{P}_0F_1=\mb{M},~~~\mb{P}_1F_1=\mb{G},
$$
and the Eq.~$\eqref{VPB}_1$ becomes
\begin{equation}
(\mb{M}+\mb{G})_t+v\cdot\nabla_x(\mb{M}+\mb{G})+\nabla_x\Pi\cdot\nabla_vF_2
=\mb{L}_\mb{M}\mb{G}+2Q(\mb{G},\mb{G}). \label{M+G}
\end{equation}
Taking the inner product of the equation \eqref{M+G} and the
collision invariants $\xi_i(v)$ $(i=0,1,2,3,4)$  with respect to
$v$ over ${\mb R}^3$, one has the following  system for the fluid
variables $(\rho, u, \theta)$:
\begin{equation}
\left\{
\begin{array}{lll}
\di \rho_{t}+{\rm div}_x (\rho u)=0, \\
\di (\rho u)_t+{\rm div}_x  (\rho u\otimes u)
+\nabla_x  p- n_2 \nabla_x \Pi =-\int v\otimes v\cdot\nabla_x\mb{G}dv,  \\
\di [\rho(e+\f{|u|^2}{2})]_t+{\rm div}_x  [\rho
u(e+\f{|u|^2}{2})+pu]-\nabla_x \Pi\cdot\int
vF_2dv=-\int\f12|v|^2v\cdot\nabla_x\mb{G}dv,
\end{array}
\right. \label{ns}
\end{equation}
where
\begin{equation}\label{macro-n2}
n_2=n_2(x,t)=\int F_2(x,t,v) dv.
\end{equation}
Yet, the above fluid-type system~\eqref{ns} is not self-contained and the equation for the microscopic component ${\mb{G}}$
is needed, which can be  derived by applying the projection operator $\mb{P}_1$ into  Eq.~\eqref{ns}:
\begin{equation}
\mb{G}_t+\mb{P}_1(v\cdot\nabla_x\mb{M})+\mb{P}_1(v\cdot\nabla_x\mb{G})+\mb{P}_1(\nabla_x\Pi\cdot\nabla_vF_2)
=\mb{L}_\mb{M}\mb{G}+ 2Q(\mb{G}, \mb{G}). \label{F1-non-fluid}
\end{equation}
Recall that the linearized collision operator $\mb{L}_\mb{M}$ defined by \eqref{L_M} is dissipative on $\mathfrak{N}_1^\bot$, and its inverse
$\mb{L}_\mb{M}^{-1}$  is a bounded operator on $\mathfrak{N}_1^\bot$. Thus, it follows from \eqref{F1-non-fluid} that
\begin{equation}
\mb{G}= \mb{L}_\mb{M}^{-1}[\mb{P}_1(v\cdot\nabla_x\mb{M})] +\Gamma
\label{(1.12)}
\end{equation}
with
\begin{equation}\label{Pi}
\Gamma=\mb{L}_\mb{M}^{-1}[\mb{G}_t+\mb {P}_1(v\cdot\nabla_x\mb{G})+\mb
{P}_1(\nabla_x\Pi\cdot\nabla_vF_2)-2Q(\mb{G}, \mb{G})].
\end{equation}
Substituting  \eqref{(1.12)} into \eqref{ns}, we finally obtain the compressible Navier-Stokes-type equations for  the macroscopic fluid quantities $(\rho,u,\theta)$
\begin{equation}\label{F1-fluid}
\left\{
\begin{array}{l}
\di \rho_{t}+{\rm div}_x (\rho u)=0,\\
\di (\rho u)_t+{\rm div}_x  (\rho u\otimes u) +\nabla_x  p-n_2\nabla_x\Pi\\
\di \qquad\qquad\quad =-\int v\otimes v\cdot\nabla_x\big( \mb{L}_\mb{M}^{-1}[\mb{P}_1(v\cdot\nabla_x\mb{M})] \big)dv-\int v\otimes v\cdot\nabla_x\Gamma dv,  \\
\di [\rho(\t+\f{|u|^2}{2})]_t+{\rm div}_x[\rho
u(\t+\f{|u|^2}{2})+pu]-\nabla_x\Pi\cdot\int vF_2dv\\
\di\qquad\qquad\qquad  =-\int\f12|v|^2v\cdot\nabla_x\big( \mb{L}_\mb{M}^{-1}[\mb{P}_1(v\cdot\nabla_x\mb{M})] \big) dv
-\int\f12|v|^2v\cdot\nabla_x\Gamma dv,
\end{array}
\right.
\end{equation}
A direct computation gives rise to
$$
-\int v_i v_j\cdot\nabla_{x_j}\big( \mb{L}_\mb{M}^{-1}[\mb{P}_1(v\cdot\nabla_x\mb{M})] \big)dv=\sum_{j=1}^3\Big[\mu(\theta)(u_{ix_j}+u_{jx_i}-\frac23\delta_{ij}{\rm div}_x u)\Big]_{x_j},
$$
and
$$
\begin{array}{ll}
\di -\int\f12|v|^2v\cdot\nabla_x\big( \mb{L}_\mb{M}^{-1}[\mb{P}_1(v\cdot\nabla_x\mb{M})] \big) dv=\sum_{j=1}^3(\k(\t)\t_{x_j})_{x_j}+\sum_{i,j=1}^3\Big\{\mu(\theta)u_i(u_{ix_j}+u_{jx_i}-\frac23\delta_{ij}{\rm div}_x u)\Big\}_{x_j},
\end{array}
$$
where the viscosity coefficient $\mu(\t)>0$ and the heat
conductivity coefficient $\k(\t)>0$ are smooth functions of the
temperature $\t$. Here,  we renormalize the gas constant $R$ to be
$\f{2}{3}$ so that $e=\t$ and $p=\f23\rho\t$.

Now we decompose $F_2$ in the equation $\eqref{VPB}_2$, which is one of main contributions of the present paper.  Roughly speaking, the macroscopic part of $F_2$ satisfies the diffusive equation with damping term and the microscopic couplings, which ensure the strong dissipation of the electric forces and further guarantee the stability of wave patterns.  More precisely, the equation $\eqref{VPB}_2$ is a system of Vlasov-Poisson-Boltzmann type, which in virtue of the decomposition $F_1=\mb{M}+\mb{G}$ becomes
\begin{equation}
F_{2t}+v\cdot\nabla_xF_2+\nabla_x\Pi\cdot\nabla_vF_1
 =\mb{N}_\mb{M} F_2+2Q(F_2,\mb{G}),             \label{F2}
\end{equation}
with the linearized operator $\mb{N}_\mb{M}$ defined by
$$
\mb{N}_{\mb{M}} h=2Q(h,\mb{M}).
$$
The null space $\mathfrak{N}_2$ of $\mb{N}_\mb{M}$ is
spanned by the single macroscopic variable:
$$
\chi_0(v)=\f{\mb{M}}{\sqrt\rho},
$$
which is totally different from  the linearized operator $\mb{L}_{\mb{M}}$ due to the quite different collision structures, whose null space $\mathfrak{N}_1$ is spanned by five macroscopic variables $\chi_j(v)~(j=0,1,2,3,4),$ one can refer Sotirov-Yu  \cite{AY} for the delicate analysis of the structure of the
linearized operator $\mb{N}_{\mb{M}}$ for the gas mixture without the electric effects. Furthermore, there exists a positive constant $\wt\sigma_2>0$ such
that it holds for any function $g(v)\in \mathfrak{N}_2^\bot$ (cf. \cite{AY}) that
$$
\langle g,\mb{N}_\mb{M}g\rangle \le -\wt\sigma_2\langle
\nu(|v|)g,g\rangle ,
$$
where $\nu(|v|)\sim(1+|v|)$ is the collision frequency for the hard sphere collision.
Consequently, the linearized collision operator $\mb{N}_\mb{M}$ is dissipative on $\mathfrak{N}_2^\bot$, and its inverse
$\mb{N}_\mb{M}^{-1}$ is a bounded operator on $\mathfrak{N}_2^\bot$.

Then we introduce a new micro-macro decomposition around the local Maxellian $\mb{M}(x,t,v)$ associated with $F_1$  as follows:
$$
\mb{P}_dg=\langle g,\mb{M}\rangle\frac{\mb{M}}{\rho},\qquad \mb{P}_c g=g-\mb{P}_dg,\qquad \forall g.
$$
Then the solution $F_2(t,x,v)$  to the VPB equation  $\eqref{VPB}_2$ can be decomposed into
\begin{equation}\label{F2-New}
F_2(t,x,v)=\f{\mb{M}}{\r}n_2+\mb{P}_c F_2,\quad
\mb{P}_dF_2=\f{\mb{M}}{\r}n_2,~~~\mb{P}_cF_2=F_2-\f{\mb{M}}{\r}n_2.
\end{equation}
Taking the inner product of the equation \eqref{F2} and the
collision invariants $\xi_0(v)=1$ with respect to $v$ over $\mb{R}^3$, one has the following conservation law:
\begin{equation}\label{F2n2}
\di n_{2t}+{\rm div}_x (\int v F_2dv)=0. \\
\end{equation}
Substituting the micro-macro decomposition of $F_2$ in
\eqref{F2-New} into the above conservation law yields
\begin{equation}\label{F2n2-1}
\di n_{2t}+{\rm div}_x (un_2)+{\rm div}_x (\int v \mb{P}_cF_2dv)=0.
\end{equation}
Applying the projection operator $\mb{P}_c$ to the equation
\eqref{F2}, one has the non-fluid part equation of $F_2$:
\begin{equation}\label{F2-pc}
\partial_t(\mb{P}_c F_2)-\mb{N}_{\mb{M}}(\mb{P}_c F_2)+\mb{P}_c(v\cdot\nabla_x F_2)+\mb{P}_c (\nabla_x\Pi\cdot\nabla_vF_1)+(\f{\mb{M}}{\rho})_t~n_2=2Q(F_2,\mb{G}),
\end{equation}
where one has used the facts that
$$
\mb{P}_c(\partial_tF_2)=\partial_t(\mb{P}_c F_2)+(\f{\mb{M}}{\rho})_t~n_2,
\quad
\mb{N}_{\mb{M}}( F_2)=\mb{N}_{\mb{M}}(\mb{P}_c F_2).
$$
By \eqref{F2-pc} and continuity of the inverse operator $\mb{N}_\mb{M}^{-1}$ on $\mathfrak{N}_2^\bot$, the $\mb{P}_c F_2$ can be expressed by
\begin{equation}\label{F2-pc1}
\mb{P}_c F_2 = \mb{N}_{\mb{M}}^{-1}\Big[\partial_t(\mb{P}_c
F_2)+\mb{P}_c(v\cdot\nabla_x F_2)+\mb{P}_c
(\nabla_x\Pi\cdot\nabla_vF_1)+(\f{\mb{M}}{\rho})_t~n_2-2Q(F_2,\mb{G})\Big],
\end{equation}
and it follows
\begin{equation}\label{m2}
\begin{array}{ll}
\di {\rm div}_x (\int v \mb{P}_cF_2dv)={\rm div}_x \Big(\int v
\mb{N}_{\mb{M}}^{-1}\Big[\partial_t(\mb{P}_c
F_2)+\mb{P}_c(v\cdot\nabla_x F_2)+\mb{P}_c
(\nabla_x\Pi\cdot\nabla_vF_1)\\
\di \qquad\qquad\qquad\qquad\qquad\qquad
+(\f{\mb{M}}{\rho})_t~n_2-2Q(F_2,\mb{G})\Big]dv\Big).
\end{array}
\end{equation}
The second and third terms on the right hand side of \eqref{m2} can be further transformed as follows. By the decomposition $F_2=\f{\mb{M}}{\r}n_2+\mb{P}_c F_2$,  it holds
\begin{equation}\label{m21}
\begin{array}{ll}
\di {\rm div}_x \Big(\int v \mb{N}_{\mb{M}}^{-1}\big[\mb{P}_c(v\cdot\nabla_x F_2)\big]\Big) ={\rm div}_x \Big(\nabla_{x_j}(\f{n_2}{\r})\int v \mb{N}_{\mb{M}}^{-1}\big[\mb{P}_c(v_j\mb{M}) )\big]\Big)\\
\di \qquad\qquad\qquad\qquad  +{\rm div}_x \Big(\f{n_2}{\r}\int v \mb{N}_{\mb{M}}^{-1}\big[\mb{P}_c(v\cdot\nabla_x\mb{M}) )\big]\Big)+{\rm div}_x \Big(\int v \mb{N}_{\mb{M}}^{-1}\big[\mb{P}_c(v\cdot\nabla_x (\mb{P}_cF_2))\big]\Big)\\
\di =-{\rm div}_x \Big(\k_1(\t)\nabla_x(\f{n_2}{\r})\Big)+{\rm div}_x
\Big(\f{n_2}{\r}\int v
\mb{N}_{\mb{M}}^{-1}\big[\mb{P}_c(v\cdot\nabla_x\mb{M})
)\big]\Big)+{\rm div}_x \Big(\int v
\mb{N}_{\mb{M}}^{-1}\big[\mb{P}_c(v\cdot\nabla_x
(\mb{P}_cF_2))\big]\Big),
\end{array}
\end{equation}
where one has used the fact
$$
\int v_i \mb{N}_{\mb{M}}^{-1}\big[\mb{P}_c(v_j\mb{M})\big]=\int v_i \mb{N}_{\mb{M}}^{-1}\big[(v_j-u_j)\mb{M}) )\big]=\int (v_i-u_i) \mb{N}_{\mb{M}}^{-1}\big[(v_j-u_j)\mb{M}) )\big]=-\k_1(\t)\delta_{ij},
$$
with $i,j=1,2,3$ and $\k_1(\t)>0$ being a smooth function of the temperature $\t$.\\
On the other hand, by the decomposition $F_1=\mb{M}+\mb{G}$ it holds
\begin{equation}\label{m22}
\begin{array}{ll}
\di {\rm div}_x \Big(\int v \mb{N}_{\mb{M}}^{-1}\big[\mb{P}_c(\nabla_x\Pi\cdot\nabla_vF_1)\big]\Big)={\rm div}_x \Big(\int v \mb{N}_{\mb{M}}^{-1}\big[\nabla_x\Pi\cdot\nabla_vF_1\big]\Big)\\
\di =-{\rm div}_x \Big(\int v \mb{N}_{\mb{M}}^{-1}\big[\nabla_x\Pi\cdot\f{v-u}{R\t}\mb{M})\big]\Big)+{\rm div}_x \Big(\int v \mb{N}_{\mb{M}}^{-1}\big[\nabla_x\Pi\cdot\nabla_v\mb{G}\big]\Big)\\
\di ={\rm div}_x \Big(\f{\k_1(\t)}{R\t}\nabla_x\Pi\Big)+{\rm div}_x \Big(\int v \mb{N}_{\mb{M}}^{-1}\big[\nabla_x\Pi\cdot\nabla_v\mb{G}\big]\Big).
\end{array}
\end{equation}
Substituting \eqref{m2} into \eqref{F2n2-1} and making use of \eqref{m21} and \eqref{m22}, we can derive the governing equation on the fluid-part of $F_2$
\begin{equation}\label{F2n2-2}
\begin{aligned}
&\di n_{2t}+{\rm div}_x(un_2)+{\rm div}_x \Big(\f{\k_1(\t)}{R\t}\nabla\Pi\Big)
 -{\rm div}_x \Big(\k_1(\t)\nabla(\f{n_2}{\r})\Big)
 \\
\di
=&-{\rm div}_x \Big(\f{n_2}{\r}\int v \mb{N}_{\mb{M}}^{-1}\big[\mb{P}_c(v\cdot\nabla_x\mb{M}) )\big]\Big)-{\rm div}_x \Big(\int v \mb{N}_{\mb{M}}^{-1}\big[\mb{P}_c(v\cdot\nabla_x (\mb{P}_cF_2))\big]\Big)
\\
&\di
 -{\rm div}_x \Big(\int v \mb{N}_{\mb{M}}^{-1}\big[\nabla_x\Pi\cdot\nabla_v\mb{G}\big]\Big)
 -{\rm div}_x \Big(\int v \mb{N}_{\mb{M}}^{-1}\Big[\partial_t(\mb{P}_c F_2)+(\f{\mb{M}}{\rho})_t~n_2-2Q(F_2,\mb{G})\Big]dv\Big).
\end{aligned}
\end{equation}
Note that the equation \eqref{F2n2-2} is a diffusive equation with the damping term and higher order source terms, which is new and observed first time. The equation \eqref{F2n2-2} is one of the main contributions of the present paper such that we can handle with the electric potential force term in order to prove the stability of basic wave patterns to the bipolar VPB system \eqref{VPB0}.

To make a conclusion, the VPB system \eqref{VPB} can be decomposed to consist of the fluid-dynamical system \eqref{F1-fluid} and non-fluid type equation \eqref{F1-non-fluid} for the distribution function $F_1$, and the fluid-dynamical equation \eqref{F2n2-2} and the non-fluid type equation \eqref{F2-pc} for $F_2$ coupled with the Poisson equation
$$
\Delta \Pi=2n_2.
$$
Note here that the micro-macro decomposition for $F_1$ is similar to the
one for the Boltzmann equation in \cite{Liu-Yu, Liu-Yang-Yu} but with some additional electric potential
force terms, and the new micro-macro decomposition for $F_2$ is made in such a way that it is possible to get the dissipative property of the electric field and further to study the stability toward wave patterns for the bipolar VPB
system \eqref{VPB}. It should be remarked that this decomposition is quite universal and give a unified framework for the stability analysis towards the wave patterns.
As an application of the new decomposition, in the following sections, we prove the nonlinear stability of viscous shock waves and rarefaction waves to the 1D bipolar VPB system as the first step. Furthermore, it can be expected to prove the nonlinear stability of the other elementary wave patterns, such as contact discontinuity, boundary layers and so on,  to the system \eqref{VPB}, which are left in the future works.

\

Now we list some lemmas on the estimates and dissipative properties of collision operators in the weighted $L^2$ space for later use. The following lemmas are based on the celebrated H-theorem. The first lemma is from \cite{Liu-Yu}.

\begin{lemma}\label{Lemma 4.1} There exists a positive
constant $C$ such that
$$
\int\f{\nu(|v|)^{-1}Q(f,g)^2}{\wt{\mb{M}}}dv\leq
C\left\{\int\f{\nu(|v|)f^2}{\wt{\mb{M}}}dv\cdot\int\f{g^2}{\wt{\mb{M}}}dv+
\int\f{f^2}{\wt{\mb{M}}}dv\cdot\int\f{\nu(|v|)g^2}{\wt{\mb{M}}}dv\right\},
$$
where $\wt{\mb{M}}$ can be any Maxwellian so that the above
integrals are well-defined.
\end{lemma}

Based on Lemma \ref{Lemma 4.1}, the following two lemmas are taken
from \cite{Liu-Yang-Yu-Zhao}. Their proofs are straightforward by
using Cauchy inequality.

\begin{lemma}\label{Lemma 4.2} If $\t/2<\t_ *<\t$, then there exist two
positive constants $\wt\sigma=\wt\sigma(\r,u,\t;\break
\r_ *,u_ *,\t_ *)$ and
$\eta_0=\eta_0(\r,u,\t;\r_ *,u_ *,\t_ *)$ such that if
$|\r-\r_ *|+|u-u_ *|+|\t-\t_ *|<\eta_0$, we have for
$g_i( v)\in  \mathfrak{N}_i^\bot~(i=1,2)$,
$$
-\int\f{g_1\mb{L}_\mb{M}g_1}{\mb{M}_ *}dv\geq
\wt\sigma\int\f{\nu(|v|)g_1^2}{\mb{M}_ *}dv,\qquad -\int\f{g_2\mb{N}_\mb{M}g_2}{\mb{M}_ *}dv\geq
\wt\sigma\int\f{\nu(|v|)g_2^2}{\mb{M}_ *}dv.
$$
\end{lemma}

 \begin{lemma}\label{Lemma 4.3} Under the assumptions in Lemma \ref{Lemma 4.2}, we
have  for each $g_i( v)\in  \mathfrak{N}_i^\bot~(i=1,2)$,
$$
 \int\f{\nu(|v|)}{\mb{M}_ *}|\mb{L}_\mb{M}^{-1}g_1|^2dv
\leq \wt\sigma^{-2}\int\f{\nu(|v|)^{-1}g_1^2}{\mb{M}_ *}dv,~~{\rm
and}~~
\int\f{\nu(|v|)}{\mb{M}_ *}|\mb{N}_\mb{M}^{-1}g_2|^2dv\leq
\wt\sigma^{-2}\int\f{\nu(|v|)^{-1}g_2^2}{\mb{M}_ *}dv.
$$
\end{lemma}

\begin{remark}
In Lemmas \ref{Lemma 4.2}-\ref{Lemma 4.3}, $\eta_0$ may not be sufficiently small positive constant. However, in the proof of Theorem \ref{thm} in the following sections, the smallness of $\eta_0$ is crucially used to close the a priori assumptions \eqref{assumption-ap-s} and \eqref{assumption-ap}.
\end{remark}

%%%%%%%
\section{Stability of Boltzmann Shock Profiles for the Bipolar VPB System}
\label{result0}
\setcounter{equation}{0}

In this section, we want to use the micro-macro
decomposition introduced in the previous section to prove the
nonlinear stability of Boltzmann shock profiles to the 1D bipolar
Vlasov-Poisson-Boltzmann system \eqref{VPB0}
\begin{equation}\label{VPB10}
\left\{
\begin{array}{ll}
\di
F_{At} +v_1\partial_x F_A+\partial_x\Pi\partial_{v_1} F_A= Q(F_A,F_A+F_B), \\
\di
F_{Bt} +v_1\partial_x F_B-\partial_x\Pi\partial_{v_1} F_B= Q(F_B,F_A+F_B), \\
\di \partial_{xx} \Pi=\int (F_A-F_B)dv,
\end{array}
\right.
\end{equation}
with $x\in \mb{R}^1, t\in \mb{R}^+, v=(v_1,v_2,v_3)^t\in \mb{R}^3$
and the initial values and the far-field states given by
\begin{equation}\label{VPB10-i}
\begin{array}{ll}
\di F_A(t=0,x,v)=F_{A0}(x,v)\rightarrow \mathbf{M}_{[\rho_\pm,u_\pm,\t_\pm]}(v),\quad {\rm as}~~x\rightarrow\pm\infty,\\
\di F_B(t=0,x,v)=F_{B0}(x,v)\rightarrow \mathbf{M}_{[\rho_\pm,u_\pm,\t_\pm]}(v),\quad {\rm as}~~x\rightarrow\pm\infty,\\
\di \Pi_x\rightarrow 0,\quad {\rm as}~~x\rightarrow\pm\infty,
\end{array}
\end{equation}
where $u_\pm=(u_{1\pm},0,0)^t$ and $\rho_\pm>0, u_{1\pm}, \t_\pm>0$ are
prescribed constant states such that the two states
$(\rho_\pm,u_\pm,\t_\pm)$ are connected by the superposition of  1-shock and 3-shock wave
solutions to the corresponding 1D Euler system
\begin{equation}\label{euler}
\left\{
\begin{array}{l}
\di\rho_t+(\rho u_1)_x=0,\\
\di(\rho u_1)_t+(\rho u_1^2+p)_x=0,\\
\di(\rho u_i)_t+(\rho u_1u_i)_x=0,~i=2,3,\\
\di[\rho(e+\f{|u|^2}{2})]_t+[\rho u_1(e+\f{|u|^2}{2})+pu_1]_x=0.
\end{array}
\right.
\end{equation}
As in the previous section, set
$$
F_1=\f{F_A+F_B}{2},\qquad F_2=\f{F_A-F_B}{2}.
$$
Then one has the following system
\begin{equation}\label{VPB1}
\left\{
\begin{array}{ll}
\di
F_{1t} +v_1\partial_x F_1+\partial_x\Pi\partial_{v_1} F_2= 2Q(F_1,F_1), \\
\di
F_{2t} +v_1\partial_x F_2+\partial_x\Pi\partial_{v_1} F_1= 2Q(F_2,F_1), \\
\di \partial_{xx} \Pi=2\int F_2dv=2n_2,
\end{array}
\right.
\end{equation}
with the initial values $(F_{10}, F_{20})$ satisfying
\begin{equation}
\begin{array}{ll}
\di F_1(t=0,x,v)=F_{10}(x,v)\rightarrow \mathbf{M}_{[\rho_\pm,u_\pm,\t_\pm]}(v),\quad {\rm as}~~x\rightarrow\pm\infty,\\
\di F_2(t=0,x,v)=F_{20}(x,v)\rightarrow 0,\quad {\rm as}~~x\rightarrow\pm\infty,\\
\di \Pi_x\rightarrow 0,\quad {\rm as}~~x\rightarrow\pm\infty.
\end{array}
\end{equation}
Remark that the equation $\eqref{VPB1}_1$ is just the Boltzmann equation coupled with the Vlasov-Poisson electric potential terms.
By the decomposition for $F_1=\mb{M}+\mb{G}$, one can derive from \eqref{F1-fluid} and \eqref{F1-non-fluid}  the following
fluid-part system
\begin{equation}\label{F1-f}
\left\{
\begin{array}{l}
\di \rho_{t}+(\rho u_1)_x=0,\\
\di (\rho u_1)_t+(\rho u_1^2 +p)_x-(\f{\Pi_x ^2}{4})_x=\f{4}{3}(\mu(\t)
u_{1x})_x-\int v_1^2\Gamma_xd v,  \\
\di (\rho u_i)_t+(\rho u_1u_i)_x=(\mu(\t)
u_{ix})_x-\int v_1 v_i\Gamma_xd v,~ i=2,3,\\
\di [\rho(\t+\f{|u|^2}{2})+\f{\Pi_x ^2}{4}]_t+[\rho
u_1(\t+\f{|u|^2}{2})+pu_1]_x=(\k(\t)\t_x)_x+\f{4}{3}(\mu(\t)u_1u_{1x})_x\\
\di\qquad\qquad +\sum_{i=2}^3(\mu(\t)u_iu_{ix})_x
-\int\f12 v_1| v|^2\Gamma_xd v,
\end{array}
\right.
\end{equation}
and the non-fluid equation for $F_1$:
\begin{equation}
\mb{G}_t+\mb{P}_1( v_1\mb{M}_x)+\mb{P}_1(
v_1\mb{G}_x)+\mb{P}_1(\Pi_x \partial_{v_1}F_2)
=\mb{L}_\mb{M}\mb{G}+2Q(\mb{G}, \mb{G}), \label{F1-G}
\end{equation}
where $\mb{G}$ can be expressed explicitly by \eqref{(1.12)} and \eqref{Pi}. Note that the fluid-system \eqref{F1-f} is the compressible Navier-Stokes type system strongly coupled with microscopic terms and electric field terms.
By the new decomposition to $F_2$,
$$
F_2=\frac{n_2}{\rho}\mb{M}+\mb{P}_c F_2,
$$
from the equations \eqref{m22} and \eqref{F2-pc}, the macroscopic part $n_2$ satisfies the following equation
\begin{equation}\label{n21}
n_{2t}+(\int v_1F_2dv)_x=0,
\end{equation}
or equivalently,
\begin{equation}\label{n2}
\begin{array}{ll}
\di n_{2t}+(u_1n_2)_x+ \Big(\f{\k_1(\t)}{R\t}\Pi_x \Big)_x- \Big(\k_1(\t)(\f{n_2}{\r})_x\Big)_x\\
\di
 =-\Big(\f{n_2}{\r}\int v_1 \mb{N}_{\mb{M}}^{-1}\big[\mb{P}_c(v_1\mb{M}_x) )\big]dv\Big)_x-\Big(\int v_1 \mb{N}_{\mb{M}}^{-1}\big[\mb{P}_c(v_1 (\mb{P}_cF_2)_x)\big]dv\Big)_x\\
\di
 - \Big(\int v_1 \mb{N}_{\mb{M}}^{-1}\big[\Pi_x \mb{G}_{v_1}\big]dv\Big)_x -\Big(\int v_1 \mb{N}_{\mb{M}}^{-1}\Big[\partial_t(\mb{P}_c F_2)+(\f{\mb{M}}{\rho})_t~n_2
 -2Q(F_2,\mb{G})\Big]dv\Big)_x,
\end{array}
\end{equation}
and the microscopic part $\mb{P}_d F_2$ satisfies the equation
\begin{equation}\label{F2-pc-1}
\partial_t(\mb{P}_c F_2)-\mb{N}_{\mb{M}}(\mb{P}_c F_2)+\mb{P}_c(v_1 F_{2x})+\mb{P}_c (\Pi_x \partial_{v_1}F_1)+(\f{\mb{M}}{\rho})_t~n_2=2Q(F_2,\mb{G}).
\end{equation}
Then the stability analysis towards the Boltzmann shock profiles for the bipolar VPB system \eqref{VPB} will be carried out for the equivalent system \eqref{F1-f}-\eqref{F2-pc-1}.

\subsection{Existence of Boltzmann Shock Profiles}
In this subsection, we state the existence and list the properties of Boltzmann shock profile. Consider the viscous shock profile to the 1D slab symmetric Boltzmann equation
\begin{equation}\label{B}
F_{1t}+v_1F_{1x}=2Q(F_1,F_1).
\end{equation}
We recall the Riemann problem for
the compressible Euler equation \eqref{euler} with the Riemann initial data
\begin{equation}
\begin{array}{ll}
(\r,m,E)(x,0)=\left\{\begin{array}{ll}
(\r_{-},m_-,E_{-}),~x<0,\\
(\r_{+},m_+,E_{+}),~x>0,
\end{array}\right.
\end{array}\label{Rdata}
\end{equation}
where $m_{\pm}=\r_{\pm}u_{\pm}$, $E_{\pm}=\r_{\pm}(\t_{\pm}+\frac{|u_{\pm}|^2}{2})$ and $u_\pm=(u_{1pm}, 0, 0)^t$.
It is well-known that the system \eqref{euler} has three eigenvalues:
$\lambda_1=u_1-\f{\sqrt{10\t}}{3}, \lambda_2=u_1, \lambda_3=u_1+\f{\sqrt{10\t}}{3}$
where the second characteristic field is linear degenerate and the
other two are genuinely nonlinear.
In the present section, we focus our attention on the situation where
the Riemann solution of \eqref{euler}, \eqref{Rdata} consists of two shock waves
(and three constant states),
that is, there exists an intermediate
state $(\r_\#,m_\#=\r_\#u_\#,E_\#=\r_\#(\t_\#+\frac{|u_\#|^2}{2}))$ with $u_\#=(u_{1\#}, 0, 0)^t$ such that $(\r_{-},m_{-},E_{-})$
connects with $(\r_\#,m_\#,E_\#)$ by the 1-shock wave with the shock speed
$s_1$ and $(\r_\#,m_\#,E_\#)$ connects with $(\r_+,m_+,E_+)$
by the 3-shock wave with the shock speed $s_3$. Here
the shock speeds $s_1$ and $s_3$ are
constants determined by R-H conditions \eqref{RH-E} and \eqref{RH-E1} and the entropy conditions \eqref{Lax-E}.
By the standard arguments (e.g. \cite{Smoller}), for each
$(\r_{-},m_{-},E_{-})$, we can see our
situation takes place provided $(\r_+,m_+,E_+)$
is located on a curved surface in a neighborhood of $(\r_{-},m_{-},E_{-})$.
In what follows, the neighborhood of $(\r_{-},m_{-},E_{-})$ is denoted by
$\Omega_-$.
To describe the
strengths of the shock waves for later use, we set
$$
\delta^{S_1}=|\r_\#-\r_-|+|m_\#-m_-|+|E_\#-E_-|, \quad
\delta^{S_3}=|\r_\#-\r_+|+|m_\#-m_+|+|E_\#-E_+|
$$
and $\delta=\min\{\delta^{S_1},\delta^{S_3}\}$.
When we choose
$|(\r_+-\r_-,m_+-m_-,E_+-E_-)|$ small in our situation
for the fixed $(\r_{-},m_{-},E_{-})$,
we note that it holds that
\begin{equation*}
\delta^{S_1}+\delta^{S_3} \le C|(\r_+ -\r_-,m_+ -m_-,E_+ -E_-)|,
\end{equation*}
where $C$ is a positive constant depending only on $(\r_-,m_-,E_-)$.
Then,
if it also holds
\begin{equation}
\delta^{S_1}+\delta^{S_3} \le C\delta,\quad {\rm as}~~\delta^{S_1}+\delta^{S_3} \to 0
\label{delta}
\end{equation}
for a positive constant $C$,
we call the strengths of the shock waves are ``small with same order".
In what follows, we always assume \eqref{delta}.

Then we recall the $i$-shock profile $ F_1^{S_i}(x- s_it, v)~(i=1,3)$ of the
Boltzmann equation \eqref{B} in Eulerian coordinates with its existence and
properties given in the papers by Caflisch-Nicolaenko \cite{Caflish-Nicolaenko} and Liu-Yu
\cite{Liu-Yu, Liu-Yu-1}. Note that the compressibility of the Boltzmann shock profile is first proved
in \cite{Liu-Yu-1} , which is crucial for the stability analysis towards the Boltzmann shock profile. First of all, the $i$-shock profile $ F_1^{S_i}(x- s_it, v)$ satisfies
\begin{equation}
\left\{
\begin{array}{ll}
\di -s_i(F_1^{S_i})^\prime+v_1(F_1^{S_i})^\prime=2Q(F_1^{S_i},F_1^{S_i}),~~~i=1,3,\\[5mm]
F_1^{S_1}(-\i, v)=\mathbf{M}_{[\r_-,u_-,\t_-]}( v),\qquad F_1^{S_3}(-\i, v)=\mathbf{M}_{[\r_\#,u_\#,\t_\#]}( v),\\[5mm]
F_1^{S_1}(+\i, v)=\mathbf{M}_{[\r_\#,u_\#,\t_\#]}( v),\qquad F_1^{S_3}(+\i, v)=\mathbf{M}_{[\r_+,u_+,\t_+]}( v),\\[5mm]
\end{array}
\right.\label{BS}
\end{equation}
where $^\prime=\f{d}{d\vartheta_i}$, $\vartheta_i=x- s_it,$ $u_\pm=(u_{1\pm},0,0)^t$, $u_\#=(u_{1\#},0,0)^t$
and
$(\r_{\pm},u_{\pm},\t_\pm)$, $(\r_\#,u_\#,\t_\#)$ satisfy Rankine-Hugoniot conditions
\begin{equation}
\left\{
\begin{array}{ll}
\di- s_1(\r_\#-\r_-)+(\r_\#u_{1\#}-\r_-u_{1-})=0,\\
\di -
s_1(\r_\#u_{1\#}-\r_-u_{1-})+(\r_\#u_{1\#}^2+p_\#-\r_-u_{1-}^2-p_-)=0,\\
\di -
s_1(\r_\#E_\#-\r_-E_-)+(\r_\#u_{1\#}E_\#+p_\#u_{1\#}-\r_-u_{1-}E_--p_-u_{1-})=0,
\end{array}
\right. \label{RH-E1}
\end{equation}
\begin{equation}
\left\{
\begin{array}{ll}
\di- s_3(\r_+-\r_\#)+(\r_+u_{1+}-\r_\#u_{1\#})=0,\\
\di -
s_3(\r_+u_{1+}-\r_\#u_{1\#})+(\r_+u_{1+}^2+p_+-\r_\#u_{1\#}^2-p_\#)=0,\\
\di -
s_3(\r_+E_+-\r_\#E_\#)+(\r_+u_{1+}E_++p_+u_{1+}-\r_\#u_{1\#}E_\#-p_\#u_{1\#})=0,
\end{array}
\right. \label{RH-E}
\end{equation}
and Lax entropy conditions
\begin{equation}
 \l_{1\#}<s_1<\l_{1-}, \qquad \l_{3+}< s_3<\l_{3\#},
 \label{Lax-E}
\end{equation}
with $s_i$ being $i$-shock wave speed and
$\l_i=u_1+(-1)^{\f{i+1}{2}}\f{\sqrt{10\t}}{3}~(i=1,3)$ being the $i$-th characteristic
eigenvalue of the Euler equations in the Eulerian coordinate and
$\l_{1-}=u_{1-}-\f{\sqrt{10\t_-}}{3}$, $\l_{i\#}=u_{1\#}+(-1)^{\f{i+1}{2}}\f{\sqrt{10\t_\#}}{3}~(i=1,3)$ and $\l_{3+}=u_{1+}+\f{\sqrt{10\t_+}}{3}$.
 By the micro-macro decomposition to the Boltzmann equation \eqref{B}  around the local Maxwellian
$\mb{M}^{S_i}~(i=1,3)$ (cf. \cite{Liu-Yu, Liu-Yang-Yu}), it holds that
$$
F_1^{S_i}(x-s_it,v)=\mathbf{M}^{S_i}(x-s_it,v)+\mathbf{G}^{S_i}(x-s_it,v),
$$
where
$$
\mathbf{M}^{S_i}(x-s_it,v)=\mathbf{M}_{[\r^{S_i},u^{S_i},\t^{S_i}]}(x-s_it,v)=\f{\r^{S_i}(x-s_it)}{\sqrt{(2\pi
R\t^{S_i}(x-s_it))^3}}e^{-\f{|v-u^{S_i}(x-s_it)|^2}{2R\t^{S_i}(x-s_it)}},
$$
with
\begin{equation*}
\left(
\begin{array}{c}
 \r^{S_i}\\
 \rho^{S_i} u^{S_i}_j\\
 \rho^{S_i}(\t^{S_i}+\frac{|u^{S_i}|^2}{2})
\end{array}
\right)=\int_{\mathbf{R}^3}\left(
\begin{array}{c}
 1\\
 v_j\\
\frac{|v|^2}{2}
\end{array}
\right)F_1^{S_i}(x-s_it, v)d v,~~j=1,2,3.\label{shock-fluid}
\end{equation*}
With respect to the inner product
$\langle\cdot,\cdot\rangle_{\mathbf{M}^{S_i}}$ defined in
\eqref{product},  we can now define the macroscopic projection
$\mb{P}^{S_i}_0$ and microscopic projection $\mb{P}^{S_i}_1$ by
\begin{equation*}
 \mb{P}^{S_i}_0g = {\di\sum_{j=0}^4\langle g,\chi^{S_i}_j\rangle_{\mathbf{M}^{S_i}}\chi^{S_i}_j,}
 \qquad
 \mb{P}^{S_i}_1g= g-\mb{P}^{S_i}_0g,
\label{PS3}
\end{equation*}
where $\chi^{S_i}_j~(j=0,1,2,3,4)$ are the corresponding pairwise
orthogonal base defined in \eqref{orthogonal-base} by replacing
$(\r,u,\t, \mb{M})$ by $(\r^{S_i},u^{S_i},\t^{S_i}, \mb{M}^{S_i})$.
Under the above micro-macro decomposition, the solution
$F_1^{S_i}=F_1^{S_i}(x- s_it, v)$ satisfies
$$
\mb{P}^{S_i}_0F_1^{S_i}=\mb{M}^{S_i},~\mb{P}^{S_i}_1F^{S_i}=\mb{G}^{S_i},
$$
and the Boltzmann equation $(\ref{B})$ becomes
\begin{equation*}
(\mb{M}^{S_i}+\mb{G}^{S_i})_t+ v_1(\mb{M}^{S_i}+\mb{G}^{S_i})_x
=2\big[Q(\mb{M}^{S_i},\mb{G}^{S_i})+Q(\mb{G}^{S_i},\mb{M}^{S_i})\big]+2Q(\mb{G}^{S_i},\mb{G}^{S_i}).
\label{(1.6)}
\end{equation*}
Correspondingly, we have the following fluid-type system for the
fluid components of shock profile:
\begin{equation}
\begin{array}{l}
\left\{
\begin{array}{l}
\di \r^{S_i}_{t}+(\r^{S_i} u^{S_i}_1)_x=0,\\
\di (\r^{S_i} u^{S_i}_1)_t+[\r^{S_i} (u^{S_i}_1)^2
+p^{S_i}]_x=\f{4}{3}\big(\mu(\t^{S_i})
u^{S_i}_{1x}\big)_x-\int v_1^2\Gamma^{S_i}_xd v,  \\
\di (\r^{S_i} u^{S_i}_j)_t+(\r^{S_i}
u^{S_i}_1u^{S_i}_j)_x=\big(\mu(\t^{S_i})
u^{S_i}_{jx}\big)_x-\int v_1 v_j\Gamma^{S_i}_xd v,~ j=2,3,\\
\di [\rho^{S_i}(\t^{S_i}+\f{|u^{S_i}|^2}{2})]_t+[\r^{S_i}
u^{S_i}_1(\t^{S_i}+\f{|u^{S_i}|^2}{2})+p^{S_i}u^{S_i}_1]_x=\big(\k(\t^{S_i})\t^{S_i}_x\big)_x\\
\di\qquad
+\f{4}{3}\big(\mu(\t^{S_i})u^{S_i}_1u^{S_i}_{1x}\big)_x+\sum_{j=2}^3\big(\mu(\t^{S_i})u^{S_i}_ju^{S_i}_{jx}\big)_x
-\int\f12 v_1| v|^2\Gamma^{S_i}_xd v.
\end{array}
\right.
\end{array}\label{shock-Euler}
\end{equation}
In fact, from the invariance of the equation \eqref{BS} by changing
$v_j$ with $-v_j$ for $j=2,3$ and the fact that $u_{j\pm}=0,$ we have $\di
u^{S_i}_j=\int v_1 v_j\Gamma^{S_i}_xd v\equiv0$ for $j=2,3.$
And the  equation for the non-fluid component
${\mb{G}}^{S_i}~(i=1,3)$ is
\begin{equation}\label{shock-non-fluid}
\mb{G}^{S_i}_t+\mb{P}^{S_i}_1( v_1\mb{M}^{S_i}_x)+\mb{P}^{S_i}_1( v_1\mb{G}^{S_i}_x)
=\mb{L}_{\mb{M}^{S_i}}\mb{G}^{S_i}+2Q(\mb{G}^{S_i}, \mb{G}^{S_i}).
\end{equation}
Here $\mb{L}_{\mb{M}^{S_i}}$ is the linearized collision operator of
$Q( F_1^{S_i}, F_1^{S_i})$
 with respect to the local Maxwellian $\mb{M}^{S_i}$:
$$
\mb{L}_{\mb{M}^{S_i}} g=2[Q(\mb{M}^{S_i}, g)+
Q(g,\mb{M}^{S_i})].
$$
Thus
\begin{equation}
\begin{array}{l}
\di \mb{G}^{S_i}= \mb{L}_{\mb{M}^{S_i}}^{-1}[\mb{P}^{S_i}_1( v_1\mb{M}^{S_i}_x)] +\Gamma^{S_i}_x ,\qquad\quad \Gamma^{S_i}_x =\mb{L}_{\mb{M}^{S_i}}^{-1}[(\mb{G}^{S_i}_t+\mb
{P}^{S_i}_1( v_1\mb{G}^{S_i}_x))-2Q(\mb{G}^{S_i}, \mb{G}^{S_i})].\label{PiS3}
\end{array}
\end{equation}
Now we recall  the properties of the shock profile $ F_1^{S_i}(x-
s_it,v)~(i=1,3)$ that are given or can be induced by  Liu-Yu  \cite{Liu-Yu-1} in Theorem
6.8.

\begin{lemma}\label{Lemma-shock} (\cite{Liu-Yu-1})  For $i=1,3$, if the shock wave strength $\d^{S_i}$ is small enough,
then the Boltzmann equation \eqref{B} admits a $i$-shock profile
solution $ F_1^{S_i}(x- s_it,v)$ uniquely up to a shift satisfying
the following properties:
\begin{itemize}
\item [(1)] The shock profile $ F_1^{S_i}(x- s_it,v)$ converges to its far fields
exponentially fast with an exponent proportional to the magnitude of
the shock wave strength, that is
\begin{equation*}
\left\{
\begin{array}{ll}
\di |(\r^{S_1}-\r_{-},u^{S_1}_1-u_{1-},\t^{S_1}-\t_-)|\leq
C\d^{S_1} e^{-c\d^{S_1}|\vartheta_1|},&\di {\rm
as}~~\vartheta_1<0,\\
\di |(\r^{S_1}-\r_\#,u^{S_1}_1-u_{1\#},\t^{S_1}-\t_\#)|\leq
C\d^{S_1} e^{-c\d^{S_1}|\vartheta_1|},&\di {\rm
as}~~\vartheta_1>0,\\
\di |(\r^{S_3}-\r_{+},u^{S_3}_1-u_{1+},\t^{S_3}-\t_+)|\leq
C\d^{S_3} e^{-c\d^{S_3}|\vartheta_3|},&\di{\rm
as}~~\vartheta_3>0,\\
\di |(\r^{S_3}-\r_\#,u^{S_3}_1-u_{1\#},\t^{S_3}-\t_\#)|\leq
C\d^{S_3} e^{-c\d^{S_3}|\vartheta_3|},&\di {\rm
as}~~\vartheta_3<0,\\
\di\big(\int
\f{\nu(|v|)|\mb{G}^{S_i}|^2}{\mb{M}_0}dv\big)^{\f12}\leq
C(\d^{S_i})^2 e^{-c \d^{S_i}|\vartheta_i|},~~i=1,3,
\end{array}
\right.
\end{equation*}
with  $\d^{S_i}$ being the $i$-shock strength and  $\mb{M}_0$ being
the global Maxwellian which is close to the shock profile
with its precise definition given in Theorem 6.8, \cite{Liu-Yu-1}.

\item[(2)] Compressibility of $i$-shock profile:
$$
(\l^{S_i}_i)_{\vartheta_i}< 0,\qquad
\l^{S_i}_i=u^{S_i}_1+(-1)^{\frac{i+1}{2}}\f{\sqrt{10\t^{S_i}}}{3}.
$$

\item[(3)] The following properties hold:
$$
\begin{array}{ll}
\di \r^{S_i}_{\vartheta_i}\sim u^{S_i}_{1\vartheta_i}\sim
\t^{S_i}_{\vartheta_i}\sim(\l^{S_i}_i)_{\vartheta_i}\sim \big(\int
\f{\nu(|v|)|\mb{G}^{S_i}|^2}{\mb{M}_0}dv\big)^{\f12},
\end{array}
$$
where $A\sim B$ denotes the equivalence of the quantities $A$ and
$B$, and
$$
\left\{
\begin{array}{ll}
\di u^{S_i}_j\equiv 0,~~\int v_1v_j\Gamma^{S_i}_{\vartheta_i} dv\equiv0, ~~j=2,3,\\
\di \di |\partial^k_{\vartheta_i}(\r^{S_i},u^{S_i}_{1},\t^{S_i})|\leq
C(\d^{S_i})^{k-1}|(\r^{S_i}_{\vartheta_i},
u^{S_i}_{1\vartheta_i},
\t^{S_i}_{\vartheta_i})|,~k\ge 2,\\[3mm]
\di \big(\int
\f{\nu(|v|)|\partial^k_{\vartheta_i}\mb{G}^{S_i}|^2}{\mb{M}_0}dv\big)^{\f12}\leq
C(\d^{S_i})^k\big(\int
\f{\nu(|v|)|{\mb{G}}^{S_i}|^2}{\mb{M}_0}dv\big)^{\f12},~~k\ge 1,\\
 \di
|\int v_1\xi_j(v)\Gamma^{S_i}_{\vartheta_i}dv|\leq
C\d^{S_i}|u^{S_i}_{1\vartheta_i}|,~~j=1,2,3,4,
\end{array}
\right.
$$
where $\xi_j(v)~(j=1,2,3,4)$ are the collision invariants
defined in \eqref{collision-invar}.
\end{itemize}
\end{lemma}

Let $ U(x,t)=(\r,m,E+\f{\Pi_x ^2}{4})^t$ with $m=\r u$ and $ E=\r(\t+\f{u^2}{2})$ being the solution of the system \eqref{F1-f} and $
\bar{U}(x,t)=(\bar{\r},\bar{m},\bar{E})^t$ being the linear superposition of 1-shock and 3-shock profile with
\begin{equation}
\left\{\begin{array}{ll}
{\bar \r}=\r^{S_1}(x-s_1t)+\r^{S_3}(x-s_3t)-\r_\#,\quad {\bar E}=E^{S_1}(x-s_1t)+E^{S_3}(x-s_3t)-E_\#,\\[3pt]
{\bar m_1}=m_1^{S_1}(x-s_1t)+m_1^{S_3}(x-s_3t)-m_{1\#}, ~~
{\bar m_i}=0,   ~~~~i=2,3 ,
\end{array}\right.
\end{equation}
where $(\cdot)^t$ denotes the transpose of the vector $(\cdot)$. Note that in order to keep the conservative form of the system \eqref{F1-f}, an additional term $\f{\Pi_x ^2}{4}$, which means the electric potential energy should be put on the total energy $E+\f{\Pi_x ^2}{4}$, which is quite different from the classic Boltzmann equation.

Since the present paper is concerning with general initial
perturbation, that is, the integral $\di \int_{-\infty}^\infty
(U-\bar{U})(x,0)\,dx$ may not be zero, to use the
anti-derivative technique, we need to find an ansatz $\tilde{U}$
such that $\di \int_{-\infty}^\infty (U(x,t)-\tilde{U}(x,t))\,dx=0$,
meanwhile, $\tilde U$ and $\bar U $ are time-asymptotically equivalent, i.e., $|\tilde{U}-\bar{U}|$ tends to zero as $t\to \infty$. As observed by Liu \cite{Liu-1985}, the general perturbations will not only produce the shift on the viscous shock wave itself, but also the linear or nonlinear diffusion waves in the transverse fields.  Let
$$
A(\r,m_1,E)=\left(\begin{array}{ccc}
0&1&0\\
-\f{m^2_1}{\r^2}+\f{m^2}{3\r^2}&\f{4m_1}{3\r}&\f23\\
-\f{5m_1E}{3\r^2}+\f{2m^2m_1}{3\r^3}&\frac{5E}{3\r}-\f{2m_1^2}{3\r^2}-\f{m^2}{3\r^2}
&\f{5m_1}{3\r}\end{array}\right)
$$
be the Jacobi matrix of the flux
$(m_1,\f{2}{3}E+\f{m_1^2}{\r}-\f{m^2}{3\r},\f{5m_1E}{3\r}-\f{m_1m^2}{3\r^2})^t$ in the Euler system \eqref{euler}
with respect to $(\r,m_1,E)$. Then the second right eigenvector of
the matrix $A(\r,m_1,E)$ at the intermediate state $(\r_\#,m_{1\#},E_\#)$ is
$$
r_2=(1,u_{1\#},\frac{u_{1\#}^2}{2})^t.
$$
Furthermore, a direct computation yields
that the three vectors $r_1=(\r_\#-\r_-,m_\#-m_-,E_\#-E_-)^t$, $r_2$ above and
$r_3=(\r_+-\r_\#,m_+-m_\#,E_+-E_\#)^t$ are linearly independent in $\mb{R}^3$
if $\delta^{S_1}+\delta^{S_3}$ is suitably small. So if the initial mass $\di \int
(U(x,0)-\bar{U}(x,0))\,dx$ is not zero, we can distribute the initial
mass along the three independent directions $r_1,r_2$ and $r_3$ as in Liu \cite{Liu-1985},
that is,
\begin{equation}
\int (U(x,0)-\bar{U}(x,0))\,dx =\sum_{i=1}^3\alpha_ir_i,
\label{initial}
\end{equation}
where $\alpha_i~(i=1,2,3)$ are constants uniquely determined by the
initial data. The excessive mass $\a_1r_1$ in the first
characteristic field can be eliminated by the translated 1-viscous
shock wave with a shift $\a_1$, i.e., $\r^{S_1}(x-s_1t+\a_1)$. Similarly,
we can eliminate $\a_3r_3$ by replacing $\r^{S_3}(x-s_3t)$ by
$\r^{S_3}(x-s_3t+\a_3)$. So the remaining problem is how to remove the
excessive mass in the second characteristic field, i.e., $\a_2r_2$.
Motivated by Huang-Matsumura \cite{Huang-matsumura},
 the desired ansatz  $\tilde{U}=(\tilde{\r},\tilde{m},\tilde{E})^t$ is constructed
 in the following form:
\begin{equation}
\begin{array}{ll}\label{ansatz}
\tilde{\r}=\r^{S_1}(x-s_1t+\alpha_1)+\r^{S_3}(x-s_3t+\alpha_3)-\r_\#+\Theta(x,t), \quad \tilde{m}_i= 0 ,~i=2,3,\\
\tilde{m}_1=m_1^{S_1}(x-s_1t+\alpha_1)+m_1^{S_3}(x-s_3t+\alpha_3)-m_{1\#}+u_{1\#}\Theta(x,t)-a\Theta_x(x,t),
\\
\di\tilde{E}=E^{S_1}(x-s_1t+\alpha_1)+E^{S_3}(x-s_3t+\alpha_3)-E_\#+\f{1}{2}u_{1\#}^2\Theta(x,t)-au_{1\#}\Theta_x(x,t) ,
\end{array}
\end{equation}
where $\alpha_i~(i=1,3)$ are the shifts of the $i-$viscous shock wave and $\Theta$ is the linear diffusion wave
\begin{equation}
\Theta(x,t)=\frac{\a_2}{\sqrt{4\pi a(1+t)}}e^{-\frac{(x-u_{1\#}t)^2}{4a(1+t)}},
\quad a=\frac{3\k(\t_\#)}{5\r_\#}>0
\label{Theta}
\end{equation}
satisfying the heat equation
\begin{equation}
\Theta_t+u_{1\#}\Theta_x=a\Theta_{xx}, \quad \Theta|_{t=-1}=\alpha_2\delta(x),\quad
\int_{-\infty}^\infty \Theta(x,t)\,dx=\alpha_2,
\label{heat}
\end{equation}
and the terms $-a\Theta_x$ and $-au_{1\#}\Theta_x$ are the coupled diffusion waves as first introduced by Szepessy-Xin \cite{Szepessy-Xin} which does not carry the initial mass, but get rid of some bad error terms decaying not enough with respect to the time $t$. Then it follows from \eqref{initial} that
\begin{equation}
\begin{array}{ll}
\di\int_{-\infty}^\infty
(U(x,0)-\tilde{U}(x,0))\,dx
=\int_{-\infty}^\infty
(U(x,0)-\bar{U}(x,0))\,dx+\int_{-\infty}^\infty
(\bar{U}(x,0)-\tilde{U}(x,0))\,dx=0,
\end{array}
\end{equation}
where we have used the fact that $\int_{-\i}^\i \Theta \,dx=\a_2$. Thus
$\tilde{U}(x,t)$ is the desired ansatz.

By Lemma \ref{Lemma-shock} ,  for suitably small $\delta$ and $\alpha_2$, we have
the wave interaction estimates between two viscous shock waves
\begin{equation}
|\r^{S_1}-\r_\#|\cdot|\r^{S_3}-\r_\#| \le C \delta^{S_1}\delta^{S_3}
(e^{-c\delta^{S_1}(|x|+t)+c\delta^{S_1}|\alpha_1|}+e^{-c\delta^{S_3}(|x|+t)+c\delta^{S_3}|\alpha_3|}) \leq C\delta^2e^{-c\delta(|x|+t)},
\label{interaction1}
\end{equation}
and between the $i-$viscous shock waves ($i=1, 3$) and the diffusion waves
\begin{equation}
\begin{array}{rll}
|\r^{S_i}-\r_\#|\cdot |\Theta| &\le
\di{C|\alpha_2|\delta^{\frac32}e^{-c\delta(|x|+t)}
+ C\frac{|\alpha_2|}{(1+t)^{\frac32}}e^{-\frac{c(x-u_{1\#}t)^2}{1+t}}
+C(\delta+|\a_2|)e^{-c(|x|+t)},}
\end{array}\label{interaction2}
\end{equation}
where we used the fact that $\delta^{S_1}\alpha_1$ and
$\delta^{S_3}\alpha_3$ are uniformly bounded by \eqref{initial} for small
$\delta^{S_1}$ and $\delta^{S_3}$ as long as the initial perturbation stays
bounded. Set
\begin{equation}\label{Q}
Q=\Big\{q(t,x): |q|\le
C(\delta^2+|\a_2|\delta^{\frac32})e^{-c\delta(|x|+t)}+
C\frac{|\alpha_2|}{(1+t)^{\frac32}}e^{-\frac{c(x-u_{1\#}t)^2}{1+t}}+C(\delta+|\a_2|)e^{-c(|x|+t)}
\Big\}.
\end{equation}
From now on, we denote  $q$ as a function belonging to the set $Q$ if without confusions. Applying \eqref{ansatz}, we can calculate directly that
\begin{equation*}
\tilde u_1
\di\sim u^{S_1}_1+u^{S_3}_1-u_{1\#}-\frac{a}{\r_\#}\Theta_x,~~~~{\rm and}~~~~
\frac{\tilde m_1^2}{\tilde\r}\sim\frac{(m_1^{S_1})^2}{\r^{S_1}}+\frac{(m_1^{S_3})^2}{\r^{S_3}}-\frac{m_{1\#}^2}{\r_\#}+u^2_{1\#}\Theta-2u_{1\#}a\Theta_x,
\end{equation*}
where $A\sim B$ means that $A=B+q$ with $q\in Q.$
Similarly, we can calculate $\tilde p$, $\tilde\t$, $\tilde u_{1x}$, $\mu(\tilde\t)\tilde u_{1x}$, $\frac{\tilde m_1\tilde E}{\tilde\r}$,
$\frac{\tilde m_1\tilde p}{\tilde\r}$, $\tilde\t_x$, $\kappa(\tilde\t)\tilde\t_x$ and $\mu(\tilde\t)\tilde u_1\tilde u_{1x}$.
Henceforth, we can check that the ansatz $\tilde{U}=(\tilde{\r},\tilde{m},\tilde{E})^t$ defined in \eqref{ansatz} well
approximates solution of the VPB system as
\begin{equation}
\left\{\begin{array}{ll}
&\di\tilde{\r}_{t}+\tilde{m}_{1x}=0,\\
&\di\tilde{m}_{1t}+(\f{\tilde{m}_1^2}{\tilde{\r}}+\tilde{p})_{x}=
\f{4}{3}(\mu(\tilde{\t})\tilde{u}_{1x})_x-\int v_1^2(\G^{S_1}_x+\G^{S_3}_x)d v+Q_{1x},\\
&\di\tilde{E}_{t}+
(\f{\tilde{m}_1\tilde{E}}{\tilde{\r}}+\f{\tilde{p}\tilde{m}_1}{\tilde{\r}})_{x}=
(\k(\tilde{\t})\tilde{\t}_x)_x+\f{4}{3}(\mu(\tilde{\t})\tilde{u}_1\tilde{u}_{1x})_x-
\int v_1\f{| v|^2}{2}(\G^{S_1}_x+\G^{S_3}_x)d v+Q_{2x},
\end{array}\right.\label{jinsi}
\end{equation}
where
\begin{equation}
\begin{array}{ll}
\di Q_1&\di=\Big(\frac{\tilde m_1^2}{\tilde\r}-\frac{(m_1^{S_1})^2}{\r^{S_1}}-\frac{(m_1^{S_3})^2}{\r^{S_3}}+\frac{m_{1\#}^2}{\r_\#}\Big)
+\Big(\tilde p-p^{S_1}-p^{S_3}+p_\#\Big)\\[3mm]
&\di\quad-\frac{4}{3}\Big(\mu(\tilde\t)\tilde u_{1x}-\mu(\t^{S_1})u^{S_1}_{1x}-\mu(\t^{S_3})u^{S_3}_{1x}\Big)
+2u_{1\#}a\Theta_x-u_{1\#}^2\Theta-a^2\Theta_{xx},
\end{array}\label{Q1}
\end{equation}
and
\begin{equation}
\begin{array}{ll}
\di Q_2&\di=\Big(\frac{\tilde m_1\tilde E}{\tilde\r}-\frac{m_1^{S_1}E^{S_1}}{\r^{S_1}}-\frac{m_1^{S_3}E^{S_3}}{\r^{S_3}}+\frac{m_{1\#}E_\#}{\r_\#}\Big)
+\Big(\frac{\tilde m_1\tilde p}{\tilde\r}-\frac{m_1^{S_1}p^{S_1}}{\r^{S_1}}-\frac{m_1^{S_3}p^{S_3}}{\r^{S_3}}+\frac{m_{1\#}p_\#}{\r_\#}\Big)\\[3mm]
&\di\quad-\Big(\kappa(\tilde\t)\tilde\t_x-\kappa(\t^{S_1})\t^{S_1}_x-\kappa(\t^{S_3})\t^{S_3}_x\Big)
-\frac{4}{3}\Big(\mu(\tilde\t)\tilde u_1\tilde u_{1x}-\mu(\t^{S_1})u^{S_1}u^{S_1}_{1x}-\mu(\t^{S_3})u^{S_3}u^{S_3}_{1x}\Big)\\
&\di\quad+\frac{3}{2}u_{1\#}^2a\Theta_x-\frac{1}{2}u_{1\#}^3\Theta-a^2u_{1\#}\Theta_{xx}.
\end{array}\label{Q2}
\end{equation}
Obviously, it holds that $Q_1, Q_2\sim 0$, i.e.,  $Q_1$, $Q_2$ $\in Q$.
Under the above preparation, we are now at the stage to state our main result.
We first fix any $(\r_{-},m_{-},E_{-})$, and assume that $(\r_+,m_+,E_+) \in \Omega_-$
and the Riemann solution of \eqref{euler}, \eqref{Rdata} consists of two shock waves.
Then the macrpscopic composite wave $(\tilde{\r},\tilde{m}=\tilde\r\tilde u,\tilde{E}=\tilde\r(\tilde\t+\frac{|\tilde u|^2}{2}))(x,t)$
defined in \eqref{ansatz} is well defined. Denote the perturbation around the ansatz by
\begin{equation}
\begin{array}{ll}
\di(\Phi_x,\Psi_x,W_x)(t,x)=(\phi,\psi,\omega)(t,x)=(\r-\tilde{\r},m-\tilde{m},E+\f{\Pi_x ^2}{4}-\tilde{E})(t,x),\\[2mm]
\di \wt
{\mathbf{G}}(t,x,v)=\mathbf{G}(t,x,v)-\mathbf{G}^{S_1}(x-s_1t+\a_1,v)-\mathbf{G}^{S_3}(x-s_3t+\a_3,v),\\[2mm]
\di \wt F_1(t,x,v)=F_1(t,x,v)-F^S_{\a_1,\a_3}(t,x,v),
\end{array}\label{perturb-s}
\end{equation}
with
\begin{equation}\label{su-s}
 F^S_{\a_1,\a_3}(t,x,v)=F_1^{S_1}(x-s_1t+\a_1,v)+F_1^{S_3}(x-s_3t+\a_3,v)-\mb{M}_\#
\end{equation}
 being the superposition of shifted 1-shock profile and 3-shock profile to the 1D Boltzmann equation \eqref{B} and $\mb{M}_\#= \f{\rho_\#}{\sqrt{ (2 \pi
R \t_\#)^3}} e^{-\f{| v-u_{\#}|^2}{2R\t_\#}}$ is the intermediate equillbrium state.

\subsection{Main Result}
Denote
\begin{equation}\label{et}
\begin{array}{ll}
\di \mathcal{E}(t)= \sup_{s\in[0,t]}\Big\{\|(\Phi,\Psi,W)(\cdot,s)\|^2_{H^2_x}
+\sum_{0\leq|\beta|\leq2}\int\int\f{|\partial^\beta(\widetilde{\mb{G}},\mb{P}_cF_{2})|^2}{\mb{M}_*}dv
dx+\|(\Pi_{x},n_{2},n_{2x})(\cdot,s)\|^2 \\
\di+\sum_{|\alpha|=1,0\leq|\beta|\leq1}\int\int\f{|\partial^\alpha\partial^\beta(\widetilde{\mb{G}},\mb{P}_cF_{2})|^2}{\mb{M}_*}dv
dx+\sum_{|\alpha|=2}\int\int\f{|\partial^\alpha(\widetilde{F}_{1},F_{2})|^2}{\mb{M}_*}dv dx\Big\}.\\
\end{array}
\end{equation}
with $\mb{M}_*$ being the global Maxellian chosen in Theorem \ref{thm0} and $\partial^\alpha=\partial_{x,t}^\alpha$, $\partial^\beta=\partial_{v}^\beta$.
Now we can state our main result as follows.

\begin{theorem}\label{thm0} There exist positive constants $\delta_0$ and $\varepsilon_0$ and a global Maxellian $\mb{M}_*=\mb{M}_{[\r_*,u_*,\t_*]}$,  such that if
the shock wave strength and the diffusion wave strength $\alpha_2$ satisfy $|(\r_+-\r_-,m_+-m_-,E_+-E_-)| +|\alpha_2|\leq \delta_0$  and \eqref{delta} and the initial data satisfies that
$$
\mathcal{E}(0) \leq \varepsilon_0,
$$
then the Cauchy problem of the bipolar  VPB system \eqref{VPB} admits a unique global solution $(F_1,F_2)(t,x,v)$ satisfying
$$
\mathcal{E}(t)\leq C(\mathcal{E}(0)+\delta_0^\f12)
$$
for the uniform-in-time positive constant $C$ and the time-asymptotic behaviors
$$
\|\big(\widetilde{F}_1(t,x,v), F_2(t,x,v)\big)\|_{L^\infty_x
L_v^2(\frac1{\sqrt\mb{M_*}})}+ \|(\Pi_x ,n_2)(t,x)\|_{L^\infty_x}
\rightarrow
0, ~~~~{\rm as}~~t\to \infty.
$$
Consequently, it holds that
$$
\|\big(F_A-F^S_{\a_1,\a_3}, F_B-F^S_{\a_1,\a_3}\big)\|_{L^\infty_x
L_v^2(\frac1{\sqrt\mb{M_*}})}+ \|(\Pi_x ,n_2)\|_{L^\infty_x}
\rightarrow
0, ~~~~{\rm as}~~t\to \infty.
$$
Here and in the sequel,  $f(v)\in L_v^2(\f{1}{\sqrt{\mb{M}_*}})$ means that
$\frac{f(v)}{\sqrt{\mb{M_*}}}\in L_v^2(\mb{R}^3)$.
\end{theorem}

\begin{remark}
Theorem \ref{thm0} implies that the linear superposition of two Boltzmann shock profiles in the first and third characteristic fields is nonlinearly stable time-asymptotically to the 1D bipolar VPB system \eqref{VPB} up to some suitable shifts. Note that there is no zero macroscopic mass conditions on the initial perturbations.
\end{remark}

With the above preparation, we will give the proof of the main theorem in the following section.
For this, we will first reformulate the problem.
It follows from \eqref{F1-f}, \eqref{jinsi} and \eqref{perturb-s} that $(\Phi, \Psi, W)$ solves
\begin{equation}
\left\{
\begin{array}{l}
\Phi_{t}+\Psi_{1x}=0,\\
\di \Psi_{1t}+\left(\f{m_1^2}{\r}+p-\f{\tilde {m}_1^2}{\tilde \r}-\tilde p\right)-\f{\Pi_x ^2}{4}=
\f{4}{3}(\mu(\t)u_{1x}-\mu(\tilde\t)\tilde u_{1x})-\int v_1^2(\G-\G^{S_1}-\G^{S_3})d v-Q_1,\\[3mm]
\di\Psi_{it}+\left(\f{m_1m_i}{\r}-\f{\tilde m_1\tilde m_i}{\tilde\r}\right)=
(\mu(\t)u_{ix}-\mu(\tilde\t)\tilde u_{ix})-\int v_1 v_i(\G-\G^{S_1}-\G^{S_3})d v,~~i=2,3,\\[3mm]
\di W_{t}+\left(\f{Em_1}{\r}+\f{pm_1}{\r}-\f{\tilde E\tilde m_1}{\tilde\r}-\f{\tilde p\tilde m_1}{\tilde\r}\right)
 =(\k(\t)\t_{x}-\k(\tilde\t)\tilde\t_{x})\\[2mm]
\di\quad +\f43(\mu(\t)u_1u_{1x}-\mu(\tilde\t)\tilde u_1\tilde u_{1x})
+\sum_{i=2}^3\mu(\t)u_iu_{ix}-\f12\int v_1| v|^2(\G-\G^{S_1}-\G^{S_3})d v-Q_2,
\end{array} \right.\label{new-sys}
\end{equation}
To capture the viscous effect of velocity and temperature, set
\begin{equation}
\Psi=\tilde\r \wt{\Psi}+\tilde u\Phi,
\label{Psi}
\end{equation}
\begin{equation}
W=\tilde\r \wt{W}+\tilde u\cdot\Psi+(\tilde\t-\f{|\tilde
u|^2}{2})\Phi=\tilde\r \wt{W}+\tilde \r\tilde u\cdot\wt{\Psi}
+(\tilde\t+\f{|\tilde u|^2}{2})\Phi.
\label{W}
\end{equation}
We also denote
\begin{equation}\label{tp}
(\tilde\psi, \tilde\omega)(t,x)=(\wt\Psi_x, \wt W_x)(t,x).
\end{equation}
Then we have the following linearized system for $(\Phi,\wt\Psi,\wt W):$
\begin{equation}
 \left\{
\begin{array}{l}
 \Phi_{t}+\tilde\r\wt\Psi_{1x}+\tilde u_1\Phi_{x}+\tilde\r_x\wt\Psi_1+\tilde u_{1x}\Phi=0,\\
\di \tilde\r\wt\Psi_{1t}+\tilde\r\tilde u_1\wt\Psi_{1x}-\f13\tilde\r\tilde
u_{1x}\wt\Psi_1 +\f23\tilde\r_x\wt W+\f23\tilde\r\wt
W_x+\f23\tilde\t\Phi_x
-\f{2\tilde\t\tilde\r_x}{3\tilde\r}\Phi\\
\di\qquad\qquad =\f43\mu(\tilde\t)\wt\Psi_{1xx}-\int v_1^2(\G-\G^{S_1}-\G^{S_3})d v+J_1+N_1-Q_1,\\[3mm]
\di \tilde\r\wt\Psi_{it}+\tilde\r\tilde u_1\wt\Psi_{ix}-\tilde\r\tilde
u_{1x}\wt\Psi_i=\mu(\tilde\t)\wt\Psi_{ixx}-\int v_1 v_i(\G-\G^{S_1}-\G^{S_3})d v+J_i+N_i,~~~i=2,3,\\[3mm]
 \tilde\r\wt W_{t}+\tilde\r\tilde u_1\wt W_x-\tilde\r\tilde u_{1x}\wt W
 +\di\f23 \tilde\r \tilde\t \wt\Psi_{1x}-\di \f23 \tilde\r \tilde\t_x \wt\Psi_1
=\k(\tilde\t)\wt W_{xx}\\[3mm]
\quad\di -\f12\int v_1| v|^2(\G-\G^{S_1}-\G^{S_3})d v+\tilde u_1\int v_1^2(\G-\G^{S_1}-\G^{S_3})dv+J_4+N_4-(Q_2-\tilde u_1Q_1),
\end{array}
\right.\label{P-sys}
\end{equation}
where
\begin{equation}
\begin{array}{ll}
\di J_1=\di
\left[\left(\int v_1^2(\G^{S_1}+\G^{S_3})d v-Q_1\right)_x-\f43\f{\mu(\tilde\t)}{\tilde\r}\tilde\r_x\tilde
u_{1x}\right]\f\Phi{\tilde\r}+\f43\f{\mu(\tilde\t)}{\tilde\r}\tilde u_{1x}\Phi_x\\
\di \qquad +\f43\left(\tilde u_{1x}^2-\f{\tilde\r_x\tilde\t_x}{\tilde\r}\right)\mu^\prime(\tilde\t)\wt\Psi_1+\f43\left(\f{\mu(\tilde\t)}{\tilde\r}\tilde\r_x\wt\Psi_1\right)_x+\f43\f{\mu^\prime(\tilde\t)}{\tilde\r}\tilde\r_x\tilde u_{1x}\wt
W+\f43\mu^\prime(\tilde\t)\tilde u_{1x}\wt W_x,\\
\di
J_i=\left(\f{\mu(\tilde\t)}{\tilde\r}\tilde\r_x\wt\Psi_i\right)_x-\f{\mu^\prime(\tilde\t)}{\tilde\r}\tilde\r_x\tilde
\t_{x}\wt \Psi_i,~~i=2,3,\\
\di J_4\di=\Big(\int v_1\f{|v|^2}{2}(\G^{S_1}_x+\G^{S_3}_x)d v-\tilde
u_1\int v_1^2(\G^{S_1}_x+\G^{S_3}_x)d v-Q_{2x}+\tilde u_1Q_{1x}-\f{\k(\tilde\t)}{\tilde\r}\tilde\r_x\tilde
\t_{x}\Big)\f\Phi{\tilde\r}\\
\di\qquad+\f{\k(\tilde\t)}{\tilde\r}\tilde \t_{x}\Phi_x+\left(\int v_1^2(\G^{S_1}_x+\G^{S_3}_x)dv-Q_{1x}+\f43\f{\mu(\tilde\t)}{\tilde\r}\tilde
u_{1x}\tilde \r_{x}\right)\wt\Psi_1+\f83\mu(\tilde\t)\tilde u_{1x}\wt\Psi_{1x}\\
\di\qquad
+\left[\left(\k(\tilde\t)-\f43\mu(\tilde\t)\right)\tilde u_{1x}\wt\Psi_1\right]_x
+\left(\f{\k(\tilde\t)}{\tilde\r}\tilde\r_x\wt
W\right)_x+\k^\prime(\tilde\t)\tilde \t_{x}\wt W_x,
\end{array}
\label{J4}
\end{equation}
and the nonlinear terms
\begin{equation}
N_i=O(1)|(\Phi_x,\Psi_x,W_x,\Pi_x ,\Phi_{xx},\Psi_{ixx})|^2,~~i=1,2,3,
\label{Ni}
\end{equation}
\begin{equation}
N_4=O(1)|(\Phi_x,\Psi_x,W_x,\Pi_x ,\Phi_{xx},\Psi_{xx},W_{xx})|^2.
\label{N4}
\end{equation}
From \eqref{F1-G}, \eqref{shock-non-fluid} and \eqref{perturb-s}, we can derive the equation for the non-fluid component $\wt{\mb{G}}(t,x, v)$ as
\begin{equation}
\begin{array}{ll}
\wt{\mb{G}}_{t}-\mb{L}_\mb{M}\wt{\mb{G}}=\di
-\mb{P}_1( v_1\wt{\mb{G}}_{x})-\mb{P}_1(\Pi_x \partial_{v_1}F_2)+2Q(\wt{\mb{G}},\wt{\mb{G}})
+2[Q(\wt{\mb{G}},\mb{G}^{S_1}+\mb{G}^{S_3})+Q(\mb{G}^{S_1}+\mb{G}^{S_3},\wt{\mb{G}})]\\[3mm]
\di+2[Q({\mb{G}^{S_1}},\mb{G}^{S_3})+Q(\mb{G}^{S_3},{\mb{G}^{S_1}})]-[\mb{P}_1(v_1\mb{M}_x)-\mb{P}^{S_1}_1(v_1\mb{M}^{S_1}_x)-\mb{P}^{S_3}_1(v_1\mb{M}^{S_3}_x)]+\sum_{i=1,3}R_i,\\
\end{array}\label{G}
\end{equation}
with $R_i$ given by
\begin{equation}
R_i=(\mb{L}_{\mb{M}}-\mb{L}_{\mb{M}^{S_i}})\mb{G}^{S_i}-[\mb{P}_1(v_1\mb{G}^{S_i}_x)-\mb{P}^{S_i}_1(v_1\mb{G}^{S_i}_x)]
,~~i=1,3.\label{Ri}
\end{equation}
From \eqref{VPB} and \eqref{BS}, we have the equation for $\wt F_1$ defined in \eqref{perturb-s}
\begin{equation}
\begin{array}{ll}
\wt F_{1t}+ v_1\wt F_{1x}+\Pi_x \partial_{v_1} F_2 =\mb{L}_\mb{M}\wt{\mb{G}}
+2Q(\wt{\mb G},\wt{\mb G})+(\mb{L}_\mb{M}-\mb{L}_{\mb {M}^{S_1}})(\mb {G}^{S_1})
+(\mb{L}_\mb{M}-\mb{L}_{\mb {M}^{S_3}})(\mb {G}^{S_3})\\[3mm]
\di \qquad\qquad\quad +2[Q(\wt{\mb G},\mb {G}^{S_1}+\mb {G}^{S_3})+Q(\mb {G}^{S_1}+\mb {G}^{S_3},\wt{\mb G})]
+2[Q(\mb {G}^{S_1},\mb {G}^{S_3})+Q(\mb {G}^{S_3},\mb {G}^{S_1})], \label{tidle-f}
\end{array}
\end{equation}
Consider the reformulated system \eqref{P-sys}, \eqref{G}, \eqref{tidle-f} and \eqref{n21}, \eqref{n2}, \eqref{F2-pc-1}. Since
the local existence of solution to the VPB system can be proved similarly as in \cite{Guo1}, to prove the global
existence on the time interval $[0, T]$, we only need to
close the following a priori assumption by the continuity argument:
\begin{equation}\label{assumption-ap-s}
\di \mathcal{E}(T)\leq \chi_{\scriptscriptstyle T} ^2,
\end{equation}
where $\mathcal{E}(T)$ is defined in \eqref{et} . Here and in the sequel $\chi_{\scriptscriptstyle T} $ is a small positive constant depending on the initial data
and wave strengths. Note that even $\chi_{\scriptscriptstyle T}$ is denoted by the subscription $T$, which means that the a priori assumption \eqref{assumption-ap-s} is imposed on the time interval $[0,T]$, however, $\chi_{\scriptscriptstyle T}$ is to be chosen to be independent of the time $T$.

\subsection{The Proof of Main Result}

By the continuum argument for the local solution to the system \eqref{VPB} or equivalently the system \eqref{F1-f}-\eqref{F2-pc-1}, to prove Theorem \ref{thm0}, it is sufficient to close the a priori assumption \eqref{assumption-ap-s} and verify the time-asymptotic behaviors of the solution. We start from the lower order estimates in the following Proposition.

\begin{proposition} \label{Prop3.1}
For each $(\r_-,m_-,E_-)$, there exists a positive constant $C$ such that, if
$\d+|\a_2|\leq \d_0$ for a suitably small positive constant $\d_0$, then it holds that
\begin{equation*}\label{LE-A}
\begin{array}{ll}
\di \|(\Phi,\Psi,W,\Phi_x,\Pi_x ,n_2)(\cdot,t)\|^2+\int\int\f{|(\widetilde{\mb{G}},\mb{P}_c F_2)|^2}{\mb{M}_*}(x,v,t) dxdv+\int_0^t\|(\Pi_x , \Pi_{x\tau}, n_2)\|^2d\tau\\
\di +\sum_{|\alpha^\prime|=1}\int_0^t\|\partial^{\alpha^\prime}(\Phi,\Psi,W,\wt\Psi,\wt W,n_2)\|^2d\tau+\int_0^t\|\sqrt{|u^{S_1}_{1x}|+|u^{S_3}_{1x}|+|\Theta_x|}(\Phi,\wt{\Psi}_1,\wt W)\|^2d\tau\\
\di+\int_0^t\int\int\f{\nu(|v|)|(\widetilde{\mb{G}},\mb{P}_c F_2)|^2}{\mb{M}_*}dxdvd\tau \leq C(\chi_{\scriptscriptstyle T} +\delta_0)\sum_{|\alpha^\prime|=1}\int_0^t\|\partial^{\alpha^\prime}(\phi,\psi,\o)\|^2d\tau\\
\di +C\int_0^t\|(\phi_x,\wt\psi_x,\wt\o_x)\|^2d\tau+C(\chi_{\scriptscriptstyle T} +\delta_0)\int_0^t\int\int\f{\nu(|v|)|(\widetilde{\mb{G}},\mb{P}_c F_2)_{v_1}|^2}{\mb{M}_*} dxdvd\tau\\
\di +C\sum_{|\alpha^\prime|=1}\int_0^t\int\int \f{\nu(|v|)|\partial^{\alpha^\prime}(\wt{\mb{G}}, \mb{P}_cF_2)|^2}{\mb{M}_*}dxdvd\tau+C(\mathcal{E}(0)^2+\delta_0^{\f12}).
\end{array}
\end{equation*}
\end{proposition}
The proof of Proposition \ref{Prop3.1} will be given in Appendix A.
Then we perform the higher order estimates. Firstly, we apply
$\partial_x$ to the system \eqref{P-sys} to get
\begin{equation}
\left\{
\begin{array}{ll}
\di \phi_t+\tilde\r\tilde\psi_{1x}+\tilde u_1\phi_x+\tilde\r_x\tilde\psi_1+\tilde u_{1x}\phi=-L_0,\\
\di \tilde\r\tilde\psi_{1t}+\tilde\r\tilde u_1\tilde\psi_{1x}-\f13\tilde\r\tilde
u_{1x}\tilde\psi_1 +\f23\tilde\r_x\tilde \omega+\f23\tilde\r\tilde
\omega_x+\f23\tilde\t\phi_x
-\f{2\tilde\t\tilde\r_x}{3\tilde\r}\phi\\
\di\qquad=\left(\f43\mu(\tilde\t)\tilde\psi_{1x}\right)_x-\int v_1^2(\G-\G^{S_1}-\G^{S_3})_x d v+(J_1+N_1-Q_1)_x-L_1,\\
\di \tilde\r\tilde\psi_{it}+\tilde\r\tilde u_1\tilde\psi_{ix}-\tilde\r\tilde
u_{1x}\tilde\psi_i=
(\mu(\tilde\t)\tilde\psi_{ix})_x\\
\di\qquad -\int v_1 v_i(\G-\G^{S_1}-\G^{S_3})_xd v+(J_i+N_i)_x-L_i,\qquad  i=2,3,\\
\di \tilde\r\tilde \omega_{t}+\tilde\r\tilde u_1\tilde \omega_x-\tilde\r\tilde u_{1x}\tilde  \omega
+\f23\tilde\r\tilde\t\tilde\psi_{1x}-\f23\tilde\r\tilde\t_x\tilde\psi_1=(\k(\tilde\t)\tilde \omega_{x})_x-\f12\int v_1| v|^2(\G-\G^{S_1}-\G^{S_3})_xd v
\\
\di\qquad +\left(\tilde u_1\int v_1^2(\G-\G^{S_1}-\G^{S_3})dv\right)_x+\big(J_4+N_4-(Q_2-\tilde u_1Q_1)\big)_x-L_4,
\end{array} \right.\label{sys-h-s}
\end{equation}
where
\begin{equation*}
\begin{array}{ll}
L_0=\tilde\r_x\tilde\psi_1+\tilde u_{1x}\phi+\tilde\r_{xx}\wt\Psi_1+\tilde u_{1xx}\Phi,
\end{array}\label{L0}
\end{equation*}
\begin{equation*}
\begin{array}{ll}
\di L_1=\tilde\r_x\wt\Psi_{1t}+(\tilde\r\tilde u_1)_x\tilde\psi_1-\di\f13(\tilde\r\tilde u_{1x})_x\wt\Psi_1
+\f23\tilde\r_{xx}\wt W+\f23\tilde\r_x\tilde\o+\f23\tilde\t_x\phi
 -\f23\left(\f{\tilde\t\tilde\r_x}{\tilde\r}\right)_x\Phi,\\
\di
L_i=\tilde\r_x\wt\Psi_{it}+(\tilde\r\tilde u_1)_x\tilde\psi_i-(\tilde\r\tilde u_{1x})_x\wt\Psi_i,~~~~i=2,3,\\
\di L_4=\tilde\r_x\wt W_{t}+(\tilde\r\tilde u_1)_x\tilde\o-(\tilde\r\tilde u_{1x})_x\wt W
+\di\f23(\tilde\r\tilde\t)_x\tilde\psi_1-\di\f23(\tilde\r\tilde\t_x)_x\wt\Psi_1.
\end{array}\label{L4}
\end{equation*}

To derive the estimate on the more higher order derivatives,
we apply $\partial_x$ to the system \eqref{sys-h-s} to obtain
\begin{equation}
\left\{
\begin{array}{ll}
\di \phi_{xt}+\tilde\r\tilde\psi_{1xx}+\tilde u_1\phi_{xx}+\tilde\r_x\tilde\psi_{1x}+\tilde u_{1x}\phi_x=-\wt L_0,\\[2mm]
\di \tilde\r\tilde\psi_{1xt}+\tilde\r\tilde u_1\tilde\psi_{1xx}-\f13\tilde\r\tilde
u_{1x}\tilde\psi_{1x} +\f23\tilde\r_x\tilde \omega_x+\f23\tilde\r\tilde
\omega_{xx}+\f23\tilde\t\phi_{xx}
-\f{2\tilde\t\tilde\r_x}{3\tilde\r}\phi_x\\
\di\qquad=\left(\f43\mu(\tilde\t)\tilde\psi_{1x}\right)_{xx}-\int v_1^2(\G-\G^{S_1}-\G^{S_3})_{xx} d v+(J_1+N_1-Q_1)_{xx}-\wt L_1,\\[2mm]
\di \tilde\r\tilde\psi_{ixt}+\tilde\r\tilde u_1\tilde\psi_{ixx}-\tilde\r\tilde
u_{1x}\tilde\psi_{ix}=
(\mu(\tilde\t)\tilde\psi_{ix})_{xx}-\int v_1 v_i(\G-\G^{S_1}-\G^{S_3})_{xx}d v\\
\di\qquad+(J_i+N_i)_{xx}-\wt L_i,\qquad i=2,3,\\[2mm]
\di \tilde\r\tilde \omega_{xt}+\tilde\r\tilde u_1\tilde \omega_{xx}-\tilde\r\tilde u_{1x}\tilde  \omega_x+\f23\tilde\r\tilde\t\tilde\psi_{1xx}-\f23\tilde\r\tilde\t_x\tilde\psi_{1x}
=(\k(\tilde\t)\tilde \omega_{x})_{xx}\\[2mm]
\di\qquad-\f12\int v_1| v|^2(\G-\G^{S_1}-\G^{S_3})_{xx}d v
+\left(\tilde u_1\int v_1^2(\G-\G^{S_1}-\G^{S_3})dv \right)_{xx}\\
\di\qquad+(J_4+N_4-(Q_2-\tilde u_1Q_1))_{xx}-\wt L_4,
\end{array} \right.\label{sys-h-h-s}
\end{equation}
where
\begin{equation*}
\begin{array}{ll}
\di \wt L_0=\tilde\r_x\tilde\psi_{1x}+\tilde u_{1x}\phi_x+\tilde\r_{xx}\tilde\psi_1+\tilde u_{1xx}\phi+L_{0x},\\
\di \wt L_1=\tilde\r_x\tilde\psi_{1t}+(\tilde\r\tilde u_1)_x\tilde\psi_{1x}-\di\f13\left(\tilde\r\tilde u_{1x}\right)_x\tilde\psi_1
+\f23\tilde\r_{xx}\tilde \o+\f23\tilde\r_x\tilde\o_x+\f23\tilde\t_x\phi_x-\f23(\f{\tilde\t\tilde\r_x}{\tilde\r})_x\phi+L_{1x},\\
\di
\wt L_i=\tilde\r_x\tilde\psi_{it}+(\tilde\r\tilde u_1)_x\tilde\psi_{ix}-(\tilde\r\tilde u_{1x})_x\tilde\psi_i
+L_{ix},~~~~i=2,3,\\
\di
\wt L_4=\tilde\r_x\tilde \o_{t}+(\tilde\r\tilde u_1)_x\tilde\o_x-(\tilde\r\tilde u_{1x})_x\tilde \o
+\di\f23(\tilde\r\tilde\t)_x\tilde\psi_{1x}-\di\f23(\tilde\r\tilde\t_x)_x\tilde\psi_1+L_{4x}.
\end{array}\label{L4+}
\end{equation*}

By using the above two systems and the
equations for $n_2$ and the non-fluid component $\mb{P}_c F_2$, we can establish the following proposition
for the
higher order energy estimates.

\begin{proposition}\label{Prop3.2}
Under the assumptions of Proposition \ref{Prop3.1},  it holds that
\begin{equation*}
\begin{array}{ll}
\di\sup_{0\leq t<\i} \Big[\|(\phi,\psi,\omega)(\cdot,t)\|^2_{H^1}+\|(\Pi_x ,n_2,n_{2x})(\cdot,t)\|^2+\sum_{0\leq|\beta|\leq2}\int\int\f{|\partial^\beta(\widetilde{\mb{G}},\mb{P}_c F_2)|^2}{\mb{M}_*}(x,v,t) dxdv\\
\di +\sum_{|\alpha^\prime|=1,0\leq|\beta^\prime|\leq1}\int\int \f{|\partial^{\alpha^\prime}\partial^{\beta^\prime}(\wt{\mb{G}},\mb{P}_c F_2)|^2}{\mb{M}_*}(x,v,t)dxdv+\sum_{|\alpha|=2}\int\int \f{|\partial^\alpha(\wt F_1,F_2)|^2}{\mb{M}_*}(x,v,t) dxdv\Big] \\
\di +\sum_{1\leq|\alpha|\leq2}\int_0^\i\|\partial^{\alpha}(\phi,\psi,\omega,n_2)\|^2d\tau+\int_0^\i\|(\Pi_x ,n_2)\|^2d\tau\\
\di +\sum_{1\leq|\alpha|\leq2}\int_0^\i\int\int \f{\nu(|v|)|\partial^\alpha(\wt{\mb{G}},\mb{P}_c F_2)|^2}{\mb{M}_*}dxdvd\tau+\sum_{0\leq|\beta|\leq2}\int_0^\i\int\int\f{\nu(|v|)|\partial^\beta(\widetilde{\mb{G}},\mb{P}_c F_2)|^2}{\mb{M}_*}dxdvd\tau\\
\di +\sum_{|\alpha^\prime|=1,|\beta^\prime|=1}\int_0^\i\int\int \f{\nu(|v|)|\partial^{\alpha^\prime}\partial^{\beta^\prime}(\wt{\mb{G}},\mb{P}_c F_2)|^2}{\mb{M}_*}dxdvd\tau\leq C(\mathcal{E}(0)^2+\delta_0^{\f12}).
\end{array}
\end{equation*}
\end{proposition}

The proof of Proposition \ref{Prop3.2} will be given in the Appendix A.

With Proposition \ref{Prop3.1} and Proposition \ref{Prop3.2}, we can close the a priori assumption \eqref{assumption-ap-s} and one has
\begin{equation*}
\begin{array}{ll}\label{aa}
\qquad \di \int_0^{+\infty}\int\int\frac{|(\wt F_1, F_2)_x|^2}{\mb{M}_*} dvdxd\tau\\[2mm]
\qquad\di\leq \int_0^{+\infty}\int\int\frac{|\big(\mb{M}- \mb{M}^{S_1}-\mb{M}^{S_3}+\mb{M}_\#, \f{n_2}{\rho}\mb{M}\big)_x|^2}{\mb{M}_*} dvdxd\tau
+\int_0^{+\infty}\int\int\frac{|(\wt{\mb{G}},\mb{P}_cF_2)_x|^2}{\mb{M}_*} dv dxd\tau\\[2mm]
\qquad\di\leq  C(\mathcal{E}(0)^2+\d_0^{\frac12}).
\end{array}
\end{equation*}
From the Vlasov-Poisson-Boltzmann system, we  can obtain
\begin{equation*}
\begin{array}{ll}\label{bb}
&\di \int_0^{+\infty}\Big|\frac{d}{dt}\int\int\frac{|(\wt F_1, F_2)_x|^2}{\mb{M}_*} dvdx\Big|d\tau
\leq C(\mathcal{E}(0)^2+\d_0^{\frac12}).
\end{array}
\end{equation*}
Therefore, one has
\begin{equation*}
\begin{array}{ll}
\di\int_0^{+\infty}\Big(\int\int\frac{|(\wt F_1, F_2)_x|^2}{\mb{M}_*} dv dx
+\Big|\frac{d}{dt}\int\int\frac{|(\wt F_1, F_2)_x|^2}{\mb{M}_*} dv dx\Big|
\Big)d\tau<\infty,
\end{array}
\end{equation*}
which implies that
$$
\lim_{t\rightarrow+\infty} \int\int\frac{|(\wt F_1, F_2)_x|^2}{\mb{M}_*} dv dx=0.
$$

\

By Sobolev inequality
$$
\begin{array}{ll}
\di \|\int\frac{|(\wt F_1, F_2)|^2}{\mb{M}_*} dv\|^2_{L^{\infty}_x} \leq C\int\int\frac{|(\wt F_1, F_2)|^2}{\mb{M}_*} dv dx \int\int\frac{|(\wt F_1, F_2)_x|^2}{\mb{M}_*} dv dx
\end{array}
$$
we can prove the asymptotic behavior of the solutions and complete the proof of Theorem \ref{thm0}.

%%%%%%%%%%%%%%%%%%%%%%%%%%%%%%%%%%%%
%
%
%%%%%%%%%%%%%%%%%%%%%%%%%%%%%%%%%%%
\section{Stability of Rarefaction Wave to the Bipolar VPB System}
\label{result}
\setcounter{equation}{0}

In this section, we employ the micro-macro decomposition introduced in Section 2 to prove the
nonlinear stability of planar rarefaction wave to the Cauchy problem of the bipolar VPB system \eqref{VPB10}-\eqref{VPB10-i} in spatial one-dimension
by the two states $(\rho_\pm,u_\pm,\t_\pm)$ with $u_\pm=(u_{1\pm},0,0)^t$ and $\rho_\pm>0, u_{1\pm}, \t_\pm>0$ being connected by the rarefaction wave solution to the Riemann problem of the corresponding 1D Euler system
\eqref{euler}
with the Riemann initial data
\begin{equation}\label{R-in}
(\r_0^r,u_0^r,\t_0^r)(x)=\left\{
\begin{array}{ll}
\di (\r_+,u_+,\t_+), &\di x>0, \\
\di (\r_-,u_-,\t_-), &\di x<0.
\end{array}
 \right.
\end{equation}
Correspondingly, the initial values to the transformed system \eqref{VPB1} satisfy
\begin{equation}  \label{VPB1a}
\left\{\begin{array}{ll}
\di F_1(t=0,x,v)=F_{10}(x,v)\rightarrow \mathbf{M}_{[\rho_\pm,u_\pm,\t_\pm]}(v),\quad {\rm as}~~x\rightarrow\pm\infty,\\
\di F_2(t=0,x,v)=F_{20}(x,v)\rightarrow 0,\quad {\rm as}~~x\rightarrow\pm\infty,\\
\di \Pi_x\rightarrow 0,\quad {\rm as}~~x\rightarrow\pm\infty.
\end{array}
\right.\end{equation}

\subsection{Rarefaction Wave and Main Result}
First we describe the rarefaction wave solution to the Euler system \eqref{euler}, \eqref{R-in} with the state equation
$$
p=\frac{2}{3}\rho\theta=k\rho^{\frac53}\exp(S),\qquad k=\frac{1}{2\pi e}.
$$
It is straight to calculate that the Euler system \eqref{euler} for $(\r,u_1,\t)$ has three distinct eigenvalues
$$
\lambda_i(\r,u_1,S)=u_1+(-1)^{\f{i+1}{2}}\sqrt{p_\r(\r,S)},~ i=1,3,\qquad
\lambda_2(\r,u_1,S)=u_1,~
$$
with corresponding right eigenvectors
$$
r_i(\r,u_1,S)=((-1)^{\f{i+1}{2}}\r,\sqrt{p_\r(\r,S)},0)^t,~i=1,3, \qquad\quad r_2(\r,u_1,S)=(p_{\scriptscriptstyle S},0,-p_\rho)^t,$$
such that
$$
r_i(\r,u_1,S)\cdot \nabla_{(\r,u_1,S)}\lambda_i(\r,u_1,S)\neq 0,~ i=1,3,\quad
{\rm
and}\quad
r_2(\r,u_1,S)\cdot \nabla_{(\r,u_1,S)}\lambda_2(\r,u_1,S)\equiv 0.
$$
Thus the  two $i$-Riemann invariants $\Sigma_i^{(j)}(i=1,3, j=1,2)$ can be defined by (cf. \cite{Smoller})
\begin{equation}\label{RI}
\Sigma_i^{(1)}=u_1+(-1)^{\f{i-1}{2}}\int^{\r}\f{\sqrt{p_z(z,S)}}{z}dz,\qquad
\Sigma_i^{(2)}=S,
\end{equation}
such that
$$
\nabla_{(\r,u_1,S)} \Sigma_i^{(j)}(\r,u_1,S)\cdot r_i(\r,u_1,S)\equiv0,\quad i=1,3,~ j=1,2.
$$
Given the right state $(\rho_+, u_{1+}, \theta_+)$ with $\r_+>0,\t_+>0$,  the $i$-Rarefaction wave curve $(i=1,3)$ in the phase space $(\r,u_1,\t)$ with $\r>0$ and $\t>0$ can be defined by (cf. \cite{Lax}):
\begin{equation}
 R_i (\rho_+, u_{1+}, \theta_+):=\Bigg{ \{} (\rho, u_1, \theta)\Bigg{ |}\lambda_{ix}(\r,u_1,S)>0, \Sigma_i^{(j)}(\r,u_1,S)=\Sigma_i^{(j)}(\r_+,u_{1+},S_+),~~j=1,2\Bigg{ \}}\label{(2.2)}.
\end{equation}
Without loss of generality, we consider stability of $3-$rarefaction wave to the Euler system \eqref{euler}, \eqref{R-in} in the present paper and the stability of $1-$rarefaction wave can be done similarly.  The $3-$rarefaction wave to the Euler system \eqref{euler}, \eqref{R-in} can be expressed explicitly by
 the Riemann solution to the inviscid Burgers equation:
\begin{equation}\label{bur}
\left\{\begin{array}{ll}
w_t+ww_x=0,\\
w(x,0)=\left\{\begin{array}{ll}
w_-,&x<0,\\
w_+,&x>0.
\end{array}
\right.
\end{array}
\right.
\end{equation}
If $w_-<w_+$, then the Riemann problem $(\ref {bur})$ admits a
rarefaction wave solution $w^r(x, t) = w^r(\f xt)$ given by
\begin{equation}\label{abur}
w^r(\f xt)=\left\{\begin{array}{lr}
w_-,&\f xt\leq w_-,\\
\f xt,&w_-\leq \f xt\leq w_+,\\
 w_+,&\f xt\geq w_+.
\end{array}
\right.
\end{equation}
Then the 3-rarefaction wave solution $(\r^r,u^r,\t^r)(\f xt)$ to the
compressible Euler equations \eqref{euler}, \eqref{R-in} can be defined explicitly by
\begin{eqnarray}   \label{3-rw}
\left\{
\begin{array}{l}
\di w_\pm=\lambda_3(\r_\pm,u_{1\pm},\t_\pm), \qquad w^r(\f xt)= \lambda_3(\r^r,u_1^r,\t^r)(\f xt),\\
\di
\Sigma_3^{(j)}(\r^r,u_1^r,\t^r)(\f xt)=\Sigma_3^{(j)}(\rho_\pm,u_{1\pm},\t_\pm),\quad j=1,2,\quad
u_2^r= u_3^r=0,
\end{array} \right.
\end{eqnarray}
where $\Sigma_3^{(j)}~(j=1,2)$ are the 3-Riemann invariants defined in \eqref{RI}.

By the previous micro-macro decomposition for $F_1=\mb{M}+\mb{G}$, one has the
fluid-part system \eqref{F1-f}
and the non-fluid equation  \eqref{F1-G} for $F_1$. By the new micro-macro decomposition to
$$
F_2=\frac{n_2}{\rho}\mb{M}+\mb{P}_c F_2,
$$
the macroscopic part $n_2$ satisfy the nonlinear diffusion equation
\eqref{n21}, or equivalently
\eqref{n2}, and the microscopic part $\mb{P}_c F_2$ satisfy the equation
\eqref{F2-pc-1}.

The analysis of nonlinear stability towards the rarefaction wave to the bipolar VPB system \eqref{VPB1} will be carried out for the equivalent system \eqref{F1-f}-\eqref{F2-pc-1}, which can be roughly viewed as the strong couplings of the system of viscous conservation laws with the microscopic terms and the electric field terms.

Denote the energy functional $\mathcal{E}(t)$ by
\begin{equation}\label{Et}
\begin{array}{ll}
\di \mathcal{E}(t)=
\di \sup_{\tau\in[0,t]}\bigg\{\|(\rho-\r^r,u-u^r,\t-\t^r)\|^2_{H^1(\mb{R})}+\|(\Pi_{x},n_{2},n_{2x})\|^2\\
\di
\qquad\quad +\sum_{0\leq|\beta|\leq2}\int_{\mb{R}}\int\f{|\partial^\beta(\mb{G},\mb{P}_cF_{2})|^2}{\mb{M}_*}dv
dx+\sum_{|\alpha|=1,0\leq|\beta|\leq1}\int\int\f{|\partial^\alpha\partial^\beta(\mb{G},\mb{P}_cF_{2})|^2}{\mb{M}_*}dv
dx
\\
\di\qquad\quad  +\sum_{|\alpha|=2}\int\int\f{|\partial^\alpha(F_{1},F_{2})|^2}{\mb{M}_*}dv dx\bigg\},
\end{array}
\end{equation}
where and in the sequel $\partial^\alpha=\partial_{x,t}^\alpha$, $\partial^\beta=\partial_{v}^\beta$.
Note that the global Maxellian $\mb{M}_*$ in \eqref{Et} is determined in Theorem \ref{thm}. Then the main result for the stability of rarefaction wave to bipolar VPB system \eqref{VPB1} can be stated as follows.

\begin{theorem}\label{thm}
Assume that Riemann solution of Euler equation \eqref{euler}, \eqref{R-in} consists of one $3$-rarefaction wave given by \eqref{3-rw}. Then, there exist positive constants $\delta_0$ and $\varepsilon_0$ and a global Maxellian $\mb{M}_*=\mb{M}_{[\r_*,u_*,\t_*]}$ with $\r_*>0,\t_*>0$, such that if
the wave strength $\delta=|(\r_+-\r_-,u_+-u_-,\t_+-\t_-)| \leq \delta_0$ and the initial data satisfies
\begin{equation}
\mathcal{E}(0)\leq \varepsilon_0,
\end{equation}
then Cauchy problem of the bipolar VPB system \eqref{VPB1}--\eqref{VPB1a} admits a unique global strong solution $(F_1,F_2,\Phi)$ satisfying
\begin{equation}\label{290}
\mathcal{E}(t)\leq C(\mathcal{E}(0)+\delta_0^{\f18}),
\end{equation}
with the positive constant $C$ independent of $t\geq0$, and the time-asymptotic behaviors
$$
\begin{array}{ll}
\di \|\big(F_1(t,x,v)-\mb{M}_{[\r^r,u^r,\t^r]}(t,x,v), F_2(t,x,v)\big)\|_{L^\infty_x
L_v^2(\frac1{\sqrt{\mb{M_*}}})}
\rightarrow
0, ~~~~{\rm as}~~t\to \infty,\\
\di \|(\r, u, \t)(t,x)-(\r^r,u^r,\t^r)(\frac x t)\|_{L^\infty_x}+ \|(\Pi_x,n_2)(t,x)\|_{L^\infty_x}\rightarrow
0, ~~~{\rm as}~~t\to \infty.
\end{array}
$$
Here $f(v)\in L_v^2(\f{1}{\sqrt{\mb{M}_*}})$ means that
$\frac{f(v)}{\sqrt{\mb{M_*}}}\in L_v^2(\mb{R}^3)$.

\begin{remark}
From Theorem~\ref{thm} it follows that there exists a unique global strong solution $(F_A,F_B,\Pi_x)$ to the original bipolar VPB system \eqref{VPB10}-\eqref{VPB10-i} which satisfies
$$
\|\big(F_A(t,x,v)-\mb{M}_{[\r^r,u^r,\t^r]}(t,x,v), F_B(t,x,v)-\mb{M}_{[\r^r,u^r,\t^r]}(t,x,v)\big)\|_{L^\infty_x
L_v^2(\frac1{\sqrt{\mb{M_*}}})}+ \|(\Pi_x,n_2)(t,x)\|_{L^\infty_x}
\rightarrow
0,
$$
as $t\to \infty$, which implies the time-asymptotic stability of rarefaction wave to the inviscid Euler system for the 1D bipolar VPB system \eqref{VPB}. Note that the time-asymptotic stability of rarefaction wave in Theorem \ref{thm} is independent of the approximation for the rarefaction wave fan in the following subsection 4.2.
\end{remark}
\begin{remark}
The initial values $\Pi_{0x}$ is defined through the Poisson equation $\Pi_{0xx}=2n_{20}$, while $\di n_{20}=\int F_{20}(x,v) dv$.
\end{remark}

\end{theorem}

\subsection{Approximate Rarefaction Wave}

We first construct an approximate smooth rarefaction wave to the 3-rarefaction wave defined in \eqref{3-rw}.  Motivated by \cite{MN-86},  the approximate rarefaction wave can be constructed
by the Burgers equation
\begin{equation}\label{dbur}
\left\{
\begin{array}{l}
\di \bar w_{t}+\bar w\bar w_{x}=0,\\
\di \bar w( 0,x
)=\bar w_0(x)=\f{w_++w_-}{2}+\f{w_+-w_-}{2}\tanh x,
\end{array}
\right.
\end{equation}
where the hyperbolic tangent function $\tanh x=\f{e^x-e^{-x}}{e^x+e^{-x}}$. Note that the solution $\bar w(t,x)$ of the
problem \eqref{dbur} can be given explicitly by
\begin{equation}\label{b-s}
\bar w(t,x)=\bar w_0(x_0(t,x)),\qquad x=x_0(t,x)+\bar w_0(x_0(t,x))t.
\end{equation}
And $\bar w(t,x)$ has the following properties:

\begin{lemma}[\cite{MN-86}]
\label{appr} The problem~$(\ref{dbur})$ has a unique smooth global
solution $\bar w(x,t)$ such that
\begin{itemize}
\item[(1)] $w_-<\bar w(x,t)<w_+, \ \partial_x \bar w(x,t)>0,$
 \ for  $x\in\mathbf{R}, \ t\geq 0.$
\item[(2)] The following estimates hold for all $\ t> 0$ and
p $\in[1,\infty]$:
\begin{equation*}
\begin{array}{ll}
\di \|\bar w(\cdot,t)-w^r(\frac\cdot t)\|_{L^p}\leq  C(w_+-w_-),\\[3mm]
\di
\|\partial_x \bar w(\cdot,t)\|_{L^p}\leq  C\min\{(w_+-w_-),
(w_+-w_-)^{1/p}t^{-1+1/p}\}, \\[3mm]
\di
  \|\partial^2_x \bar w(\cdot,t)\|_{L^p}\leq
C\min\{(w_+-w_-),t^{-1}\},\\[3mm]
\di
|\frac{\partial^2 \bar w(x,t)}{\partial
x^2}|\leq C\frac{\partial \bar w(x,t)}{\partial
x}.
\end{array}
\end{equation*}
\item[(3)] The approximate rarefaction wave $\bar w(x,t)$ and the original rarefaction wave $w^r(\frac xt)$ are time-asymptotically equivalent, i.e.,
$$
\lim_{t\rightarrow+\infty}\sup_{x\in \mb{R}} |\bar w(x,t)-w^r(\f x t)|=0.
$$
\end{itemize}
\end{lemma}

Correspondingly, the approximate rarefaction wave
$(\bar{\r},\bar{u},\bar{\t})(x,t)$ to the  3-rarefaction wave
$(\rho^{r},u^{r},\t^{r})(\f xt)$ in \eqref{3-rw} to compressible Euler equations
$\eqref {euler}, \eqref{R-in}$ can be defined by
\begin{eqnarray}
\left\{
\begin{array}{l}
\di w_\pm=\lambda_3(\r_\pm,u_{1\pm},\t_\pm), \qquad \bar w(x,1+t)= \lambda_3(\bar{\r},\bar{u}_1,\bar{\t})(x,t),\\
\di
\Sigma_3^{(j)}(\bar{\r},\bar{u}_1,\bar{\t})(x,t)=\Sigma_3^{(j)}(\rho_\pm,u_{1\pm},\t_\pm),\quad j=1,2,\quad
\bar u_2=\bar u_3=0,
\end{array} \right.\label{au}
\end{eqnarray}
where $\bar w(x,t)$ is the solution of Burger's equation $(\ref
{dbur})$ defined in \eqref{b-s}.
Then the approximate rarefaction wave $(\bar{\r},\bar{u},\bar{\t})(x,t)$ satisfies the Euler system
\begin{equation}\label{rare-s}
\left\{
\begin{array}{ll}
\di \bar\r_t+(\bar \r \bar u_1)_x=0,\\
\di (\bar \r \bar u_1)_t+(\bar \r\bar u_1^2+\bar p)_x=0,\\
\di  (\bar \r \bar u_i)_t+(\bar \r\bar u_1\bar u_i)_x=0,\qquad i=2,3,\\
\di (\bar\r\bar\t)_t+(\bar\r \bar u_1\bar\t)_x+\bar p\bar u_{1x}=0,
\end{array}
\right.
\end{equation}
and the following properties, which can be proven by similar arguments as used in \cite{MN-86} and is omitted for brevity:

\begin{lemma}%(see $\cite {mn86, xin93}$)
\label{appu} The approximate 3-rarefaction wave $(\bar
\r,\bar u, \bar\t)$ defined in (\ref{au}) satisfies the following
properties:
\begin{itemize}
\item[(i)] $\bar u_{1x}(x,t)>0,$
 \ for  $x\in\mathbf{R}, \ t\geq 0.$
\item[(ii)] The following estimates hold for all $t> 0$ and
q $\in[1,\infty]$:
\begin{equation*}\begin{array}{l}
\di \|(\bar \r, \bar u_1,\bar \t)(\cdot,t)-(\r^r,u^r,\t^r)(\frac \cdot t)\|_{L^q}\leq C\delta, \\
\di \| (\bar\r,\bar u_1,\bar\t)_x(\cdot,t)\|_{L^q}\leq
C_q ~\delta^{1/q}(1+t)^{-1+1/q},\\[4mm]
  \|(\bar\r,\bar u_1,\bar\t)_{xx}(\cdot,t)\|_{L^q}\leq
C_q \min\{\delta^{-1},(1+t)^{-1}\}.
\end{array}\end{equation*}
\item[(iii)] Time-asymptotically, the approximation rarefaction wave and the inviscid rarefaction wave are equivalent, i.e.,
\begin{equation*}
\lim_{t\rightarrow+\infty}\sup_{x\in\mathbf{R}}| (\bar\r, \bar
u, \bar\t)(x,t)-(\r^r,u^r,\t^r)(\f xt)|=0.
\end{equation*}
\end{itemize}
\end{lemma}

Denote the perturbation around the approximate rarefaction wave $(\bar
\r,\bar u, \bar\t)(x,t)$ by
\begin{equation}\label{perturb}
(\phi,\psi,\omega)(t,x)=(\r-\bar\r, u-\bar u,\t-\bar\t)(x,t)
\end{equation}
where $(\r, u, \t)(x,t)$ is the fluid-dynamical quantities related to the solution $F_1(t,x,v)$ of the VPB equation $\eqref{VPB1}_1$.
By \eqref{F1-f} and \eqref{rare-s}, we obtain the system for the perturbation $(\phi,\psi,\omega)$ in \eqref{perturb} as follows
\begin{equation}\label{sys-h}
\left\{
\begin{array}{ll}
\di \phi_t+\bar\r\psi_{1x}+\bar u_1\phi_x+\bar\r_x\psi_1+\bar u_{1x}\phi=-(\phi\psi_1)_x,\\[3mm]
\di \psi_{1t}+\bar u_1\psi_{1x}+\bar u_{1x}\psi_1 +\f23\omega_x
+\f{2\bar\theta}{3\bar\rho}\phi_x+\f23\r_x(\f\theta\r-\f{\bar\t}{\bar\r})
-\Pi_x\f{n_2}{\r}+\psi_1\psi_{1x}=\f4{3\r}(\mu(\t)u_{1x})_x-\f1\r\int v_1^2\G_x dv,\\[3mm]
\di \psi_{it}+\bar u_1\psi_{ix}+\psi_1\psi_{ix}=
\f1\r(\mu(\t)u_{ix})_x-\f1\r\int v_1v_i\G_xdv,i=2,3,\\[3mm]
\di \omega_{t}+\bar
u_1\omega_x+\t_x\psi_1+\f23u_{1x}\omega+\f23\bar\t\psi_{1x}+\f{\Pi_x}{\r}(\int
v_1F_2dv-u_1n_2)
=\f1\r(\k(\t)\t_{x})_x\\
\di\quad+\f{4}{3\r}\mu(\t)u_{1x}^2+\f1\r\sum_{i=2}^3\mu(\t)u_{ix}^2-\f1\r\int
v_1\f{|v|^2}2\G_xdv+\f1\r\sum_{i=1}^3 u_i\int v_1v_i\G_xdv.
\end{array} \right.
\end{equation}
Set the correction function $\bar{\mb{G}}$ as
\begin{equation}\label{bG}
\bar{\mb{G}}=\f{3}{2\t}\mb{L}_{\mb{M}}^{-1}\Big[\mb{P}_1\big(v_1(\f{|v-u|^2}{2\t}\bar\t_x+v_1\bar
u_{1x})\big)\mb{M}\Big],
\end{equation}
and let $\widetilde{\mb{G}}$ be the rest part related to the microscopic function $\mb{G}$
\begin{equation}\label{tG}
\widetilde{\mb{G}}=\mb{G}-\bar{\mb{G}},
\end{equation}
which satisfies
\begin{equation}\label{GE}
\begin{array}{ll}
\di \widetilde{\mb{G}}_t-\mb{L}_{\mb{M}}\widetilde{\mb{G}}=
-\f{3}{2\t}\Big[\mb{P}_1\big(v_1(\f{|v-u|^2}{2\t}\omega_x+v\cdot\psi_{x})\big)\mb{M}\Big]
-\mb{P}_1(v_1\mb{G}_x)-\mb{P}_1(\Pi_x\partial_{v_1}F_2)+2Q(\mb{G},\mb{G})-\bar{
\mb{G}}_t.
\end{array}
\end{equation}
Notice that in \eqref{tG}, $\bar{\mb{G}}$ is subtracted from  $\mb{G}$ when carrying out the lower order energy estimates because $\di\int|({\bar u}_x,{\bar\theta}_x)|^2dx \sim (1+t)^{-1}$ is not integrable with respect to the time $t$.

\subsection{The Proof of Main Result}
Consider the reformulated system \eqref{sys-h}, \eqref{GE} and \eqref{n21}, \eqref{n2}, \eqref{F2-pc-1}. Since
the local existence of solution to the VPB system can be proved similarly as in \cite{Guo1}, to prove the global
existence on the time interval $[0, T]$ with $T>0$ be any positive time, we only need to close the following a-priori estimates:
\begin{equation}\label{assumption-ap}
\begin{array}{ll}
\di \mathcal{N}(T)\di=\sup_{0\leq t \leq
T}\Bigg\{
\|(\phi,\psi,\omega)(t,\cdot)\|^2_{H^1}+\|(\Pi_x,n_2,n_{2x})\|^2+\sum_{0\leq|\beta|\leq 2}\int\int\f{|\partial^\beta(\wt{\mb{G}},\mb{P}_cF_2)|^2}{\mb{M}_*}dv dx\\
\di \qquad\quad +\sum_{|\alpha^\prime|=1,0\leq|\beta^\prime|\leq1}\int\int\f{|\partial^{\alpha^\prime}\partial^{\beta^\prime}(\mb{G},\mb{P}_cF_2)|^2}{\mb{M}_*}dv dx
+\sum_{|\alpha|=2}\int\int\f{|\partial^\alpha(F_1,F_2)|^2}{\mb{M}_*}dv dx\Bigg\}\leq \chi_{\scriptscriptstyle T} ^2,
\end{array}
\end{equation}
where and in the sequel $\chi_{\scriptscriptstyle T} $ is a small positive constant only depending on the initial data
and wave strengths and independent of the time $T$.  Notice that the difference of energy functionals  $\mathcal{E}(T)$ in \eqref{Et}  and $\mathcal{N}(T)$ in \eqref{assumption-ap} lie in the perturbations around the original inviscid rarefaction wave and the approximate rarefaction wave. By Lemma \ref{appu}, it can be seen easily that
$$
\mathcal{N}(T)\leq C(\mathcal{E}(T)+\d) ~~~{\rm and}~~~\mathcal{E}(T)\leq C(\mathcal{N}(T)+\d),
$$
with some uniform positive constant $C$.
Under the a priori assumption \eqref{assumption-ap}, we can prove that
\begin{equation}\label{Final-E}
\begin{array}{ll}
\di\sup_{0\leq t<\i} \Big[\|(\phi,\psi,\omega)(\cdot,t)\|^2_{H^1}+\|(\Pi_x,n_2,n_{2x})(\cdot,t)\|^2+\sum_{0\leq|\beta|\leq2}\int\int\f{|\partial^\beta(\widetilde{\mb{G}},\mb{P}_c F_2)|^2}{\mb{M}_*}(x,v,t) dxdv\\
\di +\sum_{|\alpha^\prime|=1,0\leq|\beta^\prime|\leq1}\int\int \f{|\partial^{\alpha^\prime}\partial^{\beta^\prime}(\mb{G},\mb{P}_c F_2)|^2}{\mb{M}_*}(x,v,t)dxdv+\sum_{|\alpha|=2}\int\int \f{|\partial^\alpha(F_1,F_2)|^2}{\mb{M}_*}(x,v,t) dxdv\Big] \\
\di +\int_0^\i\|\sqrt{\bar u_{1x}}(\phi,\psi_1,\omega)\|^2d\tau+\sum_{1\leq|\alpha|\leq2}\int_0^\i\|\partial^{\alpha}(\phi,\psi,\omega,n_2)\|^2d\tau+\int_0^\i\|(\Pi_x,n_2)\|^2d\tau\\
\di +\sum_{1\leq|\alpha|\leq2}\int_0^\i\int\int \f{\nu(|v|)|\partial^\alpha(\mb{G},\mb{P}_c F_2)|^2}{\mb{M}_*}dxdvd\tau+\sum_{0\leq|\beta|\leq2}\int_0^\i\int\int\f{\nu(|v|)|\partial^\beta(\widetilde{\mb{G}},\mb{P}_c F_2)|^2}{\mb{M}_*}dxdvd\tau\\
\di +\sum_{|\alpha^\prime|=1,|\beta^\prime|=1}\int_0^\i\int\int \f{\nu(|v|)|\partial^{\alpha^\prime}\partial^{\beta^\prime}(\mb{G},\mb{P}_c F_2)|^2}{\mb{M}_*}dxdvd\tau\leq C(\mathcal{N}(0)^2+\delta^{\f18})\leq C(\mathcal{E}(0)^2+\delta_0^{\f18}).
\end{array}
\end{equation}
The detailed proof of the a priori estimates \eqref{Final-E} will be given in Appendix B.

\

Therefore, from \eqref{Final-E}, we can close the a-priori assumption by choosing suitably small positive constants $\varepsilon_0, \delta_0$  and one has
\begin{equation*}
\begin{array}{ll}\label{aa}
\qquad \di \int_0^{+\infty}\int\int\frac{|(F_1-\mb{M}_{[\bar\r,\bar u,\bar\t]}, F_2)_x|^2}{\mb{M}_*} dvdxd\tau\\[2mm]
\qquad\di\leq \int_0^{+\infty}\int\int\frac{|\big(\mb{M}_x-(\mb{M}_{[\bar\r,\bar u,\bar\t]})_x,(\f{n_2}{\rho}\mb{M})_x\big)|^2}{\mb{M}_*} dvdxd\tau
+\int_0^{+\infty}\int\int\frac{|(\mb{G}_x,(\mb{P}_cF_2)_x)|^2}{\mb{M}_*} dv dxd\tau\\[2mm]
\qquad\di\leq C\int_0^{+\infty}\|(\phi,\psi,\omega,n_2)_x\|^2 d\tau
+C\d_0\int_0^{+\infty}\|(\phi,\psi,\omega)\|^2 d\tau\\[2mm]
\di\qquad\qquad +\int_0^{+\infty}\int\int\frac{|(\mb{G}_x,(\mb{P}_cF_2)_x)|^2}{\mb{M}_*} dv dxd\tau+C\d_0^2\leq C(\mathcal{N}(0)^2+\d_0^{\frac18})\leq C(\mathcal{E}(0)^2+\delta_0^{\f18}).
\end{array}
\end{equation*}
From the Vlasov-Poisson-Boltzmann system \eqref{VPB1}, we  can obtain
\begin{equation*}
\begin{array}{ll}\label{bb}
&\di \int_0^{+\infty}\Big|\frac{d}{dt}\int\int\frac{|(F_1-\mb{M}_{[\bar\r,\bar u,\bar\t]}, F_2)_x|^2}{\mb{M}_*} dvdx\Big|d\tau
\leq C(\mathcal{E}(0)^2+\d_0^{\frac18}).
\end{array}
\end{equation*}
Therefore, one has
\begin{equation*}
\begin{array}{ll}
\di\int_0^{+\infty}\Big(\int\int\frac{|(F_1-\mb{M}_{[\bar\r,\bar u,\bar\t]}, F_2)_x|^2}{\mb{M}_*} dv dx
+\Big|\frac{d}{dt}\int\int\frac{|(F_1-\mb{M}_{[\bar\r,\bar u,\bar\t]}, F_2)_x|^2}{\mb{M}_*} dv dx\Big|
\Big)d\tau<\infty,
\end{array}
\end{equation*}
which implies that
$$
\lim_{t\rightarrow+\infty} \int\int\frac{|(F_1-\mb{M}_{[\bar\r,\bar u,\bar\t]}, F_2)_x|^2}{\mb{M}_*} dv dx=0.
$$
By Sobolev inequality
$$
\begin{array}{ll}
\di \|\int\frac{|(F_1-\mb{M}_{[\bar\r,\bar u,\bar\t]}, F_2)|^2}{\mb{M}_*} dv\|^2_{L^{\infty}_x}\\
\di \leq C\Big(\int\int\frac{|(F_1-\mb{M}_{[\bar\r,\bar u,\bar\t]}, F_2)|^2}{\mb{M}_*} dv dx\Big)\cdot\Big(\int\int\frac{|(F_1-\mb{M}_{[\bar\r,\bar u,\bar\t]}, F_2)_x|^2}{\mb{M}_*} dv dx\Big),
\end{array}
$$
we can prove that
\begin{equation}\label{pr1}
\lim_{t\rightarrow+\infty} \sup_x\int\frac{|(F_1-\mb{M}_{[\bar\r,\bar u,\bar\t]}, F_2)|^2}{\mb{M}_*} dv =0.
\end{equation}
Similarly, one can prove that
$$
\lim_{t\rightarrow+\infty} \|(\Pi_x,n_2)\|=0.
$$
By Lemma \ref{appu}, it holds that
\begin{equation}\label{pr2}
\lim_{t\rightarrow+\infty} \sup_x\int\frac{|\mb{M}_{[\bar\r,\bar u,\bar\t]}-\mb{M}_{[\r^r,u^r,\t^r]}|^2}{\mb{M}_*} dv =0.
\end{equation}
Thus the time-asymptotic convergence of the solutions $(F_1,F_2)$ to the rarefaction wave $\mb{M}_{[\r^r,u^r,\t^r]}$ can be derived directly from \eqref{pr1} and \eqref{pr2}. By \eqref{Final-E} and Lemma \ref{appu}, it holds that
 $$
 \mathcal{E}(t)\leq C(\mathcal{N}(t)+\delta)\leq C(\mathcal{E}(0)+\delta_0^{\f18}),\qquad \forall t\in [0,+\infty),
 $$
which proves \eqref{290}. And the proof of Theorem \ref{thm} is complete.

\section*{Appendix A. A Priori Estimates for Stability of Boltzmann Shock Profiles}
\renewcommand{\theequation}{A.\arabic{equation}}
\setcounter{equation}{0}

By \eqref{F1-f} and \eqref{shock-Euler}, one has
\begin{equation}\label{sys-h-o-s}
\left\{
\begin{array}{l}
\di \phi_{t}+\psi_{1x}=0,\\
\di \psi_{1t}+(\r u^2_1-\tilde\r\tilde u^2_1+p-\tilde p)_x-(\f{\Pi_x ^2}{4})_x
=-\f{4}{3}\Big(\mu(\tilde\t)\tilde u_{1x}-\mu(\t^{S_1}) u^{S_1}_{1x}-\mu(\t^{S_3}) u^{S_3}_{1x}\Big)_x\\
\di\qquad\qquad \qquad  -\int v_1^2\wt{\mb{G}}_xd v- Q_{1x},\\
\di \psi_{it}+(\r u_1u_i)_x=-\int v_1 v_i{\wt{\mb{G}}}_xd v,~~i=2,3,\\
\di
\o_{t}+(\r u_1\t-\tilde\r\tilde u_1\tilde\t+\r u_1\f{|u|^2}{2}-\tilde\r \tilde u_1\f{|\tilde u|^2}{2}
+pu_{1}-\tilde p\tilde u_{1})_x=-\Big(\k(\tilde\t)\tilde\t_x-\k(\t^{S_1})\t^{S_1}_x-
\k(\t^{S_3})\t^{S_3}_x\Big)_x\\
\di\qquad-\f{4}{3}\Big(\mu(\tilde\t)\tilde u_1\tilde u_{1x}-\mu(\t^{S_1})u^{S_1}_1 u^{S_1}_{1x}-\mu(\t^{S_3})u^{S_3}_1 u^{S_3}_{1x}\Big)_x-\f{1}{2}\int v_1| v|^2\wt{\mb{G}}_xd v-Q_{2x}.
\end{array}
\right.
\end{equation}

In fact, by the a priori assumption \eqref{assumption-ap-s}, one also has from the system \eqref{sys-h-o-s} that
\begin{equation}
\|(\widetilde{\Psi}, \widetilde{W})\|_{H^2}^2, \|(\Phi,\Psi,W,\wt\Psi,\wt W)\|_{L^\i_x}^2,\|(\phi,\psi,\omega,\Pi_x ,n_2)\|_{L^\i_x}^2\leq C(\chi_{\scriptscriptstyle T} +\d_0)^2,
\end{equation}
and
\begin{equation*}
\|(\phi_{t},\psi_{t},\omega_{t})\|^2\leq C(\chi_{\scriptscriptstyle T} +\d_0)^2,
\label{(4.19)}
\end{equation*}
hence, one has
\begin{equation*}
\|(\r_{t},u_{t},\t_{t})\|^2\di
\leq C\|(\r_{t},m_{t},E_{t})\|^2 \leq C\|(\phi_{t},\psi_{t},\omega_{t},\Pi_x ,\Pi_{xt})\|^2+ C \|(\tilde\r_{t},\tilde m_{t},\tilde E_{t})\|^2\leq C(\chi_{\scriptscriptstyle T} +\d_0)^2.
\end{equation*}
For $|\alpha|=2$, it follows from \eqref{macro}  and \eqref{macro-n2} that
\begin{equation}
\begin{array}{ll}
\di \|\partial^\a\left(\rho,\rho u,\rho(\t+\f{|u|^2}{2}),n_2\right)\|^2 \leq
C\int\int\f{|\partial^\a (F_1,F_2)|^2}{\mb{M}_*}dv dx\\
\di \leq C\int\int\f{|\partial^\a (\wt F_1,F_2)|^2}{\mb{M}_*}dv dx+C\int\int\f{|\partial^\a (F_{\a_1,\a_3}^S)|^2}{\mb{M}_*}dv dx\leq C(\chi_{\scriptscriptstyle T} +\delta_0)^2,
\end{array}
\label{(4.8)-s}
\end{equation}
and
\begin{equation}
\begin{array}{ll}
\di\|\partial^\a(\r,u,\t)\|^2&\di\leq C\|\partial^\a\left(\rho,\rho
u,\rho(\t+\f{|u|^2}{2})\right)\|^2
+C\sum_{|\a^\prime|=1}\int|\partial^{\a^\prime}\left(\rho,\rho
u,\rho(\t+\f{|u|^2}{2})\right)|^4dx\\
 &\di\leq C(\chi_{\scriptscriptstyle T} +\delta_0)^2.
\end{array}
\label{3.25-s}
\end{equation}
Therefore, for $|\a|=2$, we have
\begin{equation}
\|\partial^\a(\phi,\psi,\omega,n_2)\|^2\leq
C(\chi_{\scriptscriptstyle T} +\delta_0)^2, \label{3.26-s}
\end{equation}
By \eqref{n21} and a priori assumption \eqref{assumption-ap-s}, it holds that
\begin{equation}\label{n2t-ap-s}
\|n_{2t}\|^2\leq C\Big[\|n_{2x}\|^2+\int|n_2|^2|(\r_x,u_x,\t_x)|^2 dx+\int\int\f{|(\mb{P}_c F_2)_x|^2}{\mb{M}_*} dxdv\Big]\leq C(\chi_{\scriptscriptstyle T} +\delta_0)^2,
\end{equation}
By \eqref{assumption-ap-s}, \eqref{3.26-s} and \eqref{n2t-ap-s}, for $|\alpha^\prime|=1$, it holds that
\begin{equation*}
\begin{array}{ll}
\|\partial^{\alpha^\prime}(\phi,\psi,\omega,n_2)\|^2_{L^{\i}}\leq
C(\chi_{\scriptscriptstyle T} +\delta_0)^2.
\end{array}
\end{equation*}
By \eqref{n21}, one has
\begin{equation*}
\f12\Pi_{xt}+\int v_1 F_2dv=0,
\end{equation*}
that is,
\begin{equation}\label{phi-xt-e-s}
\f12\Pi_{xt}+u_1n_2+\int v_1 \mb{P}_cF_2dv=0.
\end{equation}
 Therefore, one has
\begin{equation}\label{phi-xt-ap-s}
\|\Pi_{xt}\|^2\leq C\Big[\|n_2\|^2+\int\int\f{|\mb{P}_c F_2|^2}{\mb{M}_*} dxdv\Big]\leq C(\chi_{\scriptscriptstyle T} +\delta_0)^2.
\end{equation}
By \eqref{phi-xt-e-s}, it holds that
\begin{equation}\label{phi-xtt-e-s}
\f12\Pi_{xtt}+(u_1n_2)_t+\int v_1 (\mb{P}_c F_2)_{t}dv=0.
\end{equation}
Thus, one has
\begin{equation}\label{phixtt-ap-s}
\|\Pi_{xtt}\|^2\leq C\Big[\|n_{2t}\|^2+\int|n_2|^2|u_t|^2 dx+\int\int\f{|(\mb{P}_c F_2)_t|^2}{\mb{M}_*} dxdv\Big]\leq C(\chi_{\scriptscriptstyle T} +\delta_0)^2.
\end{equation}
By \eqref{phi-xt-ap-s}, \eqref{n2t-ap-s}, \eqref{phixtt-ap-s} and \eqref{3.26-s}, it holds that
\begin{equation*}\label{phi-xt-in}
\|(\Pi_{xt}, \Pi_{xtt})\|_{L^\infty}^2\leq C\|(\Pi_{xt}, \Pi_{xtt})\|\|(n_{2t}, n_{2tt})\|\leq  C(\chi_{\scriptscriptstyle T} +\delta_0)^2.
\end{equation*}
Moreover, it holds that
\begin{equation*}
\begin{array}{ll}
 \di\|\int\f{|(\wt{\mb{G}},\mb{P}_c F_2)|^2}{\mb{M}_*}dv \|_{L^\infty_x} \leq
C\left(\int\int\f{|(\wt{\mb{G}},\mb{P}_c F_2)|^2}{\mb{M}_*}dv
dx\right)^{\f{1}{2}}\cdot\left(\int\int\f{|(\wt{\mb{G}},\mb{P}_c F_2)_x|^2}{\mb{M}_*}dv dx\right)^{\f{1}{2}}\\
\di \qquad\qquad\qquad\qquad \leq
C(\chi_{\scriptscriptstyle T} +\delta_0)^2.
\end{array}
\label{G-infty}
\end{equation*}
Furthermore,  for $|\a^\prime|=1$, it holds that
\begin{equation}
\begin{array}{ll}
 \di\|\int\f{|\partial^{\a^\prime}
(\wt{\mb{G}}, \mb{P}_c F_2)|^2}{\mb{M}_*}d v\|_{L_x^{\i}}\leq
C\left(\int\int\f{|\partial^{\a^\prime} (\wt{\mb{G}}, \mb{P}_c F_2)|^2}{\mb{M}_*}dv
dx\right)^{\f{1}{2}}\cdot\left(\int\int\f{|\partial^{\a^\prime}
(\wt{\mb{G}}, \mb{P}_c F_2)
_x|^2}{\mb{M}_*}dv dx\right)^{\f{1}{2}}\\
\di\qquad\qquad \leq C\chi_{\scriptscriptstyle T} (\chi_{\scriptscriptstyle T} +\d_0)\le C(\chi_{\scriptscriptstyle T} +\delta_0)^2.
\end{array}
\label{(4.17+)}
\end{equation}
Finally, by noticing the facts that $F_1=\mb{M}+\mb{G}$ and
$F_2=\f{n_2}{\r}\mb{M}+\mb{P}_c F_2$, it holds that
\begin{equation*}
\begin{array}{ll}
\di\int\int\f{|\partial^{\alpha}(\wt{\mb{G}},\mb{P}_cF_2)|^2}{\mb{M}_*}dv dx\di
\le C\int\int\f{|\partial^\alpha (\wt F_1,F_2)|^2}{\mb{M}_*}dv dx
\\[3mm]
\di \qquad+C\int\int\f{|\partial^\alpha(\mb{M}-\mb{M}^{S_1}-\mb{M}^{S_3})|^2+|\partial^\a\big(\f{n_2}{\rho}\mb{M}\big)|^2}{\mb{M}_*}dv dx\di
\leq C(\chi_{\scriptscriptstyle T} +\d_0)^2.
\end{array}
\end{equation*}

{\bf Appendix A.1.  Proof of Proposition 3.1: Lower order estimates}

\renewcommand{\theequation}{A.1.\arabic{equation}}
\setcounter{equation}{0}

\

The proof of the lower order estimates in Proposition \ref{Prop3.1} includes the following five steps.

\

\underline{Step 1. Estimation on
$\|(\Phi,\wt\Psi,\wt W)(t,\cdot)\|^2$.}

\

Multiplying
$\eqref{P-sys}_1$ by $\di\f{2\Phi}{3\tilde\r}$, $\eqref{P-sys}_2$ by
$\di\f{\wt\Psi_1}{\tilde\t}$, $\eqref{P-sys}_3$ by $\wt\Psi_i$,
$\eqref{P-sys}_4$ by $\di\f{\wt W}{{\tilde\t}^2}$, respectively, and
collecting the resulted equations together, we can get
\begin{equation}
\begin{array}{ll}
\di\left(\f{\Phi^2}{3\tilde\r}+\f{\tilde\r\wt\Psi_1^2}{2\tilde\t}+\sum_{i=2}^3\f{\tilde\r\wt\Psi_i^2}{2}+\f{\tilde\r\wt
W^2}{2\tilde\t^2}\right)_t+(\cdots)_x
-\tilde\r\tilde u_{1x}\left(\f{\wt\Psi_1^2}{3\tilde\t}+\sum_{i=2}^3\wt\Psi_i^2+\f{\wt W^2}{\tilde\t^2}\right)+\f{4\mu(\tilde\t)}{3\tilde\t}\wt\Psi^2_{1x}+\sum_{i=2}^3\mu(\tilde\t)\wt\Psi_{ix}^2\\
\di +\f{\k(\tilde\t)}{\tilde\t^2}\wt
W_x^2=-\left(\f{4\mu(\tilde\t)}{3\tilde\t}\right)_x\wt\Psi_1\wt\Psi_{1x}-\sum_{i=2}^3(\mu(\tilde\t))_x\wt\Psi_i\wt\Psi_{ix}-\left(\f{\k(\tilde\t)}{\tilde\t^2}\right)_x\wt W\wt
W_x-\left(\f{\wt\Psi_1^2}{2\tilde\t^2}+\f{\wt
W^2}{\tilde\t^3}\right)(\tilde\r\tilde\t_t+\tilde\r\tilde u_1\tilde\t_x)\\
\di+K_1+\f{\wt\Psi_1}{\tilde\t}(J_1+N_1-Q_1)+\sum_{i=2}^3\wt\Psi_i(J_i+N_i) +\f{\wt
W}{\tilde\t^2}(J_4+N_4-Q_2+\tilde u_1Q_1\big)
\end{array}
\label{le1-s}
\end{equation}
with
\begin{equation*}
\begin{array}{ll}
K_1=&\di -\f{\wt\Psi_1}{\tilde\t}\int v_1^2(\G-\G^{S_1}-\G^{S_3})dv
-\sum_{i=2}^3\wt\Psi_i\int v_1v_i(\G-\G^{S_1}-\G^{S_3})dv \\
&\di-\f{\wt
W}{\tilde\t^2}\left[\int v_1\f{|v|^2}{2}(\G-\G^{S_1}-\G^{S_3})dv-\tilde
u_1\int v_1^2(\G-\G^{S_1}-\G^{S_3})dv\right].
\end{array}
\label{K1}
\end{equation*}
Here and in the sequel the notation $(\cdots)_x$ represents the term in the conservative form so that it
vanishes after integration. By \eqref{jinsi}, one has
\begin{equation}
\begin{array}{ll}
\tilde\r\tilde\t_t+\tilde\r\tilde u_1\tilde\t_x=&\di-\f23\tilde\r\tilde\t\tilde
u_{1x}+\f43\mu(\tilde\t)\tilde u_{1x}^2+(\k(\tilde\t)\tilde\t_x)_x
-\int v_1\f{|v|^2}{2}(\G^{S_1}_x+\G^{S_3}_x)dv  \\
&\di+\tilde u_1\int v_1^2(\G^{S_1}_x+\G^{S_3}_x)dv
+Q_{2x}-\tilde u_1Q_{1x}. \label{bar-theta}
\end{array}
\end{equation}
Substituting \eqref{bar-theta} into \eqref{le1-s}, carring out the similar estimates as in \cite{Wang-Wang}  and choosing $\delta_0$   and
$\chi_{\scriptscriptstyle T} $ suitably small yield that
\begin{equation} \label{le-s}
\begin{array}{ll}
\di \|(\Phi,\wt\Psi,\wt W)(t,\cdot)\|^2+\int_0^t\left[\|\sqrt{|u^{S_1}_{1x}|
+|u^{S_3}_{1x}|}(\wt\Psi,\wt W)\|^2+\|(\wt\Psi_x,\wt W_x)\|^2\right]d\tau\\
\di \leq C\int\int\f{|\wt{\mb{G}}|^2}{\mb{M}_*}(t,x,v)d v
dx
+ C(\chi_{\scriptscriptstyle T} +\d_0)\di\int_0^t\|(\phi,\psi,\o)\|^2d\tau+\s \int_0^t\|(\wt \Psi_{\tau},\wt W_{\tau})\|^2d\tau
\\
 \di ~ +C\d_0\int_0^t \|\sqrt{|u^{S_1}_{1x}|+|u^{S_3}_{1x}|}~\Phi\|^2 d\tau
+ C_\s\int_0^t\int\int\f{\nu(| v|)}{\mb{M}_*}|\wt{\mb{G}}|^2d v
dxd\tau+C(\chi_{\scriptscriptstyle T} +\delta_0)\int_0^t \|(\Pi_x,n_2,\Pi_{x\tau})\|^2d\tau\\
\di ~+C(\chi_{\scriptscriptstyle T} +\d_0)\di\sum_{|\a^\prime|=1}\int_0^t\|\partial^{\a^\prime}(\phi,\psi,\o)\|^2d\tau+C\chi_{\scriptscriptstyle T} \int_0^t\int\int \f{\nu(|v|)|\partial_{v_1}(\mb{P}_c F_2)|^2}{\mb{M}_*} dxdvd\tau+C(\mathcal{E}(0)^2+\d_0^{\f12}),
\end{array}
\end{equation}
where and in the sequel section the global Maxellian $\mb{M}_*=\mb{M}_{[\r_*,u_*,\t_*]}$ is chosen such that
\begin{equation}\label{ma1-s}
\rho_*>0,\qquad \f12\t(t,x)<\t_*<\t(t,x),
\end{equation}
and
\begin{equation}\label{ma2-s}
|\rho(x,t)-\r_*|+|u(x,t)-u_*|+|\t(x,t)-\t_*|<\eta_0
\end{equation}
with $\eta_0$ being the small positive constant in Lemma \ref{Lemma 4.2}. In fact, if the wave strength $\delta+|\alpha_2|$ is suitably small, then it is holds that
$$
\f12\sup_{x,t}\tilde\t(x,t)<\inf_{x,t}\tilde\t(x,t).
$$
Therefore, we can choose the global Mawellian $\mb{M}_*=\mb{M}_{[\rho_*,u_*,\t_*]}$ satisfying \eqref{ma1-s} and \eqref{ma2-s} provided that the solution $(\r,u,\t)(x,t)$ is near the ansatz $(\tilde\r,\tilde u,\tilde\t)(x,t)$ as in a priori assumption \eqref{assumption-ap-s}.

\

\underline{Step 2. Estimation on
$\int_0^t\int(|u_{1x}^{S_1}|+|u_{1x}^{S_3}|)\Phi^2 dxd\tau$.}

\

 In order to estimate
$\int_0^t\int(|u_{1x}^{S_1}|+|u_{1x}^{S_3}|)\Phi^2 dxd\tau$, we
borrow the ideas of  the vertical estimates in Goodman
\cite{Goodman} to carry out the following characteristic weight
estimates.  First, we
diagonalize the system \eqref{P-sys}. Let $
V=(\Phi,\wt{\Psi}_1,\wt{W})^t, $ then
\begin{equation}
V_t+A_1V_x+A_2V=A_3V_{xx}+A_4,\label{V-sys}
\end{equation}
where
\begin{equation*}
A_1=\left(
\begin{array}{ccc}
\tilde u_1&\tilde \r&0\\
\f{2\tilde\t}{3\tilde\r}&\tilde u_1&\f23\\
0&\f23\tilde\t&\tilde u_1
\end{array}
\right), A_2=\left(
\begin{array}{ccc}
\tilde u_{1x}&\tilde \r_x&0\\
-\f{2\tilde\t\tilde\r_x}{3\tilde\r^2}&-\f{\tilde u_{1x}}{3}&\f{2\tilde\r_x}{3\tilde\r}\\
0&-\f23\tilde\t_x&-\tilde u_{1x}
\end{array}
\right), A_3=\left(
\begin{array}{ccc}
0&0&0\\
0&\f{4\mu(\tilde\t)}{3\tilde\r}&0\\
0&0&\f{\k(\tilde\t)}{\tilde\r}
\end{array}
\right),
\end{equation*}
and
\begin{equation*}
A_4=\f{1}{\tilde\r}\left(\begin{array}{c}
\di 0\\
\di -\int v_1^2(\G-\G^{S_1}-\G^{S_3})dv+J_1+N_1-Q_1\\
\di \int\left(\tilde u_1v_1^2-\f{v_1|v|^2}{2}\right)(\G-\G^{S_1}-\G^{S_3})dv+J_4+N_4-(Q_2-\tilde u_1Q_1)
\end{array}
\right).
\end{equation*}
We can compute three eigenvalues of the matrix $A_1$ in the
system \eqref{V-sys}
\begin{equation*}
\tilde \l_1=\tilde u_1-\f{\sqrt{10\tilde\t}}{3},\quad \tilde \l_2=\tilde u_1,\quad
\tilde \l_3=\tilde u_1+\f{\sqrt{10\tilde\t}}{3},\label{eigenvalue}
\end{equation*}
with corresponding left and right eigenvectors given by
\begin{equation*}
l_1=\left(\tilde\t,-\f{\sqrt{10\tilde\t}}{2}\tilde\r,\tilde\r\right),\quad
l_2=\left(\tilde\t,0,-\f32\tilde\r\right),\quad
l_3=\left(\tilde\t,\f{\sqrt{10\tilde\t}}{2}\tilde\r,\tilde\r\right),
\label{l-eigenvector}
\end{equation*}
and
\begin{equation*}
r_1=\f{3}{10\tilde\r\tilde\t}\left(\tilde\r,-\f{\sqrt{10\tilde\t}}{3},\f23\tilde\t\right)^t,
r_2=\f{2}{5\tilde\r\tilde\t}(\tilde\r,0,-\tilde\t)^t,
r_3=\f{3}{10\tilde\r\tilde\t}\left(\tilde\r,\f{\sqrt{10\tilde\t}}{3},\f23\tilde\t\right)^t,
\end{equation*}
respectively. If we denote that the matrix composed with the left
and right eigenvalues by $L=(l_1,l_2,l_3)^t, R=(r_1,r_2,r_3),$
then we have
\begin{equation*}\label{diag}
LR=Id.,\qquad LA_1R=\Lambda:={\rm diag}(\tilde\l_1,\tilde\l_2,\tilde\l_3),
\end{equation*}
with $Id.$ being the $3\times3$ identity matrix.
Denote that $Z=LV$  with $Z=(Z_1,Z_2,Z_3)^t,$
then we have $ V=RZ.$
Multiplying the system \eqref{V-sys} by the matrix $L$ on the left,
we can obtain the diagonalized system for $Z$
\begin{equation}
Z_t +\Lambda Z_x-LA_3RZ_{xx}=-L(R_t+A_1R_x-A_3R_{xx})Z-LA_2R Z+2LA_3R_xZ_x+LA_4.
\label{Z-sys}
\end{equation}
Introduce the weight function
$$
\a(t,x)=\f{\r^{S_1}(t,x)}{\r_\#},\quad \b(t,x)=\f{\r^{S_3}(t,x)}{\r_\#}.
$$
From the properties of the shock profile to Boltzmann equation, we
have
$$
\l^{S_i}_{ix}<0,\qquad {\rm and}~~\r^{S_i}_x<0 \quad(i=1,3).
$$
Thus it holds that
$$
\a,\b~<1 \quad {\rm and} ~~|\a-1|,~ |\b-1|\leq  C\d \ll 1,~~{\rm if
}~~\d\ll1.
$$
Taking inner product by multiplying $\bar Z:=(Z_1,
\a^{N}Z_2, \a^{N}Z_3)^t$ with $N$ being a large positive constant to be
determined and noting that
 for $i=2,3,$
\begin{equation*}
\begin{array}{ll}
\di\a_t+\tilde\l_{i}\a_x&\di =-s_1\a_x+\tilde\l_i\a_x
 =(\l_i^{S_1}-s_1)\a_x+(\tilde\l_i-\l_i^{S_1})\a_x,
\end{array}
\label{eta1}
\end{equation*}
 we have
\begin{equation*}\label{y12}
\begin{array}{ll}
\di \left[\f{Z_1^2+\a^{N}(Z_2^2+Z_3^2)}{2}\right]_t
+(\cdots)_x -\l^{S_1}_{1x}\f{Z_1^2}{2}-\a^{N}\sum_{i=2}^3\f{\l^{S_3}_{ix}Z_i^2}{2}-N\a^{N-1}\a_x\sum_{i=2}^3\f{(\l^{S_1}_i-s_1)Z_i^2}{2}\\
\di -\bar Z\cdot
LA_3RZ_{xx}=-\bar Z\cdot L(R_t+A_1R_x-A_3R_{xx})Z-\bar Z\cdot L A_2R Z+2\bar Z\cdot LA_3R_xZ_x+\bar Z\cdot LA_4\\[3mm]
\di+(\tilde\l_{1x}-\l^{S_1}_{1x})\f{Z^2_1}{2}+\a^N\sum_{i=2}^3\f{(\tilde\l_{ix}-\l^{S_3}_{ix})Z_i^2}{2}
+N\a^{N-1}\a_x\sum_{i=2}^3\f{(\tilde\l_i-\l_i^{S_1})Z_i^2}{2}.
\end{array}
\end{equation*}
Note that
\begin{equation*}
\begin{array}{ll}
\di-\bar Z\cdot LA_3R Z_{xx} =-(\bar Z\cdot LA_3R Z_{x})_x+\bar Z_x\cdot LA_3R Z_x+\bar Z\cdot (LA_3R)_x Z_{x}\\
\quad\di =-(\bar Z\cdot LA_3R Z_{x})_x+Z_x\cdot LA_3R Z_{x}+(\bar
Z-Z)_x\cdot LA_3R Z_{x}+\bar Z\cdot (LA_3R)_x Z_{x}.
\end{array}
\end{equation*}
We can directly compute that the matrix $LA_3R$ is non-negatively
definite, and so
$$
Z_x\cdot LA_3R Z_{x}\geq 0.
$$
On the other hand, it holds that
\begin{equation*}
(\bar Z-Z)_x=(0,(\a^{N}-1)Z_2,(\a^{N}-1)Z_3)_x\ ^t=(\a^{N}-1)(0,Z_{2x},Z_{3x})^t+N\a^{N-1}\a_x(0,Z_2,Z_3)^t.
\end{equation*}
By the Lax entropy condition to 1-shock, we have
$$
\l_2^{S_1}-s_1>\l_2^{S_1}-\l_{1-}=u_1^{S_1}-(u_{1-}-\frac{\sqrt{10\t_-}}{3})\geq\frac{\sqrt{10\t_-}}{3}-C\d>\frac{\sqrt{10\t_-}}{6}
$$
and
$$
\l_3^{S_1}-s_1>\l_2^{S_1}-\l_{1-}>\frac{\sqrt{10\t_-}}{6}~~{\rm if
}~~\d\ll1.
$$
Therefore, it holds that
\begin{equation*}
\begin{array}{ll}\label{posi1}
\di |(\bar Z-Z)_x\cdot LA_3R Z_{x}|
\leq|(\a^{N}-1)(0,Z_{2x},Z_{3x})^t\cdot LA_3R Z_{x}|+|N\a^{N-1}\a_x(0,Z_2,Z_3)^t\cdot LA_3R Z_{x}|\\
\di\qquad \leq C\d|Z_x|^2+N\a^{N-1}|\a_x|\sum_{i=2}^3|Z_i||Z_x|\leq\f{N\a^{N-1}|\a_x|}{4}\sum_{i=2}^3\f{(\l_i^{S_1}-s_1)Z_i^2}{2}+C\sqrt{\d_0}|Z_x|^2,
\end{array}
\end{equation*}
if we choose $N=\frac{1}{\sqrt{\d_0}}$ with $\d_0\ll 1$.
Then one has
\begin{equation}
\begin{array}{ll}\label{posi2}
\di |\bar Z\cdot (LA_3R)_xZ_{x}|\leq C(|\l^{S_1}_{1x}|+|\l^{S_3}_{1x}|+|\T_x|+q)|Z||Z_x|\\[2mm]
\di\qquad \leq C\sqrt{\d_0}|Z_x|^2+C\sqrt{\d_0}(|\l^{S_1}_{1x}|+|\l^{S_3}_{1x}|)|Z|^2+q|Z|^2,
\end{array}
\end{equation}
where $q\in Q$ defined in \eqref{Q}.
Similar to \cite{Wang-Wang},  we can get
\begin{equation}
\begin{array}{ll}
\di \|Z(t,\cdot)\|^2+\int_0^t\int\left[|\l_{1x}^{S_1}|Z_1^2+\sum_{i=2}^3|\l_{ix}^{S_3}|Z_i^2
+N|\a_x|\sum_{i=2}^3Z_i^2\right]dxd\tau\\[3mm]
\di\le C\int\int\f{|\wt{\mb{G}}|^2}{\mb{M}_*}(t,x,v)d v
dx+C\sqrt{\d_0}\int_0^t\|(Z_{\tau},Z_x)\|^2 d\tau +C(\d_0+\chi_{\scriptscriptstyle T} )\int_0^t\|(\phi,\psi,\o)\|^2d\tau\\[3mm]
\di+C\int_0^t\int|\l_{3x}^{S_3}|(Z_1^2+Z_2^2)dxd\tau
+C\int_0^t\int\int\f{\nu(| v|)}{\mb{M}_*}|\wt{\mb{G}}|^2d v
dxd\tau\\
\di
+C\chi_{\scriptscriptstyle T} \sum_{|\a\prime|=1}\int_0^t\|\partial^{\a\prime}(\phi,\psi,\o)\|^2d\tau+C(\chi_{\scriptscriptstyle T} +\delta_0)\int_0^t \|(\Pi_x,n_2,\Pi_{x\tau})\|^2d\tau\\[3mm]
\di+C\chi_{\scriptscriptstyle T} \int_0^t\int\int \f{\nu(|v|)|\partial_{v_1}(\mb{P}_c F_2)|^2}{\mb{M}_*} dxdvd\tau +C\int_0^t\|\sqrt{|\T_x|}Z\|^2 d\tau
+C(\mathcal{E}(0)^2+\d_0^{\f12}).
\end{array}\label{z1}
\end{equation}
Taking inner product by multiplying $\wt Z:=(\b^{-N}Z_1,
\b^{-N}Z_2, Z_3)^t$ with $N=\frac{1}{\sqrt{\d_0}}$ as before, similar to \eqref{z1}, we can get
\begin{equation}
\begin{array}{ll}
\di \|Z(t,\cdot)\|^2+\int_0^t\int\left[|\l_{3x}^{S_3}|Z_3^2+\sum_{i=1}^2|\l_{ix}^{S_1}|Z_i^2
+N|\b_x|\sum_{i=2}^3Z_i^2\right]dxd\tau\\[3mm]
\di\leq C\int\int\f{|\wt{\mb{G}}|^2}{\mb{M}_*}(t,x,v)d v
dx+C\sqrt{\d_0}\int_0^t\|(Z_{\tau},Z_x)\|^2 d\tau+C(\d_0+\chi_{\scriptscriptstyle T} )\int_0^t\|(\phi,\psi,\o)\|^2d\tau\\[3mm]
\di+C\int_0^t\int|\l_{1x}^{S_1}|(Z_2^2+Z_3^2)dxd\tau
+C \int_0^t\int\int\f{\nu(| v|)}{\mb{M}_*}|\wt{\mb{G}}|^2d v
dxd\tau\\
\di
+C\chi_{\scriptscriptstyle T} \sum_{|\a\prime|=1}\int_0^t\|\partial^{\a\prime}(\phi,\psi,\o)\|^2d\tau+ C\int_0^t\|\sqrt{|\T_x|}Z\|^2 d\tau
+C(\mathcal{E}(0)^2+\d_0^{\f12})\\[3mm]
\di+C\chi_{\scriptscriptstyle T} \int_0^t\int\int \f{\nu(|v|)|\partial_{v_1}(\mb{P}_c F_2)|^2}{\mb{M}_*} dxdvd\tau+C(\chi_{\scriptscriptstyle T} +\delta_0)\int_0^t \|(\Pi_x,n_2,\Pi_{x\tau})\|^2d\tau.
\end{array}\label{z3}
\end{equation}
Combining \eqref{z1} and \eqref{z3} and choosing $\d_0$ sufficiently small, it holds that
\begin{equation}
\begin{array}{ll}
\di \|Z(t,\cdot)\|^2+\int_0^t\|\sqrt{|\l_{1x}^{S_1}|+|\l_{3x}^{S_3}|}Z\|^2 d\tau
\le C(\mathcal{E}(0)^2+\d_0^{\f12})+ C \int\int\f{|\wt{\mb{G}}|^2}{\mb{M}_*}(t,x,v)d v
dx\\[3mm]
\di+C\int_0^t\|\sqrt{|\T_x|}Z\|^2 d\tau
+C\int_0^t\int\int\f{\nu(| v|)}{\mb{M}_*}|\wt{\mb{G}}|^2d v
dxd\tau
+C(\d_0+\chi_{\scriptscriptstyle T} )\int_0^t\|(\phi,\psi,\o)\|^2d\tau\\
\di+C\chi_{\scriptscriptstyle T} \di\sum_{|\a\prime|=1}\int_0^t\|\partial^{\a\prime}(\phi,\psi,\o)\|^2d\tau
+C\sqrt{\d_0}\int_0^t\|(Z_{\tau},Z_x)\|^2 d\tau \\
\di +C\chi_{\scriptscriptstyle T} \int_0^t\int\int \f{\nu(|v|)|\partial_{v_1}(\mb{P}_c F_2)|^2}{\mb{M}_*} dxdvd\tau+C(\chi_{\scriptscriptstyle T}+\delta_0)\int_0^t \|(\Pi_x,n_2,\Pi_{x\tau})\|^2d\tau.
\end{array}\label{z1z3}
\end{equation}

\

\underline{Step 3. Estimation on
$\int_0^t \|\sqrt{|\T_x|}Z\|^2 d\tau$.}

In this step, we estimate $\int_0^t \|\sqrt{|\T_x|}Z\|^2 d\tau$ on the right hand side of \eqref{z1z3}.
Note that the linear diffusion wave $\T(t,x)$ in \eqref{Theta} and the coupled diffusion wave $\T_x(t,x)$ propagate along the constant speed $u_{1\#}$.  Therefore, somehow we can view these diffusion waves as the viscous contact wave in the linearly degenerate field as constructed in \cite{Huang-Matsumura-Xin}, \cite{Huang-Xin-Yang}. In fact, the viscous contact waves is exactly constructed through the self-similar solution to the diffusion equation. Thus we can borrow some ideas of the weighted estimates for the viscous contact wave in \cite{Huang-Li-matsumura} to estimate the term $\int_0^t \|\sqrt{|\T_x|}Z\|^2 d\tau$. Similar to the viscous contact wave in the second characteristic field, the diffusion waves here has the extra dissipation on the first and third transverse characteristic fields. By the delicate weighted estimates, we can get such estimate as
\begin{equation}\label{note2}
|\a_2|\int_0^t\int\big[h(Z_1^2+Z_3^2)+h^2Z_2^2\big] dxd\tau
\end{equation}
 with $h\sim (1+t)^{-\f12}e^{-\f{c(x-u_{1\#}t)^2}{1+t}}$ defined in \eqref{eta11}. Note that \eqref{note2} means that the diffusion wave in the second characteristic field has some extra dissipative effects on the first and third transverse characteristic fields compared with the second diffusion wave field. By the following inequality
\begin{equation}\label{note1}
\begin{array}{ll}
\di
|\T_x||Z|^2\leq C|\a_2|(1+t)^{-1}e^{-\frac{(x-u_{1\#}t)^2}{8a(1+t)}}\big(Z_1^2+Z_2^2+Z_3^2\big)\leq C|\a_2|\big[h(Z_1^2+Z_3^2)+h^2Z_2^2\big],
\end{array}
\end{equation}
we can get the desired estimate for $\int_0^t \|\sqrt{|\T_x|}Z\|^2 d\tau$ from \eqref{note2}. In the following, we focus on the proof of \eqref{note2}.

We first set the matrices $(c_{ij})_{n\times n}$ and $(b_{ij})_{n\times n}$ as
$$
LA_3R:=(c_{ij})_{n\times n},\qquad
L(R_t+A_1R_x-A_3R_{xx}+A_2R):=(b_{ij})_{n\times n}.
$$
Then from \eqref{Z-sys}, one has the equation for $Z_1$:
\begin{equation}\label{z1t}
Z_{1t}+\tilde\l_1Z_{1x}=\sum_{j=1}^3 c_{1j}Z_{jxx}-\sum_{j=1}^3b_{1j}Z_{j}+(2LA_3R_xZ_x+LA_4)_1,
\end{equation}
here $(\cdot)_i~(i=1,2,3)$ denotes the $i-$th component of the vector $(\cdot)$. Set
\begin{equation}\label{eta11}
h=\f{1}{\sqrt{16\pi a(1+t)}}\exp\left(-\f{(x-u_{1\#}t)^2}{16a(1+t)}\right)
~~{\rm and}~~\eta_1=\exp\left(\int_{-\i}^{x}h(y,t) dy\right),\quad
\di \end{equation}
with $a=\frac{3\kappa(\t_\#)}{5\r_\#}>0$. Then $h$ satisfies
$$
h_t+u_{1\#}h_x=ah_{xx}.
$$
Obviously, it holds that
$
1\leq \eta_1\leq e
$
and
\begin{equation*}\label{eta1t}
\di \eta_{1t}=\eta_1\int_{-\infty}^{x}h_t(y,t) dy=\eta_1(ah_x-u_{1\#}h),\quad \eta_{1x}=\eta_1h.
\end{equation*}
Multiplying \eqref{z1t} by $\eta_1Z_1$, we can get
\begin{equation}\label{024}
\begin{array}{ll}
\di\left(\eta_1\f{Z_1^2}{2}\right)_t-(\eta_{1t}+\tilde\l_1\eta_{1x})\f{Z_1^2}{2}
-\tilde\l_{1x}\eta_1\f{Z_1^2}{2}\\[3mm]
\di =-\sum_{j=1}^3Z_{jx}(c_{1j}\eta_1Z_1)_x+\left[-\sum_{j=1}^3b_{1j}Z_{j}+(2LA_3R_xZ_x+LA_4)_1\right]\eta_1Z_1+(\cdots)_x.
\end{array}
\end{equation}
Integrating \eqref{024} with respect to  $x,t$ and similar to \cite{Wang-Wang}, one has
\begin{equation}
\begin{array}{ll}
\di\int Z_1^2 dx+\int_0^t\int hZ_1^2dxd\tau\le
C\int\int\f{|\wt{\mb{G}}|^2}{\mb{M}_*}(t,x,v)d v
dx +C\int_0^t\left[\|\sqrt{|\l_{1x}^{S_1}|+|\l_{3x}^{S_3}|}Z\|^2+\|Z_x\|^2\right] d\tau\\[3mm]
\di\qquad
+\s\int_0^t\|Z_{1\tau}\|^2 d\tau+C_\s\int_0^t\int\int\f{\nu(| v|)}{\mb{M}_*}|\wt{\mb{G}}|^2d v
dxd\tau+C(\d_0+\chi_{\scriptscriptstyle T} )\int_0^t\|(\phi,\psi,\o)\|^2d\tau\\[3mm]
\di\qquad
+C\chi_{\scriptscriptstyle T} \di\sum_{|\a\prime|=1}\int_0^t\|\partial^{\a\prime}(\phi,\psi,\o)\|^2d\tau
+C\int_0^t\int |\T_x|(Z_2^2+Z_3^2) dxd\tau+C(\mathcal{E}(0)^2+\d_0^{\frac12})\\
\di\qquad+C\chi_{\scriptscriptstyle T} \int_0^t\int\int \f{\nu(|v|)|\partial_{v_1}(\mb{P}_c F_2)|^2}{\mb{M}_*} dxdvd\tau+C(\chi_{\scriptscriptstyle T} +\delta_0)\int_0^t \|(\Pi_x,n_2,\Pi_{x\tau})\|^2d\tau,
\end{array}\label{intz1}
\end{equation}
with sufficiently small positive constant $\sigma>0$ and positive constant $C_\sigma$. Similarly, one can derive that
\begin{equation}
\begin{array}{ll}
\di\int Z_3^2 dx+\int_0^t\int hZ_3^2dxd\tau\le C(\mathcal{E}(0)^2+\d_0^{\frac12})
+C\int\int\f{|\wt{\mb{G}}|^2}{\mb{M}_*}(t,x,v)d v
dx+\s\int_0^t\|Z_{3\tau}\|^2 d\tau\\[3mm]
\di\qquad
+C\int_0^t\left[\|\sqrt{|\l_{1x}^{S_1}|+|\l_{3x}^{S_3}|}Z\|^2+\|Z_x\|^2\right] d\tau+C(\d_0+\chi_{\scriptscriptstyle T} )\int_0^t\|(\phi,\psi,\o)\|^2d\tau\\
\di\qquad
+C_\s\int_0^t\int\int\f{\nu(| v|)}{\mb{M}_*}|\wt{\mb{G}}|^2d v
dxd\tau+C\chi_{\scriptscriptstyle T} \di\sum_{|\a\prime|=1}\int_0^t\|\partial^{\a\prime}(\phi,\psi,\o)\|^2d\tau
+C\int_0^t\int |\T_x|(Z_1^2+Z_2^2) dxd\tau\\
\di\qquad+C\chi_{\scriptscriptstyle T} \int_0^t\int\int \f{\nu(|v|)|\partial_{v_1}(\mb{P}_c F_2)|^2}{\mb{M}_*} dxdvd\tau+C(\chi_{\scriptscriptstyle T} +\delta_0)\int_0^t \|(\Pi_x,n_2,\Pi_{x\tau})\|^2d\tau.
\end{array}\label{intz3}
\end{equation}
Then we estimate $\int_0^t\int|\T_x|Z_2^2dxd\tau$. Set
$$
\di \eta_2(x,t)=\int_{-\i}^x h(y,t)dy,
$$
it is easy to check that
$$
\eta_{2t}=ah_x-u_{1\#}h,\qquad \|n\|_{L^{\i}}=1.
$$
The following lemma is from Huang-Li-Matsumura \cite{Huang-Li-matsumura} for the stability of viscous contact wave,
which plays an important role for weighted characteristic estimate to $\int_0^t\int|\T_x|Z_2^2dxd\tau$.

\begin{lemma}\label{huang}
For $0<\tau\le \i$, suppose that $Z(t,x)$ satisfies
$$
Z\in L^{\i}(0,t;L^2(\mb{R})),\quad Z_x\in L^2(0,t;L^2(\mb{R})),\quad
Z_t\in L^2(0,t;H^{-1}(\mb{R})).
$$
Then the following estimate holds for any $\tau\in(0,t]$,
\begin{equation*}
\begin{array}{ll}
\di\int_0^t\int h^2Z^2dxd\tau \le\f{1}{a}\int Z^2(x,0)dx+4\int_0^t\|Z_x\|^2 d\tau+\f{2}{a}\left(\int_0^t<Z_\tau,Zn^2>d\tau-u_{1\#}\int_0^t\int Z^2nh dxd\tau \right).
\end{array}
\end{equation*}
\end{lemma}
From \eqref{Z-sys}, we have
\begin{equation*}\label{z2t}
Z_{2t}+\tilde\l_2Z_{2x}=\sum_{j=1}^3 c_{2j}Z_{jxx}-\sum_{j=1}^3b_{2j}Z_{j}+(2LA_3R_xZ_x+LA_4)_2.
\end{equation*}
we can get that
\begin{equation*}
\begin{array}{ll}
\di\int_0^t<Z_{2\tau},Z_2n^2>d\tau-u_{1\#}\int_0^t\int Z_2^2nh dxd\tau\\[3mm]
\di=\int_0^t\int\Bigg\{
\big(\sum_{j=1}^3 c_{2j}Z_{jxx}-\sum_{j=1}^3b_{2j}Z_{j}+(2LA_3R_xZ_x+LA_4)_2\big)Z_2n^2-(\tilde\l_2Z_{2x}Z_2n^2+u_{1\#}Z_2^2nh) \Bigg\} dxd\tau.
\end{array}
\end{equation*}
Taking $Z=Z_2$ in Lemma \ref{huang} and using the above equation, one can get
\begin{equation}
\begin{array}{ll}
\di\int_0^t\int h^2Z_2^2dxd\tau\le
C_\s\int\int\f{|\wt{\mb{G}}|^2}{\mb{M}_*}(t,x,v)d v
dx+C\int_0^t\left[\|\sqrt{|\l_{1x}^{S_1}|+|\l_{3x}^{S_3}|}Z\|^2+\|Z_x\|^2\right] d\tau\\[3mm]
\di
+\s\left[\|Z_2\|^2+\int_0^t\|Z_{2\tau}\|^2 d\tau\right]+C_\s\int_0^t\int\int\f{\nu(| v|)}{\mb{M}_*}|\wt{\mb{G}}|^2d v
dxd\tau\\[3mm]
\di
+C\chi_{\scriptscriptstyle T} \di\sum_{|\a\prime|=1}\int_0^t\|\partial^{\a\prime}(\phi,\psi,\o)\|^2d\tau+C\int_0^t\int |\T_x|(Z_1^2+Z_3^2)dxd\tau+C(\mathcal{E}(0)^2+\d_0^{\frac12})\\[3mm]
\di+C\chi_{\scriptscriptstyle T} \int_0^t\int\int \f{\nu(|v|)|\partial_{v_1}(\mb{P}_c F_2)|^2}{\mb{M}_*} dxdvd\tau+C(\d_0+\chi_{\scriptscriptstyle T} )\int_0^t\|(\phi,\psi,\o,\Pi_x,n_2,\Pi_{x\tau})\|^2d\tau.
\end{array}\label{intz2}
\end{equation}
Combining \eqref{intz1}, \eqref{intz3} and \eqref{intz2},
noting that \eqref{z1z3} and \eqref{note1} and choosing $\d_0$ and $\s$ sufficiently small, one has
%\begin{equation}\label{4.75}
%\begin{array}{ll}
%\di\|Z(t,\cdot)\|^2+|\a_2|\int_0^t\int\big[h(Z_1^2+Z_3^2)+h^2Z_2^2\big] dxd\tau
%+\int_0^t\|\sqrt{|\l^{S_1}_{1x}|+|\l^{S_3}_{3x}|}Z\|^2 d\tau\\[2mm]
%\di\quad\le C\int\int\f{|\wt{\mb{G}}|^2}{\mb{M}_*}(t,x,v)d v dx
%+C(\s+\sqrt{\d_0})\int_0^t\|(Z_{\tau},Z_x)\|^2 d\tau\\[2mm]
%\di\quad+C_\s\int_0^t\int\int\f{\nu(| v|)}{\mb{M}_*}|\wt{\mb{G}}|^2d v
%dxd\tau
%+C(\d_0+\chi_{\scriptscriptstyle T} )\int_0^t\|(\phi,\psi,\o)\|^2d\tau\\
%\di\quad+C\chi_{\scriptscriptstyle T} \di\sum_{|\a\prime|=1}\int_0^t\|\partial^{\a\prime}(\phi,\psi,\o)\|^2d\tau
%+C(\mathcal{E}(0)^2+\d_0^{\f12}).
%\end{array}
%\end{equation}
%By \eqref{note1}, it holds that
\begin{equation}
\begin{array}{ll}\label{zover}
\di \|Z(t,\cdot)\|^2+\int_0^t \|\sqrt{|\l_{1x}^{S_1}|+|\l_{3x}^{S_3}|+|\T_x|}~Z\|^2 d\tau
\le C\int\int\f{|\wt{\mb{G}}|^2}{\mb{M}_*}(t,x,v)d v dx+C(\mathcal{E}(0)^2+\d_0^{\f12})\\[3mm]
\di
+C(\s+\sqrt{\d_0})\int_0^t\|(Z_{\tau},Z_x)\|^2 d\tau+C\chi_{\scriptscriptstyle T} \di\sum_{|\a\prime|=1}\int_0^t\|\partial^{\a\prime}(\phi,\psi,\o)\|^2d\tau
+C_\s\int_0^t\int\int\f{\nu(| v|)}{\mb{M}_*}|\wt{\mb{G}}|^2d v
dxd\tau\\[2mm]
\di+C(\d_0+\chi_{\scriptscriptstyle T} )\int_0^t\|(\phi,\psi,\o,\Pi_x,n_2,\Pi_{x\tau})\|^2d\tau+C\chi_{\scriptscriptstyle T} \int_0^t\int\int \f{\nu(|v|)|\partial_{v_1}(\mb{P}_c F_2)|^2}{\mb{M}_*} dxdvd\tau.
\end{array}
\end{equation}
Noting that $Z$ is a linear vector function of $\Phi$, $\wt{\Psi}$, $\wt{W}$ and combining \eqref{le-s} and \eqref{zover}, we can get
\begin{equation*}
\begin{array}{l}
\di \|(\Phi,\wt\Psi,\wt W)(t,\cdot)\|^2+\int_0^t\left[\|\sqrt{|u^{S_1}_{1x}|
+|u^{S_3}_{1x}|+|\T_x|}~(\Phi,\wt\Psi,\wt W)\|^2+\|(\wt\Psi_x,\wt W_x)\|^2\right]d\tau\\[3mm]
\di\quad \leq C\int\int\f{|\wt{\mb{G}}|^2}{\mb{M}_*}(t,x,v)d v
dx+C(\s+\sqrt{\d_0}) \int_0^t\|(\Phi,\wt\Psi,\wt W)_{\tau}\|^2d\tau+C(\chi_{\scriptscriptstyle T} +\delta_0)\int_0^t \|(\Pi_x,n_2,\Pi_{x\tau})\|^2d\tau
\\[3mm]
\di\quad+C(\chi_{\scriptscriptstyle T} +\d_0)\di\sum_{|\a\prime|=1}\int_0^t\|\partial^{\a\prime}(\phi,\psi,\o)\|^2d\tau
+ C_\s\int_0^t\int\int\f{\nu(| v|)}{\mb{M}_*}|\wt{\mb{G}}|^2d v
dxd\tau\\[2mm]
\di\quad+ C(\s+\chi_{\scriptscriptstyle T} +\sqrt{\d_0})\di\int_0^t\|(\phi,\psi,\o)\|^2d\tau+C\chi_{\scriptscriptstyle T} \int_0^t\int\int \f{\nu(|v|)|\partial_{v_1}(\mb{P}_c F_2)|^2}{\mb{M}_*} dxdvd\tau+C(\mathcal{E}(0)^2+\d_0^{\f12}).
\end{array}
\end{equation*}

\underline{Step 4. Estimation on $\|\Phi_x(t,\cdot)\|^2$.}

Note that the dissipation term does not contain the term
$\|\Phi_x\|^2$.  From $\eqref{P-sys}_2$, we have
\begin{equation}
\begin{array}{ll}
\di\f{4}{3}\f{\mu(\tilde\t)}{\tilde\r}\Phi_{xt}+\tilde\r\wt\Psi_{1t}+\tilde\r\tilde
u_1\wt\Psi_{1x}+\f23\tilde\t\Phi_x
=-\f{4\mu(\tilde\t)}{3\tilde\r}(2\tilde\r_x\wt\Psi_{1x}+\tilde\r_{xx}\wt\Psi_1)-\f{4\mu(\tilde\t)}{3\tilde\r}(\tilde u_1\Phi)_{xx}\\
\di+ \f13\tilde\r\tilde u_{1x}\wt\Psi_1 -\f23\tilde\r_x\wt W-\f23\tilde\r\wt
W_x +\f{2\tilde\t\tilde\r_x}{3\tilde\r}\Phi
-\int v_1^2(\G-\G^{S_1}-\G^{S_3})dv+J_1+N_1-Q_1.
\end{array}
\label{4.1.76}
\end{equation}
Multiplying (\ref{4.1.76}) by $\Phi_x$ yields that
\begin{equation}
\begin{array}{l}
\quad\di\left(\f{2\mu(\tilde\t)}{3
\tilde\r}\Phi_x^2+\tilde\r\Phi_x\wt\Psi_1\right)_{t}+
\f{2\tilde\t}{3}\Phi_x^2 =\left(\f{2\mu(\tilde\t)}{3
\tilde\r}\right)_{t}\Phi_x^2+(\tilde\r\wt\Psi_{1})_x^2+(\tilde\r\wt\Psi_1)_x\tilde  u_{1x}\Phi\\[2mm]
\di\quad+\Bigg[-\f{4\mu(\tilde\t)}{3\tilde\r}(2\tilde\r_x\wt\Psi_{1x}+\tilde\r_{xx}\wt\Psi_1)-\f{4\mu(\tilde\t)}{3\tilde\r}(\tilde  u_1\Phi)_{xx}
- \f23\tilde\r\tilde u_{1x}\wt\Psi_1 -\f23\tilde\r_x\wt W-\f23\tilde\r\wt
W_x \\[2mm]
\di\quad+\f{2\tilde\t\tilde\r_x}{3\tilde\r}\Phi
-\int v_1^2(\G-\G^{S_1}-\G^{S_3})dv+J_1+N_1-Q_1\Bigg]\Phi_x+(\cdots)_x.
\end{array}
\label{(4.39)}
\end{equation}
By using Lemmas \ref{Lemma 4.1}-\ref{Lemma 4.3},  one has,
\begin{equation}
\begin{array}{ll}
\di \int_0^t\int\left|\int v_1^2(\G-\G^{S_1}-\G^{S_3})d v\right|^2dxd\tau\leq
 C\d_0^2+C\d_0\int_0^t\|(\phi,\psi,\o)\|^2 d\tau\\
 \di\quad+C(\d_0+\chi_{\scriptscriptstyle T} )\int_0^t\int\int\f{\nu(| v|)}{\mb{M}_*}|\wt{\mb{G}}|^2d v
dxd\tau+C\sum_{|\a^\prime|=1}\int_0^t\int\int\f{\nu(| v|)}{\mb{M}_*}|\partial^{\a^\prime}\wt{\mb{G}}|^2d v
dxd\tau\\
\di\quad +C\chi_{\scriptscriptstyle T} \int_0^t\int\int \f{\nu(|v|)|\partial_{v_1}(\mb{P}_c F_2)|^2}{\mb{M}_*} dxdvd\tau+C(\chi_{\scriptscriptstyle T} +\delta_0)\int_0^t \|(\Pi_x,n_2,\Pi_{x\tau})\|^2d\tau.
\end{array}
\label{K2-E}
\end{equation}
Integrating (\ref{(4.39)}) with respect to $x$ and $t$, which together with
H${\rm\ddot{o}}$lder's  inequality and (\ref{K2-E}) gives that
\begin{equation*}
\begin{array}{l}
\di
\|\Phi_x(t,\cdot)\|^2+\int_0^t\|\Phi_x\|^2d\tau\leq C\|\wt\Psi_1(t,\cdot)\|^2+C\int_0^t\|(\wt\Psi_x,\wt W_x)\|^2d\tau+C\chi_{\scriptscriptstyle T} \int_0^t \|(\phi_x,\psi_{1x})\|^2 d\tau\\[3mm]
\di\quad +C\d_0\int_0^t\|\sqrt{|u^{S_1}_{1x}|+|u^{S_3}_{1x}|+|\T_{x}|}(\Phi,\wt\Psi_1,\wt W)\|^2d\tau
+C\chi_{\scriptscriptstyle T} \int_0^t\int\int \f{\nu(|v|)|\partial_{v_1}(\mb{P}_c F_2)|^2}{\mb{M}_*} dxdvd\tau\\
\di\quad +C(\chi_{\scriptscriptstyle T} +\d_0)\int_0^t \|(\psi,\o,\Pi_x,n_2,\Pi_{x\tau})\|^2 d\tau
+C(\d_0+\chi_{\scriptscriptstyle T} )\int_0^t\int\int\f{\nu(| v|)}{\mb{M}_*}|\wt{\mb{G}}|^2d v
dxd\tau
+C(\mathcal{E}(0)^2+\d_0^{\f12})\\
\di\quad +C\sum_{|\a^\prime|=1}\int_0^t\int\int\f{\nu(| v|)}{\mb{M}_*}|\partial^{\a^\prime}\wt{\mb{G}}|^2d v
dxd\tau+\int_0^t\int q|(\Phi,\wt \Psi, \wt W)|^2dxd\tau.
\end{array}
\end{equation*}

\underline{Step 5. Estimation on the non-fluid component.}

Next we do the microscopic estimates for the Vlasov-Poisson-Boltzmann system.
Multiplying the equation \eqref{G} and the equation \eqref{F2-pc-1} by
$\f{\widetilde{\mb{G}}}{\mb{M_*}}$ and $\f{\mb{P}_c
F_2}{\mb{M}_*}$ , respectively, one has
\begin{equation}\label{Gle0-s}
\begin{array}{ll}
\di\Big(\frac{|\wt{\mb{G}}|^2}{2\mb{M}_*}\Big)_{t}-\f{\wt{\mb{G}}}{\mb{M}_*}\mb{L}_\mb{M}\wt{\mb{G}}
=\Big\{-\mb{P}_1( v_1\wt{\mb{G}}_x)-\mb{P}_1( \Pi_x \partial_{v_1}F_2)+2[Q(\mb{G}^{S_1},\mb{G}^{S_3})+Q(\mb{G}^{S_3},\mb{G}^{S_1})]
\\
\di \quad+2[Q(\wt{\mb{G}},\mb{G}^{S_1}+\mb{G}^{S_3})+Q(\mb{G}^{S_1}+\mb{G}^{S_3},\wt{\mb{G}})]-\big[\mb{P}_1( v_1\mb{M}_x)-\mb{P}_1^{S_1}( v_1\mb{M}_x^{S_1})
-\mb{P}_1^{S_3}( v_1\mb{M}_x^{S_3})\big]\\
\di \quad+2Q(\wt{\mb{G}},\wt{\mb{G}})+\sum_{j=1,3}R_j\Big\}\f{\wt{\mb{G}}}{\mb{M}_*}.
\end{array}
\end{equation}
and
\begin{equation}\label{Pc-F2-le-s}
\begin{array}{ll}
\di \Big(\f{|\mb{P}_c F_2|^2}{2\mb{M}_*}\Big)_t-\f{\mb{P}_c
F_2}{\mb{M}_*}\mb{N}_{\mb{M}}(\mb{P_c}F_2)=\Big[-v_1\partial_x
F_2-(\f{n_2}{\rho}\mb{M})_t-\mb{P}_c(\Pi_x\partial_{v_1}F_1)+2Q(F_2,\mb{G})\Big]\f{\mb{P}_c F_2}{\mb{M}_*},
\end{array}
\end{equation}
where in the equation \eqref{Pc-F2-le-s}, one has used the fact that
$$
\begin{array}{ll}
\di -\mb{P}_c (v_1\partial_x
F_2)-(\f{\mb{M}}{\rho})_t n_2=-v_1\partial_x
F_2+\mb{P}_d (v_1\partial_x
F_2)-(\f{\mb{M}}{\rho})_t n_2\\
\di\qquad  =-v_1\partial_x
F_2+\f{\mb{M}}{\rho} \big(\int v_1
F_2dv\big)_x-(\f{\mb{M}}{\rho})_t n_2=-v_1\partial_x
F_2-\f{\mb{M}}{\rho} n_{2t}-(\f{\mb{M}}{\rho})_t n_2 =-v_1\partial_x
F_2-(\f{n_2}{\rho}\mb{M})_t.
\end{array}
$$
It can be computed that
\begin{equation}\label{P1-1-s}
\begin{array}{ll}
\di \mb{P}_1(\Pi_x\partial_{v_1}F_2)=\Pi_x\mb{P}_1\big[\partial_{v_1}\big(\f{n_2}{\rho}\mb{M}+\mb{P}_cF_2\big)\big]=\Pi_x\mb{P}_1\big[\partial_{v_1}\big(\mb{P}_cF_2\big)\big]\\
\di =\Pi_x\big(\mb{P}_cF_2\big)_{v_1}-\Pi_x\sum_{j=0}^4\int\big(\mb{P}_cF_2\big)_{v_1}\f{\chi_j}{\mb{M}}dv\chi_j =\Pi_x\big(\mb{P}_cF_2\big)_{v_1}+\Pi_x\sum_{j=0}^4\int(\mb{P}_cF_2) \big(\f{\chi_j}{\mb{M}}\big)_{v_1}dv\chi_j,
\end{array}
\end{equation}
and
\begin{equation}\label{Pc-1-s}
\di \mb{P}_c(\Pi_x\partial_{v_1}F_1)=\Pi_x\partial_{v_1}F_1-\Pi_x\mb{P}_d\big(\partial_{v_1}F_1\big)=\Pi_x\partial_{v_1}F_1=\Pi_x \mb{M}_{v_1}+\Pi_x \widetilde{\mb{G}}_{v_1}+\Pi_x (\mb{G}^{S_1}_{v_1}+\mb{G}^{S_3}_{v_1}).
\end{equation}
Substituting \eqref{P1-1-s} and \eqref{Pc-1-s} into \eqref{Gle0-s} and \eqref{Pc-F2-le-s}, respectively, then summing the resulting equations together and noting the following cancelation
$$
\Pi_x\big(\mb{P}_cF_2\big)_{v_1}\f{\widetilde{\mb{G}}}{\mb{M}_*}+\Pi_x \widetilde{\mb{G}}_{v_1}\f{\mb{P}_c F_2}{\mb{M}_*}=\Big(\Pi_x\f{\mb{P}_cF_2\cdot\widetilde{\mb{G}}}{\mb{M}_*}\Big)_{v_1}+\Pi_x\f{\mb{P}_cF_2\cdot\widetilde{\mb{G}}}{\mb{M}_*^2}\big(\mb{M}_*\big)_{v_1},
$$
it holds that
\begin{equation*}\label{Gle}
\begin{array}{ll}
\di \Big(\f{|\widetilde{\mb{G}}|^2+|\mb{P}_c F_2|^2}{2\mb{M}_*}\Big)_t-\f{\widetilde{\mb{G}}}{\mb{M}_*}\mb{L}_{\mb{M}}\widetilde{\mb{G}}-\f{\mb{P}_c
F_2}{\mb{M}_*}\mb{N}_{\mb{M}}(\mb{P_c}F_2)+(\cdots)_{v_1}=\Big\{-\mb{P}_1( v_1\wt{\mb{G}}_x)+2Q(\wt{\mb{G}},\wt{\mb{G}})\\
\quad \di+2[Q(\mb{G}^{S_1},\mb{G}^{S_3})+Q(\mb{G}^{S_3},\mb{G}^{S_1})]
+2[Q(\wt{\mb{G}},\mb{G}^{S_1}+\mb{G}^{S_3})+Q(\mb{G}^{S_1}+\mb{G}^{S_3},\wt{\mb{G}})]\\
\di\quad+\sum_{j=1,3}R_j-\big[\mb{P}_1( v_1\mb{M}_x)-\mb{P}_1^{S_1}( v_1\mb{M}_x^{S_1})
-\mb{P}_1^{S_3}( v_1\mb{M}_x^{S_3})\big]-\Pi_x\sum_{j=0}^4\int(\mb{P}_cF_2) \big(\f{\chi_j}{\mb{M}}\big)_{v_1}dv\chi_j\Big\}\f{\widetilde{\mb{G}}}{\mb{M}_*}\\
\di+\Big[-v_1\partial_x
F_2-(\f{n_2}{\rho}\mb{M})_t+2Q(F_2,\mb{G})-\Pi_x \mb{M}_{v_1}-\Pi_x (\mb{G}^{S_1}_{v_1}+\mb{G}^{S_3}_{v_1})\Big]\f{\mb{P}_c F_2}{\mb{M}_*}+\Pi_x\f{\mb{P}_cF_2\cdot\widetilde{\mb{G}}}{\mb{M}_*^2}\big(\mb{M}_*\big)_{v_1}.
\end{array}
\end{equation*}
Integrating the above
equation with respect to $x,t,v$ yields that
\begin{equation}\label{GE0-s}
\begin{array}{ll}
\di \int\int \f{|\widetilde{\mb{G}}|^2+|\mb{P}_c F_2|^2}{2\mb{M_*}}(x,v,t)
dxdv-\int\int \f{|\widetilde{\mb{G}}|^2+|\mb{P}_c F_2|^2}{2\mb{M_*}}(x,v,0)
dxdv\\
\di  -\int_0^t\int\int \Big[\f{\widetilde{\mb{G}}}{\mb{M}_*}\mb{L}_{\mb{M}}\widetilde{\mb{G}}+\f{\mb{P}_c
F_2}{\mb{M}_*}\mb{N}_{\mb{M}}(\mb{P_c}F_2)\Big]
dxdvd\tau
=-\int_0^t\int\int\Bigg\{\Big[-\mb{P}_1( v_1\wt{\mb{G}}_x)+2Q(\wt{\mb{G}},\wt{\mb{G}})\\
 \di+2[Q(\mb{G}^{S_1},\mb{G}^{S_3})+Q(\mb{G}^{S_3},\mb{G}^{S_1})]
+2[Q(\wt{\mb{G}},\mb{G}^{S_1}+\mb{G}^{S_3})+Q(\mb{G}^{S_1}+\mb{G}^{S_3},\wt{\mb{G}})]+\sum_{j=1,3}R_j\\
\di-\big[\mb{P}_1( v_1\mb{M}_x)-\mb{P}_1^{S_1}( v_1\mb{M}_x^{S_1})
-\mb{P}_1^{S_3}( v_1\mb{M}_x^{S_3})\big]-\Pi_x\sum_{j=0}^4\int(\mb{P}_cF_2) \big(\f{\chi_j}{\mb{M}}\big)_{v_1}dv\chi_j\Big]\f{\widetilde{\mb{G}}}{\mb{M}_*}+\Big[-v_1\partial_x
F_2\\
\di -(\f{n_2}{\rho}\mb{M})_\tau+2Q(F_2,\mb{G})-\Pi_x \mb{M}_{v_1}-\Pi_x (\mb{G}^{S_1}_{v_1}+\mb{G}^{S_3}_{v_1})\Big]\f{\mb{P}_c F_2}{\mb{M}_*} +\Pi_x\f{\mb{P}_cF_2\cdot\widetilde{\mb{G}}}{\mb{M}_*^2}\big(\mb{M}_*\big)_{v_1}\Bigg\}dxdvd\tau\\
\di :=\sum_{i=1}^{13}Y_i.
\end{array}
\end{equation}
Now we calculate the right hand side of \eqref{GE0-s}.  The estimations of $Y_i~(i=1,2,\cdots 7, 10, 12, 13)$ are standard and will be skipped for brevity.
Noting that
$$
\begin{array}{ll}
\di -v_1\partial_x F_2\f{\mb{P}_c F_2}{\mb{M}_*}=-v_1\partial_x
\big(\f{n_2}{\r}\mb{M}+\mb{P}_cF_2\big)\f{\mb{P}_c
F_2}{\mb{M}_*}=-v_1 \big(\f{n_2}{\r}\mb{M}\big)_x\f{\mb{P}_c
F_2}{\mb{M}_*}-\Big(v_1\f{|\mb{P}_c F_2|^2}{2\mb{M}_*}\Big)_x\\
\di \qquad=-v_1 \f{n_{2x}\mb{M}}{\r}\f{\mb{P}_c F_2}{\mb{M}_*}-
n_{2}v_1\big(\f{\mb{M}}{\r}\big)_x\f{\mb{P}_c
F_2}{\mb{M}_*}+(\cdots)_x\\
\di \qquad=-
\f{n_{2x}}{\r}v_1\big(\f{\mb{M}}{\mb{M}_*}-1\big)\mb{P}_c F_2-
n_{2}v_1\big(\f{\mb{M}}{\r}\big)_x\f{\mb{P}_c
F_2}{\mb{M}_*}-\f{n_{2x}}{\r}v_1\mb{P}_c F_2+(\cdots)_x,
\end{array}
$$
it follows from \eqref{phi-xt-e-s} that
\begin{equation}\label{Y8}
\begin{array}{ll}
\di Y_8=-\int_0^t\int\int v_1\partial_x F_2\f{\mb{P}_c F_2}{\mb{M}_*}dxdvd\tau \\
\di\quad =\int_0^t\int\int \Big[-
\f{n_{2x}}{\r}v_1\big(\f{\mb{M}}{\mb{M}_*}-1\big)\mb{P}_c F_2-
n_{2}v_1\big(\f{\mb{M}}{\r}\big)_x\f{\mb{P}_c
F_2}{\mb{M}_*}-\f{n_{2x}}{\r}v_1\mb{P}_c F_2\Big]dvdxd\tau\\
\di\quad = \int_0^t\int\int \Big[-
\f{n_{2x}}{\r}v_1\big(\f{\mb{M}}{\mb{M}_*}-1\big)\mb{P}_c F_2-
n_{2}v_1\big(\f{\mb{M}}{\r}\big)_x\f{\mb{P}_c
F_2}{\mb{M}_*}\Big]dvdxd\tau+\int_0^t\int \f{n_{2x}}{\r}(\f12\Pi_{x\tau}+u_1n_2) dxd\tau.
\end{array}
\end{equation}
By integration by parts, we can compute that
\begin{equation}\label{Y81}
\begin{array}{ll}
\di \int_0^t\int\f{n_{2x}}{\r}(\f12\Pi_{x\tau}+u_1n_2)
dxd\tau=-\int_0^t\int \Big[
\f{n_{2}}{\r}n_{2\tau}+\f12n_2\Pi_{x\tau}(\f1\r)_x+\big(\f{u_1}{2\r}\big)_xn_2^2\Big]
dxd\tau\\
\qquad \di=-\int\f{n_2^2}{2\r}(x,t)dx+\int\f{n_{20}^2}{2\r_0}
dx-\int_0^t\int \Big[
\big(\f{1}{2\r}\big)_\tau n^2_{2}+\f12n_2\Pi_{x\tau}(\f1\r)_x+\big(\f{u_1}{2\r}\big)_xn_2^2\Big]
dxd\tau.
\end{array}
\end{equation}
Substituting \eqref{Y81} into \eqref{Y8} and estimating the other
terms yield that
\begin{equation*}\label{Y8-E}
\begin{array}{ll}
\di Y_8+\int\f{n_2^2}{2\r}(x,t)dx \leq C\|n_{20}\|^2
+C(\chi_{\scriptscriptstyle T} +\delta_0)\sum_{|\alpha^\prime|=1}\int_0^t\|\partial^{\alpha^\prime}(\phi,\psi,\omega)\|^2
d\tau\\
\di\quad +C(\chi_{\scriptscriptstyle T} +\delta_0)\int_0^t\|(n_2,\Pi_{x\tau})\|^2d\tau
+C(\chi_{\scriptscriptstyle T} +\delta_0+\eta_0)\int_0^t\Big[\|n_{2x}\|^2+\int\int
\f{\nu(|v|)|\mb{P}_c F_2|^2}{\mb{M}_*}dxdv\Big]d\tau.
\end{array}
\end{equation*}
Similarly, it holds that
\begin{equation*}\label{Y9}
\begin{array}{ll}
\di |Y_9|= |\int_0^t\int\int
\Big[\f{n_{2\tau}}{\rho}\f{\mb{M}}{\mb{M}_*} \mb{P}_c
F_2+n_2\big(\f{\mb{M}}{\r}\big)_\tau\f{\mb{P}_c
F_2}{\mb{M}_*}\Big]dvdxd\tau|\\
\di \qquad= |\int_0^t\int\int
\Big[\f{n_{2\tau}}{\rho}\big(\f{\mb{M}}{\mb{M}_*}-1\big) \mb{P}_c
F_2+n_2\big(\f{\mb{M}}{\r}\big)_t\f{\mb{P}_c
F_2}{\mb{M}_*}\Big]dvdxd\tau|\\ \di \qquad\leq
C(\chi_{\scriptscriptstyle T} +\delta_0+\eta_0)\int_0^t\Big[\|n_{2\tau}\|^2+\int\int
\f{\nu(|v|)|\mb{P}_c F_2|^2}{\mb{M}_*}dxdv\Big]d\tau+C(\chi_{\scriptscriptstyle T} +\delta_0) \int_0^t\|(\Pi_{x\tau},n_2)\|^2d\tau\\
\di\qquad\quad +C(\chi_{\scriptscriptstyle T} +\delta_0) \sum_{|\alpha^\prime|=1}\int_0^t
\|\partial^{\alpha^\prime}(\phi,\psi,\omega)\|^2d\tau.
\end{array}
\end{equation*}
Again it follows from \eqref{phi-xt-e-s} that
\begin{equation*}\label{Y11}
\begin{array}{ll}
\di Y_{11}=\int_0^t\int\int \Pi_x \mb{M}\f{v_1-u_1}{R\t}\f{\mb{P}_c
F_2}{\mb{M}_*} dvdxd\tau\\
\di \quad=\int_0^t\int\int
\Pi_x\f{v_1-u_1}{R\t}\big(\f{\mb{M}}{\mb{M}_*}-1\big)\mb{P}_c F_2
dvdxd\tau+\int_0^t\int\int \Pi_x\f{v_1-u_1}{R\t}\mb{P}_c F_2 dvdxd\tau
\\
\di \quad=\int_0^t\int\int
\Pi_x\f{v_1-u_1}{R\t}\big(\f{\mb{M}}{\mb{M}_*}-1\big)\mb{P}_c F_2
dvdxd\tau+\int_0^t\int\f{\Pi_x}{R\t}\big(\int v_1\mb{P}_c F_2
dv\big)dxd\tau\\
\di\quad =\int_0^t\int\int
\Pi_x\f{v_1-u_1}{R\t}\big(\f{\mb{M}}{\mb{M}_*}-1\big)\mb{P}_c F_2
dvdxd\tau-\int_0^t\int\f{\Pi_x}{R\t}(\f12\Pi_{xt}+u_1n_2)dxd\tau\\
\di\quad \leq -\int\f{\Pi_x^2}{4R\t}
dx+C\|\Pi_{x0}\|^2+C(\chi_{\scriptscriptstyle T} +\delta_0+\eta_0)\int_0^t\Big[\|\Pi_x\|^2+\int\int
\f{\nu(|v|)|\mb{P}_c
F_2|^2}{\mb{M}_*}dxdv\Big]d\tau\\
\di\qquad+C(\chi_{\scriptscriptstyle T} +\delta_0)\int_0^t\|(\Pi_x,\Pi_{x\tau},n_2,\omega_\tau,\psi_{1x},\omega_x)\|^2d\tau,
\end{array}
\end{equation*}
where in the last inequality one has used the fact that
$$
\begin{array}{ll}
\di
-\int_0^t\int\f{\Pi_x}{R\t}(\f12\Pi_{x\tau}+u_1n_2)dxd\tau=-\int\f{\Pi_x^2}{4R\t} dx+\int\f{\Pi_{x0}^2}{4R\t_0}
dx+\int_0^t\int
\Big[\big(\f{1}{4R\t}\big)_\tau+\big(\f{u_1}{4R\t}\big)_x\Big]\Pi_x^2dxd\tau\\
\di\qquad\qquad  \leq -\int\f{\Pi_x^2}{4R\t}
dx+C\|\Pi_{0x}\|^2+C(\chi_{\scriptscriptstyle T} +\delta_0)\int_0^t\|(\Pi_x,\Pi_{x\tau},n_2,\omega_\tau,\psi_{1x},\omega_x)\|^2d\tau.
\end{array}
$$

Substituting the above estimations  into \eqref{GE0-s}, then choosing
$\chi_{\scriptscriptstyle T} ,\delta,\eta_0$ sufficiently small imply that
\begin{equation}\label{FE}
\begin{array}{ll}
\di \int\int \f{|(\widetilde{\mb{G}},\mb{P}_c
F_2)|^2}{\mb{M}_*}(x,v,t)dxdv+\|(\Pi_x,n_2)(\cdot,t)\|^2+\int_0^t\int\int\frac{\nu(|v|)|\big(\widetilde{\mb{G}},
\mb{P}_c
F_2\big)|^2}{\mb{M}_*}dvdxd\tau\\
\di \leq C\int\int \f{|(\widetilde{\mb{G}}, \mb{P}_c
F_2)|^2}{\mb{M}_*}(x,v,0)dxdv+C\|(\Pi_{x0},n_{20})\|^2+
C\int_0^t\|(\phi_x,\psi_x,\omega_x)\|^2d\tau+C\delta_0^{\f12}\\
\di +C(\chi_{\scriptscriptstyle T} +\delta_0) \int_0^t
\Big[\|(\Pi_{x\tau},n_2)\|^2+\sum_{|\alpha|=1}\|\partial^\alpha(\phi,\psi,\omega)\|^2\Big]d\tau\\
\di +C(\chi_{\scriptscriptstyle T} +\delta_0+\eta_0)\int_0^t\|(\Pi_x,n_{2x},n_{2\tau})\|^2
d\tau+C\sum_{|\alpha^\prime|=1}\int_0^t\int\int\frac{\nu(|v|)|\partial^{\alpha^\prime} \wt{\mb{G}}|^2}{\mb{M}_*}dvdxd\tau.
\end{array}
\end{equation}

\underline{Step 6. Estimation of electric field terms.}

Now we estimate the Poisson term in the electric fields, which is one of the key ingredients of the present paper. Multiplying the equation
\eqref{n2} by $-\Pi$ yields that
\begin{equation*}
\begin{array}{ll}
\di
(\f{\Pi_x^2}{2})_t+\f{3\k_1(\t)}{2\t}\Pi_x^2  +\f{2\k_1(\t)}{\r}n_2^2
-\f14 u_{1x}\Pi_x^2=(\cdots)_x-\Pi_x\f{n_2}{\r}\int v_1 \mb{N}_{\mb{M}}^{-1}\big[\mb{P}_c(v_1\mb{M}_x)\big]dv\\
\di\quad  -\Pi_x\int v_1 \mb{N}_{\mb{M}}^{-1}\big[\mb{P}_c(v_1
(\mb{P}_cF_2)_x)\big]dv
 - \Pi_x^2\int v_1 \mb{N}_{\mb{M}}^{-1}\big[\mb{G}_{v_1}\big]dv\\
 \di\quad  -\Pi_x\Big(\int v_1 \mb{N}_{\mb{M}}^{-1}\Big[\partial_t(\mb{P}_c
 F_2)+(\f{\mb{M}}{\rho})_t~n_2-2Q(F_2,\mb{G})\Big]dv\Big).
\end{array}
\end{equation*}
Integration the above equation with respect to $x,t$ implies that
\begin{equation}\label{Phi-1}
\begin{array}{ll}
\di \|\Pi_x\|^2(t)-\|\Pi_x\|^2(0)+\int_0^t\|(\Pi_x,n_2)\|^2d\tau\leq
C(\chi_{\scriptscriptstyle T} +\delta)\int_0^t \|(\Pi_{x},n_2,\psi_{1x},\omega_x)\|^2d\tau\\
\di   +C\int_0^t\int|\Pi_x\f{n_2}{\r}\int v_1 \mb{N}_{\mb{M}}^{-1}\big[\mb{P}_c(v_1\mb{M}_x) )\big]dv|dxd\tau+C\int_0^t\int|\Pi_x\int v_1
\mb{N}_{\mb{M}}^{-1}\big[\mb{P}_c(v_1 (\mb{P}_cF_2)_x)\big]dv|dxd\tau\\
\di
+C\int_0^t\int| \Pi_x^2\int v_1 \mb{N}_{\mb{M}}^{-1}\big[\mb{G}_{v_1}\big]dv|dxd\tau +C\int_0^t\int|\Pi_x\int v_1 \mb{N}_{\mb{M}}^{-1}\Big[\partial_\tau(\mb{P}_c
 F_2)+(\f{\mb{M}}{\rho})_\tau~n_2\Big]dv|dxd\tau\\
 \di +C\int_0^t\int|\Pi_x\int v_1 \mb{N}_{\mb{M}}^{-1} Q(F_2,\mb{G})dv|dxd\tau:=C(\chi_{\scriptscriptstyle T} +\delta_0)\int_0^t \|(\Pi_{x},n_2,\psi_{1x},\omega_x)\|^2d\tau+\sum_{i=1}^5T_i.
\end{array}
\end{equation}
Now we estimate $T_i~(i=1,2,3,4,5)$ in \eqref{Phi-1} one by one. First, it holds that
\begin{equation}\label{T1}
\begin{array}{ll}
\di T_1\leq C \int_0^t\int |n_2||\Pi_x|\big(\int\f{\nu(|v|)|\mb{N}_{\mb{M}}^{-1}[\mb{P}_c(v_1\mb{M}_x)]|^2}{\mb{M}_*} dv\big)^{\f12}\big(\int \nu(|v|)^{-1}v_1^2\mb{M}_*dv\big)^{\f12}dxd\tau\\
\di \qquad \leq C\int_0^t\int  |n_2||\Pi_x| \big(\int\f{\nu(|v|)^{-1}|\mb{P}_c(v_1\mb{M}_x)|^2}{\mb{M}_*} dv\big)^{\f12}dxd\tau\\
\di \qquad \leq C\int_0^t\int  |n_2||\Pi_x| \big(\int\f{\nu(|v|)^{-1}|v_1\mb{M}_x-\f{\mb{M}}{\rho}(\r u_1)_x|^2}{\mb{M}_*} dv\big)^{\f12}dxd\tau\\
\di \qquad \leq C\int_0^t\int  |n_2||\Pi_x| |(\r_x,u_x,\t_x)|dxd\tau
\leq
C(\chi_{\scriptscriptstyle T} +\delta_0)\int_0^t\big[\|(n_2,\Pi_x)\|^2+\|(\phi_x,\psi_x,\omega_x)\|^2\big]dxd\tau,
\end{array}
\end{equation}
and
\begin{equation*}\label{T2}
\begin{array}{ll}
\di T_2\leq C \int_0^t\int |\Pi_x|\big(\int\f{\nu(|v|)|\mb{N}_{\mb{M}}^{-1}\big[\mb{P}_c(v_1 (\mb{P}_cF_2)_x)\big]|^2}{\mb{M}_{[\r_*,u_*,2\t_*]}} dv\big)^{\f12}\big(\int \nu(|v|)^{-1}v_1^2\mb{M}_{[\r_*,u_*,2\t_*]}dv\big)^{\f12}dxd\tau\\
\di \qquad \leq C\int_0^t\int  |\Pi_x| \big(\int\f{\nu(|v|)^{-1}|\mb{P}_c(v_1 (\mb{P}_cF_2)_x)|^2}{\mb{M}_{[\r_*,u_*,2\t_*]}} dv\big)^{\f12}dxd\tau\\
\di \qquad \leq C\int_0^t\int  |\Pi_x| \big(\int\f{\nu(|v|)^{-1}|(\mb{P}_cF_2)_x)|^2}{\mb{M}_*} dv\big)^{\f12}dxd\tau \\
\di \qquad \leq \f18\int_0^t\|\Pi_x\|^2d\tau+C~\int_0^t\int
\int\f{\nu(|v|)|(\mb{P}_cF_2)_x)|^2}{\mb{M}_*} dvdxd\tau.
\end{array}
\end{equation*}
Then one has
\begin{equation*}\label{T3}
\begin{array}{ll}
\di T_3\leq C \int_0^t\int |\Pi_x|^2\big(\int\f{\nu(|v|)^{-1}|\mb{G}_{v_1}|^2}{\mb{M}_*} dv\big)^{\f12}dxd\tau \\
\di \quad \leq
C(\chi_{\scriptscriptstyle T} +\delta_0)\int_0^t\Big[\|\Pi_x\|^2+\int\int\f{\nu(|v|)^{-1}|\widetilde{\mb{G}}_{v_1}|^2}{\mb{M}_*}
dv dx\Big]d\tau,
\end{array}
\end{equation*}
\begin{equation*}\label{T4}
\begin{array}{ll}
\di T_4\leq C \int_0^t\int |\Pi_x|\big(\int\f{\nu(|v|)^{-1}|\partial_\tau(\mb{P}_c
 F_2)+(\f{\mb{M}}{\rho})_\tau~n_2|^2}{\mb{M}_*} dv\big)^{\f12}dxd\tau\\
\di\quad \leq C(\d_0+\chi_{\scriptscriptstyle T} )\int_0^t\Big[\|(\Pi_x,\Pi_{x\tau},n_2)\|^2+\|(\phi_\tau,\psi_\tau,\omega_\tau)\|^2\Big]d\tau\\
\di \qquad+\f18\int\|\Pi_x\|^2d\tau +C~\int_0^t\int
\int\f{\nu(|v|)|(\mb{P}_cF_2)_t|^2}{\mb{M}_*} dvdxd\tau,
\end{array}
\end{equation*}
and
\begin{equation}\label{T5}
\begin{array}{ll}
\di T_5\leq C \int_0^t\int |\Pi_x|\big(\int\f{\nu(|v|)^{-1}|Q(F_2,\mb{G})|^2}{\mb{M}_*} dv\big)^{\f12}dxd\tau\\
\di \leq C\int_0^t\int|\Pi_x||n_2|(\int\f{\nu(|v|)|\mb{G}|^2}{\mb{M}_*}dv)^{\f12}dxd\tau+C\int_0^t\int |\Pi_x|(\int\f{\nu(|v|)|\mb{P}_cF_2|^2}{\mb{M}_*} dv)^{\f12}(\int\f{|\mb{G}|^2}{\mb{M}_*}dv)^{\f12}dxd\tau\\
\di\qquad +C\int_0^t\int |\Pi_x|(\int\f{|\mb{P}_cF_2|^2}{\mb{M}_*}
dv)^{\f12}(\int\f{\nu(|v|)|\mb{G}|^2}{\mb{M}_*}dv)^{\f12}dxd\tau\\
\di \leq C(\chi_{\scriptscriptstyle T} +\delta_0)\Big[
\int_0^t\|(\Pi_x,n_2)\|^2d\tau+\int_0^t\int
\int\f{\nu(|v|)|\mb{P}_cF_2|^2}{\mb{M}_*} dvdxd\tau+\int_0^t\int
\int\f{\nu(|v|)|\widetilde{\mb{G}}|^2}{\mb{M}_*} dvdxd\tau\Big].
\end{array}
\end{equation}
Substituting \eqref{T1}-\eqref{T5} into \eqref{Phi-1} and choosing
$\chi_{\scriptscriptstyle T} ,\delta$ sufficiently small yield that
\begin{equation}\label{Phi-E-s}
\begin{array}{ll}
\di \|\Pi_x\|^2(t)+\int_0^t\|(\Pi_x,n_2)\|^2d\tau\leq
C\|\Pi_x\|^2(0)+C(\chi_{\scriptscriptstyle T} +\delta_0)\int_0^t\|\Pi_{x\tau}\|^2 d\tau\\
\di \quad
+C(\chi_{\scriptscriptstyle T} +\delta_0)\sum_{|\alpha|=1}\int_0^t\|\partial^\alpha(\phi,\psi,\omega)\|^2d\tau+C~\sum_{|\alpha^\prime|=1}\int_0^t\int
\int\f{\nu(|v|)|\partial^{\alpha^\prime}(\mb{P}_cF_2))|^2}{\mb{M}_*}
dvdxd\tau\\
\di \quad +C(\chi_{\scriptscriptstyle T} +\delta_0)\Big[\int_0^t\int
\int\f{\nu(|v|)|(\widetilde{\mb{G}},\mb{P}_cF_2)|^2}{\mb{M}_*} dvdxd\tau+\int_0^t\int
\int\f{\nu(|v|)| \widetilde{\mb{G}}_{v_1}|^2}{\mb{M}_*}
dvdxd\tau\Big]
\end{array}
\end{equation}
Multiplying the equation \eqref{n2} by $n_2$ yields that
\begin{equation*}\label{n2-1}
\begin{array}{ll}
\di
(\f{n_2^2}{2})_t+\f{3\k_1(\t)}{2\t}n_2^2+(\f{3\k_1(\t)}{2\t})_x\Pi_xn_2-(\k_1(\t)(\f{n_2}{\rho})_xn_2)_x+\k_1(\t)(\f{n_2}{\rho})_xn_{2x}+(u_1n_2^2)_x\\
\di \quad -u_1n_2n_{2x}
=(\cdots)_x+n_{2x}\f{n_2}{\r}\int v_1 \mb{N}_{\mb{M}}^{-1}\big[\mb{P}_c(v_1\mb{M}_x) )\big]dv+n_{2x}\int v_1 \mb{N}_{\mb{M}}^{-1}\big[\mb{P}_c(v_1
(\mb{P}_cF_2)_x)\big]dv
\\
\di\quad   +n_{2x}\Pi_x\int v_1 \mb{N}_{\mb{M}}^{-1}\big[\mb{G}_{v_1}\big]dv +n_{2x}\Big(\int v_1 \mb{N}_{\mb{M}}^{-1}\Big[\partial_t(\mb{P}_c
 F_2)+(\f{\mb{M}}{\rho})_t~n_2-2Q(F_2,\mb{G})\Big]dv\Big).
\end{array}
\end{equation*}
Integrating the above equation with respect to $x,t$, one has
\begin{equation}\label{n2-E-s}
\begin{array}{ll}
\di \|n_{2}\|^2(t)+\int_0^t\|(n_{2x},n_2)\|^2d\tau\leq C\|n_{20}\|^2+ C(\chi_{\scriptscriptstyle T} +\delta_0)\int_0^t \|\Pi_x\|^2d\tau\\
\di
+C(\chi_{\scriptscriptstyle T} +\delta_0)\sum_{|\alpha|=1}\int_0^t\|\partial^\alpha(\phi,\psi,\omega)\|^2d\tau+C~\sum_{|\alpha^\prime|=1}\int_0^t\int
\int\f{\nu(|v|)|\partial^{\alpha^\prime}(\mb{P}_cF_2))|^2}{\mb{M}_*}
dvdxd\tau\\
\di \quad +C(\chi_{\scriptscriptstyle T} +\delta_0)\Big[\int_0^t\int
\int\f{\nu(|v|)|\mb{P}_cF_2|^2}{\mb{M}_*} dvdxd\tau+\int_0^t\int
\int\f{\nu(|v|)|(\widetilde{\mb{G}}, \widetilde{\mb{G}}_{v_1})|^2}{\mb{M}_*}
dvdxd\tau\Big].
\end{array}
\end{equation}
By \eqref{phi-xtt-e-s}, it holds that
\begin{equation}\label{phixt-s}
 \int_0^t\|\Pi_{x\tau}\|^2d\tau \leq
C\int_0^t\Big[\|n_{2}\|^2+\int\int\f{\nu(|v|)|\mb{P}_c
F_2|^2}{\mb{M}_*}dxdv\Big]d\tau.
\end{equation}
By the equation \eqref{n21}, one has
\begin{equation}\label{n2t-E-s}
\begin{array}{ll}
\di \int_0^t\|n_{2\tau}\|^2d\tau=\int_0^t\|\int v_1F_{2x}dv\|^2d\tau=\int_0^t\|\int v_1\big(\f{\mb{M}}{\rho}n_2+\mb{P}_c F_2\big)_xdv\|^2d\tau\\
\di\quad \leq C\int_0^t\Big[\|n_{2x}\|^2+\int\int\f{\nu(|v|)|(\mb{P}_cF_2)_{x}|^2}{\mb{M}_*}dxdv\Big]d\tau+C(\chi_{\scriptscriptstyle T} +\delta_0)\int_0^t\|(n_2,\phi_x,\psi_x,\omega_x)\|^2d\tau.
\end{array}
\end{equation}

On the other hand, from \eqref{Psi}, \eqref{W} and \eqref{tp}, it holds that
\begin{equation*}\label{130}
\begin{array}{ll}
\di\int_0^t\|(\phi,\psi,\o)\|^2d\tau \leq C\int_0^t\|(\phi,\tilde\psi,\tilde\o)\|^2d\tau\\
\di+C\d_0\int_0^t\|\sqrt{|u^{S_1}_{1x}|+|u^{S_3}_{1x}|+|\T_x|}(\Phi,\wt\Psi,\wt W)\|^2d\tau
+\int_0^t\int q|(\Phi,\wt \Psi, \wt W)|^2dxd\tau
\end{array}
\end{equation*}
and
\begin{equation*}
\begin{array}{ll}
\di \int_0^t\|(\phi_x,\psi_x,\o_x)\|^2d\tau\leq C\int_0^t\|(\phi_x,\tilde\psi_x,\tilde\o_x)\|^2d\tau
+C\d_0
\int_0^t\|(\phi,\tilde\psi,\tilde\o)\|^2d\tau\\[3mm]
\di\quad+C\d_0\int_0^t\|\sqrt{|u^{S_1}_{1x}|+|u^{S_3}_{1x}|+|\T_x|}(\Phi,\wt\Psi,\wt W)\|^2d\tau
+\int_0^t\int q|(\Phi,\wt \Psi, \wt W)|^2dxd\tau.
\end{array}
\end{equation*}
On the other hand, from the fluid-type system  (\ref{P-sys}), we can
get \begin{equation*}
\begin{array}{ll}
\di \int_0^t\|(\Phi_{\tau},\wt\Psi_{\tau},\wt W_{\tau})\|^2d\tau\leq C
\int_0^t\|(\Phi_{xx},\wt\Psi_{xx},\wt W_{xx},\Phi_x,\wt \Psi_x, \wt W_x)\|^2d\tau\\[3mm]
\di\quad+C\d_0\int_0^t\|\sqrt{|u^{S_1}_{1x}|+|u^{S_3}_{1x}|+|\T_{x}|}(\Phi,\wt\Psi,\wt W)\|^2d\tau+C(\delta_0+\chi_{\scriptscriptstyle T} )\int_0^t\int\int\f{\nu(| v|) }{\mb{M}_*}|\wt{\mb{G}}|^2d v
dxd\tau\\[3mm]
\di\quad+C(\delta_0+\chi_{\scriptscriptstyle T} )\int_0^t\|\Pi_x\|^2d\tau+C\sum_{|\a^\prime|=1}\int_0^t\int\int\f{\nu(| v|)}{\mb{M}_*}|\partial^{\a^\prime}\wt{\mb{G}}|^2d v
dxd\tau +C\d_0^2.
\end{array}
\label{Phit}
\end{equation*}
In summary, collecting all the above lower order estimates and choosing suitably small $\chi_{\scriptscriptstyle T} $, $\delta$ and $\eta_0$, we complete the proof of Proposition \ref{Prop3.1}.

\

{\bf Appendix A.2. Proof of Proposition \ref{Prop3.2}: Higher order estimates}
\renewcommand{\theequation}{A.2.\arabic{equation}}
\setcounter{equation}{0}

\

\underline{Step 1. Estimation on $\|(\phi,\tilde\psi,\tilde\omega)\|^2$.}

Similar to the lower order estimates \eqref{le1-s}, we multiply $\eqref{sys-h}_1$ by $\di\f{2\phi}{3\tilde\r}$,
$\eqref{sys-h}_2$ by $\di \f{\tilde\psi_1}{\tilde\t}$, $\eqref{sys-h}_3$ by
$\di\tilde\psi_i$, $\eqref{sys-h}_4$ by $\di \f{\tilde\o}{\tilde\t^2}$ respectively
and adding them together to get
\begin{equation*}
\begin{array}{ll}
&\di \|(\phi,\tilde\psi,\tilde\o)(t,\cdot)\|^2+\int_0^t\|(\tilde\psi_x,\tilde\o_x)\|^2d\tau
\leq C(\mathcal{E}(0)^2+\d_0^{\f12})+\int_0^t\int q|(\Phi,\wt\Psi,\wt W)|^2 dxd\tau\\[3mm]
&\di +C\d_0\int_0^t \|\sqrt{|u_{1x}^{S_1}|+|u_{1x}^{S_3}|+|\T_x|}(\Phi,\wt\Psi,\wt W)\|^2 d\tau
+C(\chi_{\scriptscriptstyle T} +\delta_0)\int_0^t\|(\phi_x,\phi_{xx},\tilde\psi_{xx},\tilde\o_{xx})\|^2d\tau\\[3mm]
&\di
+C(\chi_{\scriptscriptstyle T} +\d_0)\int_0^t\|(\phi,\tilde\psi,\tilde\o,\Pi_x,n_2)\|^2d\tau+C\sum_{|\a^\prime|=1}\int_0^t\int\int\f{\nu(| v|)}{\mb{M}_*}|\partial^{\a^\prime}\wt{\mb{G}}|^2d v
dxd\tau
\\
&\di +C\d_0\int_0^t\|(\wt\Psi_{\tau},\wt W_{\tau})\|^2 d\tau+C(\chi_{\scriptscriptstyle T} +\d_0)\int_0^t\int\int\f{\nu(|v|)|(\partial_{v_1}(\mb{P}_c F_2), \wt{\mb{G}})|^2}{\mb{M}_*} dvdxd\tau.
\end{array}
\end{equation*}

\underline{Step 2. Estimation on $\|\phi_x(t,\cdot)\|^2$.}

To estimate the term $\di \int_0^t\|\p_x\|^2d\tau$, we
rewrite the equation $(\ref{sys-h})_2$ as
\begin{equation}
\begin{array}{l}
\di\f{4\mu(\tilde\t)}{3\tilde\r}\phi_{xt}+\tilde\r\tilde\psi_{1t}
+\tilde\r\tilde u_1\tilde\psi_{1x}+\f23\tilde\t\phi_x
=-\f{4\mu(\tilde\t)}{3\tilde\r}\big(2\tilde\r_x\tilde\psi_{1x}+\tilde\r_{xx}\tilde\psi_1
+(\tilde u_1\phi)_{xx}+L_{0x}\big)+\f43\mu'(\t)\tilde\t_x\tilde\psi_{1x}\\[3mm]
\di\qquad+\f13\tilde\r\tilde u_{1x}\tilde\psi_1
-\f23\tilde\r_x\tilde\o-\f23\tilde\r\tilde\o_x+\f{2\tilde\t\tilde\r_x}{3\tilde\r}\phi-\int  v_1^2(\G-\G^{S_1}-\G^{S_3})_x d v
+(J_1+N_1-Q_1)_x-L_1,
\end{array}
\label{(4.36+)}
\end{equation}
Multiplying the equation (\ref{(4.36+)}) by $\phi_x$
and integrating the resulting equation with respect to $t,x$, we get
\begin{equation*}
\begin{array}{ll}
\di \|\p_x(t,\cdot)\|^2+\int_{0}^t\|\p_x\|^2d\tau\leq
C\|\tilde\psi_1(t,\cdot)\|^2+ C\int_0^t\|(\tilde\psi_{1x},\tilde\o_x)\|^2d\tau\\
\di \quad+C(\chi_{\scriptscriptstyle T} +\d_0)\int_0^t\|(\phi_{xx},\tilde\psi_{xx},\tilde\o_{xx})\|^2d\tau
+C(\chi_{\scriptscriptstyle T} +\d_0)\int_0^t\|(\Pi_x,n_2)\|^2d\tau\\
\di \quad +C\d_0\int_0^t\|\sqrt{|u_{1x}^{S_1}|+|u_{1x}^{S_3}|+|\T_{x}|}(\Phi,\wt\Psi,\wt W)\|^2d\tau
+C(\mathcal{E}(0)^2+\d_0^\f12)+C\d_0\int_0^t\|(\phi,\tilde\psi,\tilde\o)\|^2 d\tau\\
\di\quad +C(\chi_{\scriptscriptstyle T} +\d_0)\sum_{|\a^\prime|=1}\int_0^t\int
\int\f{\nu(| v|)|(\partial^{\a^\prime}\wt{\mb{G}},\wt{\mb{G}},\partial_{v_1}(\mb{P}_c F_2),\partial_{v_1x}(\mb{P}_c F_2))|^2}{\mb{M}_*}
d v dxd\tau\\
\di\quad + C\sum_{|\a|=2}\int_{0}^t \int\int\f{\nu^{-1}(| v|)}{\mb{M}_*}|\partial^\a \mb{G} |^2d v dxd\tau
+\int_0^t\int q|(\Phi, \wt\Psi, \wt W)|^2 dxd\tau.
\end{array}
\end{equation*}
Now we turn to the time-derivative terms. To estimate
$\|(\p_{t},\psi_{t},\o_{t})\|^2$, we need to use the system
(\ref{sys-h-o-s}). By multiplying $(\ref{sys-h-o-s})_1$ by $\p_{t}$,
$(\ref{sys-h-o})_2$ by $\psi_{1t}$, $(\ref{sys-h-o-s})_3$ by
$\psi_{it}~(i=2,3)$ and $(\ref{sys-h-o-s})_4$ by $\o_{t}$
respectively, and adding them together, after integrating with
respect to $t$ and $x$, we have
\begin{equation*}
\begin{array}{ll}
\di \int_0^t\|(\p_{\tau},\psi_{\tau},\o_{\tau})(\tau,\cdot)\|^2d\tau\leq C\d_0\int_0^t\|\sqrt{|u_{1x}^{S_1}|+|u_{1x}^{S_3}|+|\T_{x}|}(\Phi,\wt\Psi,\wt W)\|^2 d\tau\\
\di\qquad +C\int_0^t \|(\phi_x,\tilde\psi_x,\tilde\o_x)\|^2d\tau
++C\d_0\int_0^t\|(\phi,\tilde\psi,\tilde\o)\|^2d\tau
+C\int_0^t\int\int\f{\nu(| v|)}{\mb{M}_*}|\wt{\mb{G}}_x|^2d v
dxd\tau\\
\di\qquad +\int_0^t\int q|(\Phi, \wt\Psi, \wt W)|^2 dxd\tau+C(\chi_{\scriptscriptstyle T} +\d_0)\int_0^t\|(\Pi_x,n_2)\|^2d\tau+C\d_0.
\end{array}
%\label{(4.45)}
\end{equation*}

\underline{Step 3. Estimation on
$\|(\phi_x,\tilde\psi_x,\tilde\o_x)(t,\cdot)\|^2$.}

Multiplying $\eqref{sys-h-h-s}_1$ by $\di\f{2\phi_x}{3\tilde\r}$, $\eqref{sys-h-h-s}_2$
by $\di\f{\tilde\psi_{1x}}{\tilde\t}$, $\eqref{sys-h-h-s}_3$ by $\tilde\psi_{ix}$ and
$\eqref{sys-h-h-s}_4$ by $\di\f{\tilde\o_x}{\tilde\t^2}$, adding them together,
and  integrating the resulting equation with respect to $t,x$, we have
\begin{equation}
\begin{array}{l}
\di \|(\phi_x,\tilde\psi_{x},\tilde\o_x)(t,\cdot)\|^2+\int_0^t\|(\tilde\psi_{xx},\tilde\o_{xx})\|^2d\tau
\leq C(\mathcal{E}(0)^2+\d_0)+C(\chi_{\scriptscriptstyle T} +\d_0)\int_0^t\|(\Pi_x,n_2)\|^2d\tau\\
\di\quad + C(\chi_{\scriptscriptstyle T} +\d_0)\int_0^t\|(\phi,\tilde\psi,\tilde\o)\|_{H^1}^2d\tau
+C\d_0\int_0^t\|\sqrt{|u_{1x}^{S_1}|+|u_{1x}^{S_3}|+|\T_x|}(\Phi,\wt\Psi,\wt W)\|^2 d\tau\\
\di\quad +C\sum_{|\a|=2}\int_0^t\int\int\f{\nu(| v|)}{\mb{M}_*}|\partial^\a
\wt{\mb{G}} |^2d v
dxd\tau
+\int_0^t\int q|(\Phi, \wt\Psi, \wt W)|^2 dxd\tau\\
\di\quad +C(\chi_{\scriptscriptstyle T} +\d_0)\sum_{|\a^\prime|=1}\int_0^t\int
\int\f{\nu(| v|)|(\partial^{\a^\prime}\wt{\mb{G}},\wt{\mb{G}},\partial_{v_1}(\mb{P}_c F_2),\partial_{v_1x}(\mb{P}_c F_2))|^2}{\mb{M}_*}
d v dxd\tau.
\end{array}
\label{(4.49)}
\end{equation}
To get  the estimation of $\|\p_{xx}\|^2$ , we apply
$\partial_x$ to $(\ref{sys-h-o-s})_2$, we get
\begin{equation}
\begin{array}{ll}
\di \psi_{1xt}+(\r u_1^2-\tilde\r\tilde u_1^2+p-\tilde p)_{xx}-(\f{\Pi_x^2}{4})_{xx}
=-\f{4}{3}\Big[\mu(\tilde\t)\tilde u_{1x}-\mu(\t^{S_1})u_{1x}^{S_1}-\mu(\t^{S_3})u_{1x}^{S_3} \Big]_{xx}\\
\di \qquad -\int v_1^2\wt{\mb{G}}_{xx}d v- Q_{1xx}. \label{(4.50)}
\end{array}
\end{equation}
Multiplying (\ref{(4.50)}) by $\p_{xx}$ and integrating the resulting equations,   we obtain
\begin{equation}
\begin{array}{l}
\di \int\psi_{1x}\p_{xx}(t,x)dx+\int_0^t\|\phi_{xx}\|^2d\tau\leq C\|\psi_{1xx}\|^2
+C(\chi_{\scriptscriptstyle T} +\d_0)\int_0^t\|(\p,\tilde\psi,\tilde\o)\|_{H^1}^2d\tau+C(\mathcal{E}(0)^2+\d_0^\f12)\\[3mm]
\di\quad
+C\int_0^t\|(\tilde\psi_{xx},\tilde\o_{xx})\|^2d\tau+C\d_0\int_0^t\|\sqrt{|u_{1x}^{S_1}|+|u_{1x}^{S_3}|+|\T_x|}(\Phi,\wt\Psi,\wt W)\|^2 d\tau \\[3mm]
\di\quad+C(\chi_{\scriptscriptstyle T} +\d_0)\int_0^t\|(n_2,n_{2x})\|^2d\tau
+C\sum_{|\a|=2}\int_0^t\int\int\f{\nu(| v|)}{\mb{M}_*}|\partial^\a
\wt{\mb{G}} |^2d v dxd\tau+\int_0^t\int q|(\Phi, \wt\Psi, \wt W)|^2 dxd\tau.
\end{array}
\label{(4.52)}
\end{equation}

To estimate $\|(\p_{xt},\psi_{xt},\o_{xt})\|^2$ and
$\|(\p_{tt},\psi_{tt},\o_{tt})\|^2$, we use the
system \eqref{sys-h-o-s} again. Applying $\partial_x$ first, and multiplying the four equations of \eqref{sys-h-o-s}
by $\p_{xt}$, $\psi_{1xt}$, $\psi_{ixt}$ $(i= 2,3)$,
$\o_{xt}$ respectively, then adding them together and integrating
with respect to $t$ and $x$, we have
\begin{equation*}
\begin{array}{ll}
\di \int_0^t\|(\p_{x\tau},\psi_{x\tau},\o_{x\tau})\|^2d\tau \leq C\int_0^t\|(\p_{xx},\tilde\psi_{xx},\tilde\o_{xx})\|^2d\tau
 +C(\chi_{\scriptscriptstyle T} +\d_0)\int_0^t\|(\phi,\tilde\psi,\tilde\o)\|_{H^1}^2 d\tau\\
\quad\di+C\d_0\int_0^t\|\sqrt{|u_{1x}^{S_1}|+|u_{1x}^{S_3}|+|\T_x|}(\Phi,\wt\Psi,\wt W)\|^2d\tau
+\int_0^t\int q|(\Phi, \wt\Psi, \wt W)|^2 dxd\tau\\
 \quad\di +C\sum_{|\a|=2}\int_0^t\int\int\f{\nu(| v|)}{\mb{M}_*}|\partial^\a
 \wt{\mb{G}}|^2d v dxd\tau+C(\chi_{\scriptscriptstyle T} +\d_0)\int_0^t\|(\Pi_{x\tau},n_{2\tau})\|^2d\tau +C\d_0^\f12.
\end{array}
\label{(4.53)}
\end{equation*}
Similarly, we can obtain
\begin{equation}
\begin{array}{l}
\di \int_0^t\|(\p_{\tau\tau},\psi_{\tau\tau},\o_{\tau\tau})\|^2d\tau\leq C\int_0^t\|(\p_{x\tau},\psi_{x\tau},\o_{x\tau})\|^2d\tau
 +C(\chi_{\scriptscriptstyle T} +\d_0)\sum_{|\a'|=1}\int_0^t\|\partial^{\a'}(\phi,\psi,\o)\|^2 d\tau\\
\quad\di +C(\chi_{\scriptscriptstyle T} +\d_0)\int_0^t\|(\phi,\psi,\o)\|^2 d\tau
+C\sum_{|\a|=2}\int_0^t\int\int\f{\nu(| v|)}{\mb{M}_*}|\partial^\a
 \wt{\mb{G}}|^2d v dxd\tau+C\d_0^\f12.
\end{array}
\label{(4.54)}
\end{equation}
A suitable linear combination of \eqref{(4.49)}, \eqref{(4.52)}-\eqref{(4.54)} gives
\begin{equation}
\begin{array}{l}
\di\|(\phi_x,\tilde\psi_{x},\tilde\o_x)(t,\cdot)\|^2
+\int_0^t\Big[\|(\phi_{xx},\tilde\psi_{xx},\tilde\o_{xx})\|^2
+\sum_{|\a'|=1}\|\partial^{\a'}(\phi_{\tau},\psi_{\tau},\o_{\tau})\|^2 \Big] d\tau\\[3mm]
\di\quad\le C\|(\phi_{xx},\psi_{1xx})\|^2
+ C(\chi_{\scriptscriptstyle T} +\d_0)\int_0^t\big[\|(\phi,\tilde\psi,\tilde\o)\|_{H^1}^2
+\|(\phi,\psi,\o)_{\tau}\|^2 \big] d\tau
\\[3mm]
\di\quad
+C\d_0\int_0^t \|\sqrt{|u_{1x}^{S_1}|+|u_{1x}^{S_3}|+|\T_x|}(\Phi,\wt\Psi,\wt W)\|^2d\tau
 +C\sum_{|\a|=2}\int_0^t\int\int\f{\nu(| v|)}{\mb{M}_*}|\partial^\a
\wt{\mb{G}}|^2d v dxd\tau\\
\di\quad +C(\chi_{\scriptscriptstyle T} +\d_0)\sum_{|\a^\prime|=1}\int_0^t\int
\int\f{\nu(| v|)|(\partial^{\a^\prime}\wt{\mb{G}},\wt{\mb{G}},\partial_{v_1}(\mb{P}_c F_2),\partial_{v_1x}(\mb{P}_c F_2))|^2}{\mb{M}_*}
d v dxd\tau\\
\di\quad +\int_0^t\int q|(\Phi, \wt\Psi, \wt W)|^2 dxd\tau+C(\chi_{\scriptscriptstyle T} +\d_0)\int_0^t\|(\Pi_x,\Pi_{x\tau},n_2,n_{2x},n_{2\tau})\|^2d\tau+C(\mathcal{E}(0)^2+\d_0^\f12).
\end{array}
\label{(4.55)}
\end{equation}

\underline{Step 4. Estimation on the non-fluid component.}

To close the above estimate, we need to estimate the derivatives
on the non-fluid component\\ $\partial^\alpha\partial^\b(\wt{\mathbf{G}}, \mb{P}_cF_2)(1\leq
|\alpha|\leq 2, 0\leq|\beta|\leq 2)$. First, applying $\partial_x$ to \eqref{F1-G} and \eqref{F2-pc-1}, respectively, we have
\begin{equation}\label{Gx-s}
\begin{array}{ll}
\di\wt{\mb{G}}_{xt}-(\mb{L}_\mb{M}\wt{\mb{G}})_x =\Big\{
-\mb{P}_1( v_1\wt{\mb{G}}_x)-\mathbf{P}_1( \Pi_xF_{2v_1})
+2[Q(\wt{\mb{G}},\mb{G}^{S_1}+\mb{G}^{S_3})+Q(\mb{G}^{S_1}+\mb{G}^{S_3},\wt{\mb{G}})]\\[3mm]
\di\qquad\qquad\qquad +2Q(\wt{\mb{G}},\wt{\mb{G}})+2[Q(\mb{G}^{S_1},\mb{G}^{S_3})+Q(\mb{G}^{S_3},\mb{G}^{S_1})]\\
\di\qquad\qquad\qquad -\Big[\mb{P}_1( v_1\mb{M}_x)-\mb{P}^{S_1}_1( v_1\mb{M}^{S_1}_x)
-\mb{P}^{S_3}_1( v_1\mb{M}^{S_3}_x)\Big]+\sum_{j=1,3}R_j
\Big\}_x.
\end{array}
\end{equation}
and
\begin{equation}\label{Pc-F2x-s}
\begin{array}{ll}
\di
(\mb{P}_c F_2)_{xt}+\big[\mb{P}_c(v_1 F_{2x})\big]_x+\big[\mb{P}_c (\Pi_x\partial_{v_1}F_1)\big]_x+\big[(\f{\mb{M}}{\rho})_t~n_2\big]_x\\
\di\qquad=\mb{N}_{\mb{M}}(\mb{P}_c F_2)_x+2Q(\mb{P}_c F_2,\mathbf{M}_x)+2Q(F_{2x},\mb{G})+2Q(F_{2},\mb{G}_x).
\end{array}
\end{equation}
Multiplying \eqref{Gx-s} and \eqref{Pc-F2x-s} by
$\frac{\wt{\mathbf{G}}_x}{\mathbf{M}_*}$ and
$\frac{(\mathbf{P}_cF_2)_x}{\mathbf{M}_*}$, respectively, and
integrating with respect to $x, v$ and $ t$, and also using Lemmas
\ref{Lemma 4.1}-\ref{Lemma 4.3}, we obtain
\begin{equation}\label{Gx-E}
\begin{array}{l}
\displaystyle
\int\int\frac{|(\wt{\mathbf{G}}, \mb{P}_cF_2)_x|^2}{2\mathbf{M}_*}(x,v,t) d v dx+\|n_{2}(\cdot,t)\|^2
+\frac{\widetilde{\sigma}}{2}\int_0^t\int\int\frac{\nu(| v|)}{\mathbf{M}_*}|(\wt{\mathbf{G}}, \mb{P}_cF_2)_x|^2d v dxd\tau\\
\leq\displaystyle C\int\int\frac{|(\wt{\mathbf{G}},
\mb{P}_cF_2)_x|^2}{2\mathbf{M}_*}(x,v,0) d v
dx+C\|n_{20}\|^2+C\int_0^t\int\int\frac{\nu(| v|)|(\wt{\mathbf{G}},
\mb{P}_cF_2)_{xx}|^2}{\mathbf{M}_*}d v
dxd\tau\\
\di +C(\chi_{\scriptscriptstyle T} +\delta_0)\sum_{|\alpha^\prime|=1}\int_0^t
\|\partial^{\alpha^\prime}(\phi,\psi,\omega)\|^2d\tau
+C\int_0^t\|(\phi_{xx},\psi_{xx},\omega_{xx},n_{2x},n_{2xx})\|^2 d\tau\\
\di +C(\chi_{\scriptscriptstyle T} +\delta_0)\int_0^t \|(\phi_{x\tau},\psi_{x\tau},\omega_{x\tau})\|^2
d\tau+C\delta_0+C(\chi_{\scriptscriptstyle T} +\delta_0+\eta_0)\int_0^t\|n_2\|^2d\tau
\\
\di
+C(\chi_{\scriptscriptstyle T} +\delta_0)\sum_{0\leq|\beta^\prime|\leq1}\int_0^t\int\int\frac{\nu(|
v|)|\partial^{\beta^\prime}(\widetilde{\mb{G}},
\mb{P}_cF_2)|^2}{\mathbf{M}_*}d v dxd\tau.
\end{array}
\end{equation}
Note that in the above inequality, we have used the fact that
$$
\begin{array}{ll}
\di -\int_0^t\int\int \Pi_{xx}\mb{M}_{v_1}\f{(\mb{P}_c
F_2)_x}{\mb{M}_*}dxdvd\tau=2\int_0^t\int
n_2\mb{M}\f{v_1-u_1}{R\t}\f{(\mb{P}_c F_2)_x}{\mb{M}_*}dxd\tau\\
\di =2\int_0^t\int\int
n_2\f{v_1-u_1}{R\t}\big(\f{\mb{M}}{\mb{M}_*}-1\big)(\mb{P}_c
F_2)_xdxdvd\tau+2\int_0^t\int\int n_2\f{v_1-u_1}{R\t}(\mb{P}_c
F_2)_xdxdvd\tau\\
\di =2\int_0^t\int\int
n_2\f{v_1-u_1}{R\t}\big(\f{\mb{M}}{\mb{M}_*}-1\big)(\mb{P}_c
F_2)_xdxdvd\tau+2\int_0^t\int \f{n_2}{R\t}\big(\int v_1\mb{P}_c
F_2dv\big)_xdxd\tau\\
\di =2\int_0^t\int\int
n_2\f{v_1-u_1}{R\t}\big(\f{\mb{M}}{\mb{M}_*}-1\big)(\mb{P}_c
F_2)_xdxdvd\tau-2\int_0^t\int
\f{n_2}{R\t}\big[n_{2\tau}+(u_1n_2)_x\big]dxd\tau\\
\di\leq -\int\f{n_2^2}{R\t}(x,t)
dx+C\|n_{20}\|^2+C(\chi_{\scriptscriptstyle T} +\delta_0+\eta_0)\int_0^t\Big[\|n_2\|^2+\int\int
\f{\nu(|v|)|(\mb{P}_c
F_2)_x|^2}{\mb{M}_*}dxdv\Big]d\tau\\
\di \qquad +C(\chi_{\scriptscriptstyle T} +\delta_0)\int_0^t \|(\omega_\tau,\psi_{1x},\omega_x)\|^2
d\tau.
\end{array}
$$
Similarly, one has
$$
\begin{array}{ll}
\di -\int_0^t\int\int \Pi_{xt}\mb{M}_{v_1}\f{(\mb{P}_c
F_2)_\tau}{\mb{M}_*}dxdvd\tau=\int_0^t\int
\Pi_{xt}\mb{M}\f{v_1-u_1}{R\t}\f{(\mb{P}_c F_2)_x}{\mb{M}_*}dxd\tau\\
\di =\int_0^t\int\int
\Pi_{x\tau}\f{v_1-u_1}{R\t}\big(\f{\mb{M}}{\mb{M}_*}-1\big)(\mb{P}_c
F_2)_\tau dxdvd\tau+\int_0^t\int \f{\Pi_{x\tau}}{R\t}\big(\int v_1\mb{P}_c
F_2dv\big)_\tau dxd\tau\\
\di =\int_0^t\int\int \Pi_{x\tau}
\f{v_1-u_1}{R\t}\big(\f{\mb{M}}{\mb{M}_*}-1\big)(\mb{P}_c
F_2)_xdxdvd\tau-\int_0^t\int
\f{\Pi_{x\tau}}{R\t}\big[\f12\Pi_{x\tau\tau}+(u_1n_2)_\tau\big]dxd\tau\\
\di\leq -\int\f{\Pi_{x\tau}^2}{4R\t}(x,t)
dx+C\|\Pi_{0xt}\|^2+C(\chi_{\scriptscriptstyle T} +\delta_0+\eta_0)\int_0^t\Big[\|\Pi_{x\tau}\|^2+\int\int
\f{\nu(|v|)|(\mb{P}_c
F_2)_\tau|^2}{\mb{M}_*}dxdv\Big]d\tau\\
\di \qquad +C(\chi_{\scriptscriptstyle T} +\delta_0)\int_0^t \|(\omega_\tau,\psi_{1x},\omega_x)\|^2
d\tau,
\end{array}
$$
and further
\begin{equation}\label{Gt-E}
\begin{array}{l}
\displaystyle
\int\int\frac{|(\wt{\mathbf{G}}, \mb{P}_cF_2)_t|^2}{2\mathbf{M}_*}(x,v,t)d v dx+\|\Pi_{xt}(\cdot,t)\|^2+\frac{\widetilde{\sigma}}{2}\int_0^t\int\int\frac{\nu(| v|)}{\mathbf{M}_*}|(\wt{\mathbf{G}}, \mb{P}_cF_2)_\tau|^2d v dxd\tau\\
\leq\displaystyle C\int\int\frac{|(\wt{\mathbf{G}},
\mb{P}_cF_2)_t|^2}{2\mathbf{M}_*}(x,v,0)d v
dx+C\|\Pi_{0xt}\|^2+C\int_0^t\int\int\frac{\nu(| v|)|(\wt{\mathbf{G}},
\mb{P}_cF_2)_{x\tau}|^2}{\mathbf{M}_*}d v
dxd\tau\\
\di +C(\chi_{\scriptscriptstyle T} +\delta_0)\int_0^t \|(\phi_x,\varphi_{\tau\tau},\psi_{\tau\tau},\omega_{\tau\tau})\|^2 d\tau+C\delta_0+C\int_0^t\|(\phi_{x\tau},\psi_{x\tau},\omega_{x\tau},n_{2\tau},n_{2x\tau})\|^2 d\tau\\
\di +C(\chi_{\scriptscriptstyle T} +\delta_0)\sum_{|\alpha^\prime|=1}\int_0^t
\|\partial^{\alpha^\prime}(\phi,\psi,\omega)\|^2d\tau
+C(\chi_{\scriptscriptstyle T} +\delta_0+\eta_0)\int_0^t\|\Pi_{x\tau}\|^2d\tau \\
\di
+C(\chi_{\scriptscriptstyle T} +\delta_0)\sum_{0\leq|\beta^\prime|\leq1}\int_0^t\int\int\frac{\nu(|
v|)|\big(\partial^{\beta^\prime}(\widetilde{\mb{G}}, \mb{P}_cF_2),
(\wt{\mathbf{G}},\mb{P}_cF_2)_x\big)|^2}{\mathbf{M}_*}d v dxd\tau.
\end{array}
\end{equation}
The combination of  \eqref{Gx-E} and \eqref{Gt-E}, and choosing $\chi_{\scriptscriptstyle T} ,\delta$ suitably small  yield that
\begin{equation*}\label{xt-e}
\begin{array}{ll}
\displaystyle
\sum_{|\alpha|=1}\int\int\frac{|\partial^\alpha(\wt{\mathbf{G}}, \mb{P}_cF_2)|^2}{\mathbf{M}_*}(x,v,t)d v dx+\sum_{|\alpha|=1}\int_0^t\int\int\frac{\nu(| v|)|\partial^\alpha(\wt{\mathbf{G}}, \mb{P}_cF_2)|^2}{\mathbf{M}_*}d v dxd\tau\\
\di +\|(n_2,\Pi_{xt})(\cdot,t)\|^2 \leq C\int |\phi_{xx}\psi_{1x}(x,t)|dx+C\sum_{|\alpha|=2}\int_0^t\int\int\frac{\nu(| v|)|\partial^\alpha(\wt{\mathbf{G}},
\mb{P}_cF_2)|^2}{\mathbf{M}_*}d v
dxd\tau\\
\di +C(\chi_{\scriptscriptstyle T} +\sqrt\delta_0)\int_0^t\Big[\|\sqrt{\bar u_{1x}}(\phi,\psi_1,\omega)\|^2+
\sum_{|\alpha^\prime|=1}\|\partial^{\alpha^\prime}(\phi,\psi,\omega)\|^2+\|\Pi_x\|^2\Big]d\tau\\
\di

+C(\chi_{\scriptscriptstyle T} +\delta_0)\sum_{0\leq|\beta^\prime|\leq1}\int_0^t\int\int\frac{\nu(|
v|)|\partial^{\beta^\prime}(\widetilde{\mb{G}}, \mb{P}_cF_2\big)|^2}{\mathbf{M}_*}d v dxd\tau+C(\mathcal{E}(0)^2+\delta_0)\\
\di+C(\chi_{\scriptscriptstyle T} +\delta_0)\int_0^t\int\int\frac{\nu(|
v|)|(\wt{\mb{G}}, \mb{P}_cF_2)_{vx}|^2}{\mathbf{M}_*}d v dxd\tau+C(\chi_{\scriptscriptstyle T} +\sqrt\delta_0+\eta_0)\int_0^t\|(n_2,\Pi_{x\tau})\|^2d\tau.
\end{array}
\end{equation*}
Applying $\partial_{v_j}~(j=1,2,3)$ to the equation \eqref{G}  and noting that \eqref{P1-1-s}, one has
\begin{equation}\label{Gv}
\begin{array}{ll}
\di
\widetilde{\mb{G}}_{v_jt}-\mb{L}_{\mb{M}}\widetilde{\mb{G}}_{v_j}=-\partial_{v_j}\mb{P}_1(v_1\wt{\mb{G}}_x)- \Pi_x\big(\mb{P}_cF_2\big)_{v_1v_j}+\Pi_x\sum_{k=0}^4\int(\mb{P}_cF_2) \big(\f{\chi_k}{\mb{M}}\big)_{v_1}dv(\chi_k)_{v_j}\\
\di \quad+(\partial_{v_j}\mb{L}_{\mb{M}}\widetilde{\mb{G}}-\mb{L}_{\mb{M}}\widetilde{\mb{G}}_{v_j}) +2\partial_{v_j}Q(\wt{\mb{G}},\wt{\mb{G}})+2\partial_{v_j}[Q(\wt{\mb{G}},\mb{G}^{S_1}+\mb{G}^{S_3})+Q(\mb{G}^{S_1}+\mb{G}^{S_3},\wt{\mb{G}})]\\[3mm]
\di\quad +2\partial_{v_j}[Q(\mb{G}^{S_1},\mb{G}^{S_3})+Q(\mb{G}^{S_3},\mb{G}^{S_1})]-\Big[\mb{P}_1( v_1\mb{M}_x)-\mb{P}^{S_1}_1( v_1\mb{M}^{S_1}_x)
-\mb{P}^{S_3}_1( v_1\mb{M}^{S_3}_x)\Big]_{v_j}+\partial_{v_j}\sum_{k=1,3}R_k.
\end{array}
\end{equation}
Similarly, applying $\partial_{v_j}~(j=1,2,3)$ to the equation \eqref{F2-pc-1}  and noting that \eqref{Pc-1-s}  imply
that
\begin{equation}\label{F2-pc-v}
\begin{array}{ll}
\di \partial_t(\mb{P}_c F_2)_{v_j}-\mb{N}_{\mb{M}}(\mb{P}_c F_2)_{v_j}+\partial_{v_j}\mb{P}_c(v_1 F_{2x})+\Pi_x\mb{M}_{v_1v_j}+\Pi_x\widetilde{\mb{G}}_{v_1v_j}\\
\di \qquad+\Pi_x(\mb{G}^{S_1}_{v_1v_j}+\mb{G}^{S_3}_{v_1v_j})+(\f{\mb{M}}{\rho})_{tv_j}~n_2 =2Q(\mb{P}_cF_2,\mb{M}_{v_j})+2\partial_{v_j}Q(F_2,\mb{G}).
\end{array}
\end{equation}
Recalling the following two facts
\begin{equation*}\label{Fact1}
\partial_{v_j}\mb{L}_{\mb{M}}g-\mb{L}_{\mb{M}}(g_{v_j})=2Q(\partial_{v_j}\mb{M},g)+2Q(g,\partial_{v_j}\mb{M}),
\end{equation*}
\begin{equation*}\label{Fact2}
\big(\mb{L}_{\mb{M}}^{-1}g\big)_{v_j}=\mb{L}_{\mb{M}}^{-1}(g_{v_j})+\mb{L}_{\mb{M}}^{-1}\big(\partial_{v_j}\mb{L}_{\mb{M}}g-\mb{L}_{\mb{M}}(g_{v_j})\big).
\end{equation*}
Multiplying the equations \eqref{Gv} and \eqref{F2-pc-v} by
$\f{\widetilde{\mb{G}}_{v_j}}{\mb{M}_*}$ and $\f{(\mb{P}_cF_2)_{v_j}}{\mb{M}_*}$, respectively, and then integrating with
respect to $x,v,t$ imply that
\begin{equation}\label{GvE1}
\begin{array}{ll}
\di \int\int \f{|(\widetilde{\mb{G}}, \mb{P}_c F_2)_{v_j}|^2}{\mb{M_*}}(x,v,t)
dxdv+\int_0^t\int\int\f{\nu(|v|)|(\widetilde{\mb{G}}, \mb{P}_c F_2)_{v_j}|^2}{\mb{M_*}}
dxdvd\tau
\\[3mm]
 \di
\leq C\int\int \f{|(\widetilde{\mb{G}}, \mb{P}_c
F_2)_{v_j}|^2}{\mb{M_*}}(x,v,0)
dxdv+C(\chi_{\scriptscriptstyle T} +\delta_0)
\int_0^t\Big[\|n_2\|^2+\sum_{|\alpha|=1}\|\partial^\alpha(\phi,\psi,\omega)\|^2\Big]
d\tau\\
\di +C\int_0^t\|(\psi_x,\omega_x,\Pi_x,
n_{2x})\|^2d\tau+C\delta_0+C\int_0^t\int\int\f{\nu(|v|)|(\widetilde{\mb{G}},\mb{P}_c F_2,
\wt{\mb{G}}_x,(\mb{P}_c F_2)_x)|^2}{\mb{M_*}} dxdvd\tau.
\end{array}
\end{equation}

Applying $\partial_{v_jv_k}~(j,k=1,2,3)$ to the equations \eqref{G} and \eqref{F2-pc} and using the similar analysis as in obtaining \eqref{GvE1}, one has for $|\beta|=2$
\begin{equation*}\label{GvE2}
\begin{array}{ll}
\di \int\int \f{|\partial^\beta(\widetilde{\mb{G}}, \mb{P}_c F_2)|^2}{\mb{M_*}}(x,v,t)
dxdv+\int_0^t\int\int\f{\nu(|v|)|\partial^\beta(\widetilde{\mb{G}}, \mb{P}_c F_2)|^2}{\mb{M_*}}
dxdvd\tau
\\[3mm]
 \di
\leq C\int\int \f{|\partial^\beta(\widetilde{\mb{G}}, \mb{P}_c F_2)|^2}{\mb{M_*}}(x,v,0)
dxdv+C\int_0^t\|(\psi_x,\omega_x,\Pi_x,n_{2x})\|^2d\tau\\
\di +C\sum_{0\leq|\beta^\prime|\leq 1}\int_0^t\int\int\f{\nu(|v|)|\partial^{\beta^\prime}(\widetilde{\mb{G}},\mb{P}_cF_2, \wt{\mb{G}}_x, (\mb{P}_cF_2)_x)|^2}{\mb{M_*}}
dxdvd\tau\\
\di
+C(\chi_{\scriptscriptstyle T} +\delta_0)
\sum_{|\alpha|=1}\int_0^t\int|\partial^\alpha(\phi,\psi,\omega)|^2
dxd\tau+C\delta_0.
\end{array}
\end{equation*}
By \eqref{le1-s}, \eqref{Phi-E-s}, \eqref{n2-E-s}, it holds that
\begin{equation*}\label{A1}
\begin{array}{ll}
\di \int_0^t\|(\psi_x,\omega_x,\Pi_x,n_{2x})\|^2d\tau\leq   C\|(\phi,\psi,\omega,\phi_x,\Pi_x,n_{2})(\cdot,0)\|^2+C \d_0^{\f12}
+ C(\chi_{\scriptscriptstyle T} +\sqrt\delta_0) \int_0^t\|(\psi_{1xx},\phi_{xx})\|^2d\tau\\
\di
+C(\chi_{\scriptscriptstyle T} +\delta_0)\int_0^t\int\int\f{\nu(|v|)|(\widetilde{\mb{G}},\mb{P}_c F_2,\wt{\mb{G}}_{v_1})|^2}{\mb{M}_*}
dv dxd\tau+C\sum_{1\leq|\alpha|\leq2}
\int_0^t\int\int\f{\nu(|v|)|\partial^\alpha(\wt{\mb{G}},\mb{P}_c F_2)|^2}{\mb{M}_*}
dv
dxd\tau.
\end{array}
\end{equation*}
Similarly, one has for $|\alpha^\prime|=1$ and $|\beta^\prime|=1$
\begin{equation*}\label{GvE11}
\begin{array}{ll}
\di \int\int \f{|\partial^{\alpha^\prime}\partial^{\beta^\prime}(\wt{\mb{G}},\mb{P}_c F_2)|^2}{\mb{M_*}}(x,v,t)
dxdv+\int_0^t\int\int\f{\nu(|v|)|\partial^{\alpha^\prime}\partial^{\beta^\prime}(\wt{\mb{G}},\mb{P}_c F_2)|^2}{\mb{M_*}}
dxdvd\tau
\\[3mm]
 \di
\leq C\int\int
\f{|\partial^{\alpha^\prime}\partial^{\beta^\prime}(\wt{\mb{G}},\mb{P}_c
F_2)|^2}{\mb{M_*}}(x,v,0)
dxdv+C\sum_{1\leq|\alpha|\leq2}
\int_0^t\int\int\f{\nu(|v|)|\partial^\alpha(\wt{\mb{G}},\mb{P}_c
F_2)|^2}{\mb{M_*}} dxdvd\tau\\
\di+C\int_0^t\Big[\sum_{|\alpha|=2}\|\partial^\alpha(\varphi,\psi,\omega,n_2)\|^2+\|(n_2,n_{2x},n_{2\tau},\Pi_{x\tau})\|^2\Big]d\tau+C(\chi_{\scriptscriptstyle T} +\delta_0)
\sum_{|\alpha|=1}\int_0^t\|\partial^\alpha(\phi,\psi,\omega)\|^2
d\tau\\
\di +C(\chi_{\scriptscriptstyle T} +\delta_0)\sum_{0\leq|\beta|\leq
1}\int_0^t\int\int\f{\nu(|v|)|\partial^\beta(\widetilde{\mb{G}},
\mb{P}_c F_2)|^2}{\mb{M_*}}
dxdvd\tau
+C\delta_0.
\end{array}
\end{equation*}
By \eqref{n2-E-s}, \eqref{phixt-s}, \eqref{n2t-E-s} and \eqref{n2-he}, it holds that
\begin{equation*}\label{A2}
\begin{array}{ll}
\di \int_0^t\Big[\sum_{|\alpha|=2}\|\partial^\alpha(\varphi,\psi,\omega,n_2)\|^2+
\|(n_2,n_{2x},n_{2\tau},\Pi_{x\tau})\|^2\Big]d\tau\leq C\|(\Pi_{x0},\psi_{0x},\omega_{0x},\phi_{0xx},n_{20},n_{20x})\|^2\\
\di
+C\int|\varphi_{xx}\psi_{1x}(x,t)|dx+C~\sum_{1\leq|\alpha|\leq2}\int_0^t\int
\int\f{\nu(|v|)|\partial^{\alpha}(\wt{\mb{G}},\mb{P}_cF_2))|^2}{\mb{M}_*}
dvdxd\tau\\
\di+C(\chi_{\scriptscriptstyle T} +\sqrt\delta_0)\int_0^t\Big[\sum_{|\alpha|=1}\|\partial^\alpha(\phi,\psi,\omega)\|^2+\|\sqrt{\bar u_{1x}}(\phi,\psi_1,\omega)\|^2+\|\Pi_x\|^2\Big]d\tau\\
\di +C\int_0^t\int
\int\f{\nu(|v|)|\mb{P}_cF_2|^2}{\mb{M}_*} dvdxd\tau +C(\chi_{\scriptscriptstyle T} +\sqrt\delta_0)\int_0^t\int
\int\f{\nu(|v|)|(\widetilde{\mb{G}}, \widetilde{\mb{G}}_{v_1},\wt{\mb{G}}_{vx})|^2}{\mb{M}_*}
dvdxd\tau.
\end{array}
\end{equation*}
In summary, one has
\begin{equation*}\label{G3}
\begin{array}{ll}
\di \int\int \Big[\sum_{|\alpha|=1, 0\leq |\beta|\leq1} \f{|\partial^{\alpha}\partial^{\beta}(\wt{\mb{G}},\mb{P}_c F_2)|^2}{\mb{M_*}}(x,v,t) +\sum_{1\leq |\beta|\leq2} \f{|\partial^{\beta}(\wt{\mb{G}},\mb{P}_c F_2)|^2}{\mb{M_*}}(x,v,t)\Big]
dxdv
\\
\di + \int_0^t\int\int\Big[\sum_{|\alpha|=1, 0\leq |\beta|\leq1}\f{\nu(|v|)|\partial^{\alpha}\partial^{\beta}(\wt{\mb{G}},\mb{P}_c F_2)|^2}{\mb{M_*}}
+\sum_{1\leq |\beta|\leq2} \f{\nu(|v|)|\partial^{\beta}(\wt{\mb{G}},\mb{P}_c F_2)|^2}{\mb{M_*}}
\Big] dxdvd\tau\\
\di \leq \mathcal{E}^2(0)+C\int |\phi_{xx}\psi_{1x}(x,t)|dx+C\d_0^{\f12}+C\sum_{|\alpha|=0,2} \int_0^t\int\int\f{\nu(|v|)|\partial^{\alpha}(\wt{\mb{G}},\mb{P}_c F_2)|^2}{\mb{M_*}}
dxdvd\tau\\
\di + C(\chi_{\scriptscriptstyle T} +\sqrt\delta_0) \int_0^t\Big[\sum_{|\alpha|=1}\|\partial^\alpha(\phi,\psi,\omega)\|^2+\|\sqrt{\bar u_{1x}}(\phi,\psi_1,\omega)\|^2+\|\Pi_x\|^2+\|(\psi_{1xx},\phi_{xx})\|^2\Big]d\tau\\
\di +C(\chi_{\scriptscriptstyle T} +\sqrt\delta_0+\eta_0)\int_0^t\|(n_2,\Pi_{x\tau})\|^2d\tau.
\end{array}
\end{equation*}

\underline{Step 5. Estimation on $\|n_{2x}(t,\cdot)\|^2$.}

Now we estimate $\|n_{2x}\|^2.$ Applying $\partial_x$ to the equation \eqref{n2}, one has
\begin{equation}\label{n2x}
\begin{array}{ll}
\di n_{2xt}+(u_1n_2)_{xx}+ \Big(\f{\k_1(\t)}{R\t}\Pi_x\Big)_{xx}- \Big(\k_1(\t)(\f{n_2}{\r})_x\Big)_{xx}\\
\di
 =-\Big(\f{n_2}{\r}\int v_1 \mb{N}_{\mb{M}}^{-1}\big[\mb{P}_c(v_1\mb{M}_x) )\big]dv\Big)_{xx}-\Big(\int v_1 \mb{N}_{\mb{M}}^{-1}\big[\mb{P}_c(v_1 (\mb{P}_cF_2)_x)\big]dv\Big)_{xx}\\
\di
 - \Big(\int v_1 \mb{N}_{\mb{M}}^{-1}\big[\Pi_x\mb{G}_{v_1}\big]dv\Big)_{xx} -\Big(\int v_1 \mb{N}_{\mb{M}}^{-1}\Big[\partial_t(\mb{P}_c F_2)+(\f{\mb{M}}{\rho})_t~n_2
 -2Q(F_2,\mb{G})\Big]dv\Big)_{xx}.
\end{array}
\end{equation}
Multiplying the equation \eqref{n2x} by $n_{2x}$ and integrating the resulting equation with respect to $x,t$ yield that
\begin{equation}\label{n2x-2}
\begin{array}{ll}
\di \|n_{2x}\|^2(t)+\int_0^t\|(n_{2x},n_{2xx})\|^2d\tau\leq +C(\chi_{\scriptscriptstyle T} +\delta_0)\int_0^t\|(\Pi_x,n_2,\psi_x,\phi_x,\omega_x,\psi_{1xx},\phi_{xx})\|^2d\tau\\
\di +C(\chi_{\scriptscriptstyle T} +\delta_0)\int_0^t\int\int \f{\nu(|v|)|\big((\mb{P}_cF_2)_{x}, (\mb{P}_cF_2)_{t},\widetilde{\mb{G}}_{v_1},\wt{\mb{G}}_{v_1x},\widetilde{\mb{G}},\wt{\mb{G}}_x\big)|^2)}{\mb{M}_*}dxdvd\tau\\
\di +C\int_0^t\int\int \f{\nu(|v|)|(\mb{P}_cF_2)_{xx}, \mb{P}_cF_2)_{xt})|^2)}{\mb{M}_*}dxdvd\tau+C(\|n_{2x}\|^2(0)+\delta_0).
\end{array}
\end{equation}

By the equation \eqref{n21}, it holds that
\begin{equation}\label{n2xt}
\begin{array}{ll}
\di \int_0^t\|n_{2x\tau}\|^2d\tau=\int_0^t\|\int v_1F_{2xx}dv\|^2d\tau=\int_0^t\|\int v_1\big(\f{\mb{M}}{\rho}n_2+\mb{P}_c F_2\big)_{xx}dv\|^2d\tau\\
\di\quad \leq C\int_0^t\Big[\|n_{2xx}\|^2+\int\int\f{\nu(|v|)|(\mb{P}_cF_2)_{xx}|^2}{\mb{M}_*}dxdv\Big]d\tau\\
\di \qquad+C(\chi_{\scriptscriptstyle T} +\delta_0)\int_0^t\|(n_2,n_{2x},\phi_{xx},\psi_{xx},\omega_{xx},\phi_x,\psi_x,\omega_x)\|^2d\tau,
\end{array}
\end{equation}
and
\begin{equation}\label{n2tt}
\begin{array}{ll}
\di \int_0^t\|n_{2\tau\tau}\|^2d\tau=\int_0^t\|\int v_1F_{2x\tau}dv\|^2d\tau=\int_0^t\|\int v_1\big(\f{\mb{M}}{\rho}n_2+\mb{P}_c F_2\big)_{x\tau}dv\|^2d\tau\\
\di\quad \leq C\int_0^t\Big[\|n_{2x\tau}\|^2+\int\int\f{\nu(|v|)|(\mb{P}_cF_2)_{x\tau}|^2}{\mb{M}_*}dxdv\Big]d\tau\\
\di \qquad+C(\chi_{\scriptscriptstyle T} +\delta_0)\int_0^t\Big[\|(n_2,\phi_{x\tau},\psi_{x\tau},\omega_{x\tau})\|^2+\sum_{|\alpha^\prime|=1}\|\partial^\prime(\phi,\psi,\omega,n_2)\|^2\Big]d\tau.
\end{array}
\end{equation}
By \eqref{n2x-2}, \eqref{n2xt} and \eqref{n2tt}, one has
\begin{equation}\label{n2-he}
\begin{array}{ll}
\di \|n_{2x}(\cdot,t)\|^2+\int_0^t\Big[\|n_{2x}\|^2+\sum_{|\alpha|=2}\|\partial^\alpha n_2\|^2\Big] d\tau\leq C\|n_{20x}\|^2+C(\chi_{\scriptscriptstyle T} +\delta_0)\sum_{1\leq|\alpha|\leq2}\int_0^t\|\partial^\alpha (\phi,\psi,\omega)\|^2 d\tau\\
\di\qquad +C\sum_{|\alpha|=2}\int_0^t\int\int \f{\nu(|v|)|\partial^\alpha(\mb{P}_c F_2)|^2}{\mb{M}_*}dxdvd\tau+C
(\chi_{\scriptscriptstyle T} +\delta_0)\int_0^t\|(\Pi_x,n_2, n_{2\tau})\|^2d\tau
\\
\di\qquad +C(\chi_{\scriptscriptstyle T} +\delta_0)\int_0^t\int\int \f{\nu(|v|)|\big(\mb{P}_cF_2)_{x}, (\mb{P}_cF_2)_{\tau},\widetilde{\mb{G}}_{v_1},\wt{\mb{G}}_{v_1x},\widetilde{\mb{G}}\big)|^2)}{\mb{M}_*}dxdvd\tau.\\\end{array}
\end{equation}
By \eqref{phi-xtt-e-s}, it holds that
\begin{equation}\label{phixtt}
\begin{array}{ll}
\di \int_0^t\|\Pi_{x\tau\tau}\|^2d\tau \leq C\int_0^t\Big[\|n_{2\tau}\|^2+\int\int\f{\nu(|v|)|(\mb{P}_c F_2)_{\tau}|^2}{\mb{M}_*}dxdv\Big]d\tau+C(\chi_{\scriptscriptstyle T} +\delta_0)\int_0^t\|(n_2,\phi_\tau,\psi_\tau,\omega_\tau)\|^2d\tau.
\end{array}
\end{equation}

\underline{Step 6. Highest order estimates.}

 Finally, we estimate the highest
order derivatives, that is,  $\int_0^t\int\int \f{\nu(| v|)|\partial^\a \wt{\mb{G}}
|^2}{\mb{M}_*}d v dxd\tau$ with $|\a|=2$ and $ \|\phi_{xx}\|^2$ in (\ref{(4.55)}). To
do so, it is sufficient to study  $ \int\int \f{|\partial^\a
\wt{F}_1|^2}{\mb{M}_*}d v dx(|\a|=2)$ in view of \eqref{(4.8)-s}, \eqref{3.25-s} and \eqref{3.26-s}.
Applying
$\partial^\alpha(|\alpha|=2)$ on the Vlasov-Poisson-Boltzmann equation $\eqref{VPB}$, one has
\begin{equation}\label{F1-a2-s}
\begin{array}{ll}
\di(\partial^{\a}\wt F_1)_t+ v_1(\partial^{\a}\wt F_1)_x+\Pi_x\partial_{v_1}(\partial^\alpha F_2)+\partial^\alpha\Pi_x\partial_{v_1}F_2+\sum_{|\alpha^\prime|=1,\alpha^\prime\leq\alpha}C_{\alpha}^{\alpha^\prime}\partial^{\alpha^\prime}\Pi_x\partial_{v_1}\partial^{\alpha-\alpha^\prime}F_2\\
\di\qquad =
\partial^{\a}\mb{L}_\mb{M}\wt{\mb{G}}+2\partial^{\a}Q(\wt{\mb G},\wt{\mb G})
+\partial^{\a}\Big\{(\mb{L}_\mb{M}-\mb{L}_{\mb {M}^{S_1}})(\mb {G}^{S_1})+(\mb{L}_\mb{M}-\mb{L}_{\mb {M}^{S_3}})(\mb {G}^{S_3})\\[2mm]
\di\qquad
+2[Q(\wt{\mb G},\mb {G}^{S_1}+\mb {G}^{S_3})+Q(\mb {G}^{S_1}+\mb {G}^{S_3},\wt{\mb G})]
+2[Q(\mb {G}^{S_1},\mb {G}^{S_3})+Q(\mb {G}^{S_3},\mb {G}^{S_1})]
\Big\}.
\end{array}
\end{equation}
and
\begin{equation}\label{F2-a2-s}
\begin{array}{ll}
\di (\partial^\alpha F_2)_t+ v_1(\partial^\alpha F_2)_x+\Pi_x\partial_{v_1}(\partial^\alpha F_2)+\partial^\alpha\Pi_x\partial_{v_1}F_1\\
\di  \qquad+\sum_{|\alpha^\prime|=1,\alpha^\prime\leq\alpha}C_{\alpha}^{\alpha^\prime}\partial^{\alpha^\prime}\Pi_x\partial_{v_1}\partial^{\alpha-\alpha^\prime}F_1=\partial^\alpha
\mb{N}_\mathbf{M}(\mathbf{P}_cF_2)+2\partial^\alpha Q(F_2, \mathbf{G}).
\end{array}
\end{equation}
Multiplying \eqref{F1-a2-s} and \eqref{F2-a2-s} by $\frac{\partial^{\a}\wt F_1}{\mb{M}_*}=\frac{\partial^{\a}(\mb{M}-\mb{M}^{S_1}-\mb{M}^{S_3})}{\mb{M}_*}
+\frac{\partial^{\a}\wt{\mb{G}}}{\mb{M}_*}$ and $\frac{\partial^\alpha
F_2}{\mathbf{M}_*}=\frac{\partial^\alpha
(\f{n_2}\r\mathbf{M})}{\mathbf{M}_*}+\frac{\partial^\alpha
(\mathbf{P}_c F_2)}{\mathbf{M}_*}$, respectively,  we obtain
\begin{equation}\label{h1}
\begin{array}{l}
\quad\displaystyle (\frac{|\partial^\alpha
\wt F_1|^2}{2\mathbf{M}_*})_t-\frac{\partial^\alpha
\wt{\mathbf{G}}}{\mathbf{M}_*}\mb{L}_\mathbf{M}\partial^\alpha \wt{\mathbf{G}} =
\frac{\partial^\alpha
\wt F_1}{\mathbf{M}_*}\Big\{\Pi_x\partial_{v_1}(\partial^\alpha F_2)+\partial^\alpha\Pi_x\partial_{v_1}F_2+\sum_{|\alpha^\prime|=1,\alpha^\prime\leq\alpha}C_{\alpha}^{\alpha^\prime}\partial^{\alpha^\prime}\Pi_x\partial_{v_1}\partial^{\alpha-\alpha^\prime}F_2\\
\di\qquad +(\partial^{\a}\mb{L}_\mb{M}\wt{\mb{G}}-\mb{L}_\mathbf{M}\partial^\alpha \wt{\mathbf{G}})+2\partial^{\a}Q(\wt{\mb G},\wt{\mb G})
+\partial^{\a}\Big[(\mb{L}_\mb{M}-\mb{L}_{\mb {M}^{S_1}})(\mb {G}^{S_1})+(\mb{L}_\mb{M}-\mb{L}_{\mb {M}^{S_3}})(\mb {G}^{S_3})\\[2mm]
\di\qquad
+2[Q(\wt{\mb G},\mb {G}^{S_1}+\mb {G}^{S_3})+Q(\mb {G}^{S_1}+\mb {G}^{S_3},\wt{\mb G})]
+2[Q(\mb {G}^{S_1},\mb {G}^{S_3})+Q(\mb {G}^{S_3},\mb {G}^{S_1})]
\Big]\Big\}\\
\qquad{\displaystyle+\frac{\partial^{\a}(\mb{M}-\mb{M}^{S_1}-\mb{M}^{S_3})}{\mb{M}_*}\mb{L}_\mathbf{M}\partial^\alpha
\wt{\mathbf{G}}+( v_1\frac{|\partial^\alpha \wt F_1|^2}{2\mathbf{M}_*})_x},
\end{array}
\end{equation}
and
\begin{equation}\label{h2}
\begin{array}{l}
\quad\displaystyle (\frac{|\partial^\alpha
F_2|^2}{2\mathbf{M}_*})_t+\frac{\partial^\alpha
(\mathbf{P}_cF_2)}{\mathbf{M}_*}\mb{N}_\mathbf{M}\partial^\alpha(\mathbf{P}_cF_2) =\Pi_x\frac{\partial^\alpha
F_2}{\mathbf{M}_*}\partial_{v_1}(\partial^\alpha \wt F_1)+\Pi_x\frac{\partial^\alpha
F_2}{\mathbf{M}_*}\partial_{v_1}(\partial^\alpha F_{\a_1,\a_3}^S)\\
\di\qquad  +
\frac{\partial^\alpha
F_2}{\mathbf{M}_*}\Big[\partial^\alpha\Pi_x\partial_{v_1}F_1+\sum_{|\alpha^\prime|=1,\alpha^\prime\leq\alpha}C_{\alpha}^{\alpha^\prime}\partial^{\alpha^\prime}\Pi_x\partial_{v_1}\partial^{\alpha-\alpha^\prime}F_1+\sum_{|\alpha^\prime|=1,\alpha^\prime\leq \alpha}2C_{\alpha}^{\alpha^\prime}Q(\partial^{\alpha^\prime}\mathbf{P}_cF_2,\partial^{\alpha-\alpha^\prime}
\mathbf{M})\\
\di\qquad +2Q(\mathbf{P}_cF_2,\partial^\alpha \mathbf{M})\Big]+\frac{\partial^\alpha\big(
\f{\mathbf{M}}{\rho}n_2\big)}{\mathbf{M}_*}\mb{N}_\mathbf{M}\partial^\alpha
(\mathbf{P}_cF_2)+2\partial^\alpha
Q(F_2,\mathbf{G})\frac{\partial^\alpha F_2}{\mb{M}_*}+( v_1\frac{|\partial^\alpha F_2|^2}{2\mathbf{M}_*})_x,
\end{array}
\end{equation}
Adding \eqref{h1} and \eqref{h2} together and noting that
$$
\Pi_x\frac{\partial^\alpha
\wt F_1}{\mathbf{M}_*}\partial_{v_1}(\partial^\alpha F_2)+\Pi_x\frac{\partial^\alpha
F_2}{\mathbf{M}_*}\partial_{v_1}(\partial^\alpha \wt F_1)=(\Pi_x\frac{\partial^\alpha
\wt F_1\partial^\alpha F_2}{\mathbf{M}_*})_{v_1}+\Pi_x\frac{\partial^\alpha
\wt F_1\partial^\alpha F_2}{(\mathbf{M}_*)^2} (\mathbf{M}_*)_{v_1},
$$
and then integrating the resulting equation over $x, v, t$ imply that
\begin{equation*}\label{h3}
\begin{array}{ll}
\displaystyle \int\int \frac{|\partial^\alpha
\wt F_1|^2+|\partial^\alpha F_2|^2}{2\mathbf{M}_*}(x,v,t)dxdv- \int\int \frac{|\partial^\alpha
\wt F_{10}|^2+|\partial^\alpha F_{20}|^2}{2\mathbf{M}_*}dxdv\\
\di \quad-\int_0^t\int\int\Big[ \frac{\partial^\alpha
\wt{\mathbf{G}}}{\mathbf{M}_*}\mb{L}_\mathbf{M}\partial^\alpha \wt{\mathbf{G}}+\frac{\partial^\alpha
(\mathbf{P}_cF_2)}{\mathbf{M}_*}\mb{N}_\mathbf{M}\partial^\alpha(\mathbf{P}_cF_2) \Big]dxdvd\tau \\
\di =\int_0^t\int\int\Bigg\{ \Pi_x\frac{\partial^\alpha
\wt F_1\partial^\alpha F_2}{(\mathbf{M}_*)^2} (\mathbf{M}_*)_{v_1}+
\frac{\partial^\alpha
\wt F_1}{\mathbf{M}_*}\partial^\alpha\Pi_x\partial_{v_1}F_2+ \frac{\partial^\alpha
\wt F_1}{\mathbf{M}_*}\sum_{|\alpha^\prime|=1,\alpha^\prime\leq\alpha}C_{\alpha}^{\alpha^\prime}\partial^{\alpha^\prime}\Pi_x\partial_{v_1}\partial^{\alpha-\alpha^\prime}F_2\\
\di \quad+ \frac{\partial^\alpha
\wt F_1}{\mathbf{M}_*}(\partial^{\a}\mb{L}_\mb{M}\wt{\mb{G}}-\mb{L}_\mathbf{M}\partial^\alpha \wt{\mathbf{G}})+ \frac{\partial^\alpha
\wt F_1}{\mathbf{M}_*}\partial^{\a}\Big[(\mb{L}_\mb{M}-\mb{L}_{\mb {M}^{S_1}})(\mb {G}^{S_1})+(\mb{L}_\mb{M}-\mb{L}_{\mb {M}^{S_3}})(\mb {G}^{S_3})\\[2mm]
\di\quad
+2[Q(\wt{\mb G},\mb {G}^{S_1}+\mb {G}^{S_3})+Q(\mb {G}^{S_1}+\mb {G}^{S_3},\wt{\mb G})]
+2[Q(\mb {G}^{S_1},\mb {G}^{S_3})+Q(\mb {G}^{S_3},\mb {G}^{S_1})]
\Big]+ \frac{\partial^\alpha
\wt F_1}{\mathbf{M}_*}2\partial^\alpha Q(\wt{\mb{G}},\wt{\mathbf{G}})\\
\displaystyle\quad +\frac{\partial^{\a}(\mb{M}-\mb{M}^{S_1}-\mb{M}^{S_3})}{\mb{M}_*}\mb{L}_\mathbf{M}\partial^\alpha
\wt{\mathbf{G}}+\Pi_x\frac{\partial^\alpha
F_2}{\mathbf{M}_*}\partial_{v_1}(\partial^\alpha F_{\a_1,\a_3}^S)+\frac{\partial^\alpha
F_2}{\mathbf{M}_*}\sum_{|\alpha^\prime|=1,\alpha^\prime\leq\alpha}C_{\alpha}^{\alpha^\prime}\partial^{\alpha^\prime}\Pi_x\partial_{v_1}\partial^{\alpha-\alpha^\prime}F_1 \\
\di\quad  + \frac{\partial^\alpha
F_2}{\mathbf{M}_*}\partial^\alpha\Pi_x\partial_{v_1}F_1 + \frac{\partial^\alpha
F_2}{\mathbf{M}_*}\sum_{|\alpha^\prime|=1,\alpha^\prime\leq \alpha}2C_{\alpha}^{\alpha^\prime}Q(\partial^{\alpha^\prime}\mathbf{P}_cF_2,\partial^{\alpha-\alpha^\prime}
\mathbf{M})+ \frac{\partial^\alpha
F_2}{\mathbf{M}_*}2Q(\mathbf{P}_cF_2,\partial^\alpha \mathbf{M})\\
\di \quad +\frac{\partial^\alpha\big(
\f{\mathbf{M}}{\rho}n_2\big)}{\mathbf{M}_*}\mb{N}_\mathbf{M}\partial^\alpha
(\mathbf{P}_cF_2) +  2\partial^\alpha
Q(F_2,\mathbf{G})\frac{\partial^\alpha F_2}{\mb{M}_*}\Bigg\}dxdvd\tau.
\end{array}
\end{equation*}
Then we can get
\begin{equation*}
\begin{array}{l}
 \displaystyle \sum_{|\alpha|=2}\Big[\|\partial^\alpha\Pi_x(\cdot,t)\|^2+\int\int
\frac{|\partial^\alpha (\wt F_1,F_2)|^2}{2\mathbf{M}_*}(x,v,t)d v
dx\Big]+\frac{\widetilde{\sigma}}{2}\sum_{|\alpha|=2}\int_0^t\int\int\frac{\nu(|
v|)}{\mathbf{M}_*}|\partial^\alpha
(\wt{\mathbf{G}},\mb{P}_c F_2)|^2d v dxd\tau\\
\displaystyle\leq C \sum_{|\alpha|=2}\Big[\|\partial^\alpha\Pi_{x0}\|^2+\int\int
\frac{|\partial^\alpha (\wt F_{10},F_{20})|^2}{2\mathbf{M}_*}d v
dx\Big]+C(\chi_{\scriptscriptstyle T} +\delta_0)\sum_{|\alpha^\prime|=1}\int_0^t \|\partial^{\alpha^\prime}(\phi,\psi,\omega,n_2)\|^2d\tau
\\
\displaystyle +C(\chi_{\scriptscriptstyle T} +\delta_0)\int_0^t\int\int \Big[\f{\nu(|v|)|\big(\widetilde{\mb{G}}, \mb{P}_c F_2\big)|^2}{\mb{M}_*}+\sum_{|\alpha^\prime|=1}\frac{\nu(| v|)}{\mathbf{M}_*}|\partial^{\alpha^\prime}
(\wt{\mathbf{G}}, \mb{P}_c F_2)|^2\Big]d v
dxd\tau\\
\di +C(\eta_0+\chi_{\scriptscriptstyle T} +\delta_0)\sum_{|\alpha|=2}\int_0^t\|\partial^\alpha(\phi,\psi,\omega,\Pi_x,n_2)\|^2d\tau+C(\chi_{\scriptscriptstyle T} +\delta_0)\int_0^t\|(\Pi_{x\tau}, n_2)\|^2d\tau+ C\delta_0.
\end{array}
\end{equation*}
Note that by \eqref{n2-E-s}, \eqref{n2t-E-s} and \eqref{phixtt}, for
$|\alpha|=2$, it holds that
\begin{equation*}\label{phi-x-a}
\begin{array}{ll}
\di \int_0^t\|\partial^\alpha \Pi_x\|^2d\tau\leq
\int_0^t\|(n_{2x},n_{2\tau}, \Pi_{x\tau\tau})\|^2 d\tau\leq C\|n_{20}\|^2+ C(\chi_{\scriptscriptstyle T} +\delta_0)\int_0^t \|\Pi_x\|^2d\tau\\
\di
+C(\chi_{\scriptscriptstyle T} +\delta_0)\sum_{|\alpha|=1}\int_0^t\|\partial^\alpha(\phi,\psi,\omega)\|^2dxd\tau+C~\sum_{|\alpha^\prime|=1}\int_0^t\int
\int\f{\nu(|v|)|\partial^{\alpha^\prime}(\mb{P}_cF_2)|^2}{\mb{M}_*}
dvdxd\tau\\
\di \quad +C(\chi_{\scriptscriptstyle T} +\delta_0)\Big[\int_0^t\int
\int\f{\nu(|v|)|\mb{P}_cF_2|^2}{\mb{M}_*} dvdxd\tau+\int_0^t\int
\int\f{\nu(|v|)|(\widetilde{\mb{G}},
\widetilde{\mb{G}}_{v_1})|^2}{\mb{M}_*} dvdxd\tau\Big].
\end{array}
\end{equation*}
Therefore, collecting all the above estimates together, we can get Proposition \ref{Prop3.2}.

\section*{Appendix B. A Priori Estimates for Stability of Rarefaction Wave}
\renewcommand{\theequation}{B.\arabic{equation}}
\setcounter{equation}{0}

The proof of Theorem \ref{thm} is shown by the continuum argument for the local solution to the system \eqref{VPB1} or equivalently the system \eqref{F1-f}-\eqref{F2-pc-1}, which can be proved similarly as in \cite{Guo1, Yang-Zhao}. Therefore, to prove Theorem \ref{thm}, it is sufficient to close the a-priori assumption \eqref{assumption-ap} and verify the a priori estimates \eqref{Final-E} and the time-asymptotic behaviors of the solution.
By \eqref{F1-f} and \eqref{rare-s}, one has
\begin{equation}\label{sys-h-o}
\left\{
\begin{array}{ll}
\di \phi_t+\bar\r\psi_{1x}+\bar u_1\phi_x+\bar\r_x\psi_1+\bar u_{1x}\phi=-(\phi\psi_1)_x,\\[3mm]
\di \psi_{1t}+\bar u_1\psi_{1x}+\bar u_{1x}\psi_1 +\f23\omega_x
+\f{2\bar\theta}{3\bar\rho}\phi_x+\f23\r_x(\f\theta\r-\f{\bar\t}{\bar\r})
-\Pi_x\f{n_2}{\r}+\psi_1\psi_{1x}=-\f1\r\int v_1^2\mb{G}_x dv,\\[3mm]
\di \psi_{it}+\bar u_1\psi_{ix}=
-\f1\r\int v_1v_i\mb{G}_xdv,\quad i=2,3,\\[3mm]
\di \omega_{t}+\bar
u_1\omega_x+\t_x\psi_1+\f23u_{1x}\omega+\f23\bar\t\psi_{1x}+\f{\Pi_x}{\r}(\int
v_1F_2dv-u_1n_2)\\
\di\qquad\qquad\quad
=-\f1\r\int
v_1\f{|v|^2}2\mb{G}_xdv+\f1\r\sum_{i=1}^3 u_i\int v_1v_i\mb{G}_xdv.
\end{array} \right.
\end{equation}

In fact, by the a-priori assumption \eqref{assumption-ap}, one also has from the system \eqref{sys-h-o} that
\begin{equation*}
\|(\phi,\psi,\omega)\|_{L^\i_x}^2+
\|(\phi_{t},\psi_{t},\omega_{t})\|^2\leq C(\chi_{\scriptscriptstyle T} +\d_0)^2,
\end{equation*}
hence, one has
\begin{equation*}
\|(\r_{t},u_{t},\t_{t})\|^2\di
\leq C\|(\phi_{t},\psi_{t},\omega_{t})\|^2+ C \|(\bar\r_{t},\bar u_{t},\bar\t_{t})\|^2\leq C(\chi_{\scriptscriptstyle T} +\d_0)^2.
\end{equation*}
For $|\alpha|=2$, it follows from \eqref{macro}  and \eqref{macro-n2} that
\begin{equation*}
\|\partial^\a\left(\rho,\rho u,\rho(\t+\f{|u|^2}{2}),n_2\right)\|^2 \leq
C\int\int\f{|\partial^\a (F_1,F_2)|^2}{\mb{M}_*}dv dx\leq C(\chi_{\scriptscriptstyle T} +\delta_0)^2,
\label{(4.8)}
\end{equation*}
and
\begin{equation}
\begin{array}{ll}
\di\|\partial^\a(\r,u,\t)\|^2&\di\leq C\|\partial^\a\left(\rho,\rho
u,\rho(\t+\f{|u|^2}{2})\right)\|^2
+C\sum_{|\a^\prime|=1}\int|\partial^{\a^\prime}\left(\rho,\rho
u,\rho(\t+\f{|u|^2}{2})\right)|^4dx\\
 &\di\leq C(\chi_{\scriptscriptstyle T} +\delta_0)^2.
\end{array}
\label{3.25}
\end{equation}
Therefore, for $|\a|=2$, we have
\begin{equation}
\|\partial^\a(\phi,\psi,\omega,n_2)\|^2\leq
C(\chi_{\scriptscriptstyle T} +\delta_0)^2. \label{3.26}
\end{equation}
By \eqref{n21} and a-priori assumption \eqref{assumption-ap}, it holds that
\begin{equation}\label{n2t-ap}
\|n_{2t}\|^2\leq C\Big[\|n_{2x}\|^2+\int|n_2|^2|(\r_x,u_x,\t_x)|^2 dx+\int\int\f{|(\mb{P}_c F_2)_x|^2}{\mb{M}_*} dxdv\Big]\leq C(\chi_{\scriptscriptstyle T} +\delta_0)^2.
\end{equation}
By \eqref{assumption-ap}, \eqref{3.26} and \eqref{n2t-ap}, for $|\alpha^\prime|=1$, it holds that
\begin{equation*}
\begin{array}{ll}
\|\partial^{\alpha^\prime}(\phi,\psi,\omega,n_2)\|^2_{L^{\i}}\leq
C(\chi_{\scriptscriptstyle T} +\delta_0)^2.
\end{array}
\end{equation*}
By \eqref{phi-xt-e-s}, one has
\begin{equation}\label{phi-xt-ap}
\|\Pi_{xt}\|^2\leq C\Big[\|n_2\|^2+\int\int\f{|\mb{P}_c F_2|^2}{\mb{M}_*} dxdv\Big]\leq C(\chi_{\scriptscriptstyle T} +\delta_0)^2.
\end{equation}
By \eqref{phi-xt-e-s}, it holds that
\begin{equation}\label{phi-xtt-e}
\f12\Pi_{xtt}+\int v_1 F_{2t}dv=0.
\end{equation}
Thus, one has
\begin{equation}\label{phixtt-ap}
\|\Pi_{xtt}\|^2\leq C\Big[\|n_{2t}\|^2+\int|n_2|^2|(\r_t,u_t,\t_t)|^2 dx+\int\int\f{|(\mb{P}_c F_2)_t|^2}{\mb{M}_*} dxdv\Big]\leq C(\chi_{\scriptscriptstyle T} +\delta_0)^2.
\end{equation}
By \eqref{phi-xt-ap}, \eqref{n2t-ap}, \eqref{phixtt-ap} and \eqref{3.26}, it holds that
\begin{equation*}\label{phi-xt-in}
\|(\Pi_{xt}, \Pi_{xtt})\|_{L^\infty}^2\leq C\|(\Pi_{xt}, \Pi_{xtt})\|\|(n_{2t}, n_{2tt})\|\leq  C(\chi_{\scriptscriptstyle T} +\delta_0)^2.
\end{equation*}
Moreover, it holds that
\begin{equation*}
\begin{array}{ll}
 \di\|\int\f{|(\wt{\mb{G}},\mb{P}_c F_2)|^2}{\mb{M}_*}dv \|_{L^\infty_x} \leq
C\left(\int\int\f{|(\wt{\mb{G}},\mb{P}_c F_2)|^2}{\mb{M}_*}dv
dx\right)^{\f{1}{2}}\cdot\left(\int\int\f{|(\wt{\mb{G}},\mb{P}_c F_2)_x|^2}{\mb{M}_*}dv dx\right)^{\f{1}{2}}\\
\di \qquad\qquad\qquad\qquad \leq
C(\chi_{\scriptscriptstyle T} +\delta_0)^2.
\end{array}
\label{G-infty}
\end{equation*}
Furthermore,  for $|\a^\prime|=1$, it holds that
\begin{equation}
\begin{array}{ll}
 \di\|\int\f{|\partial^{\a^\prime}
(\mb{G}, \mb{P}_c F_2)|^2}{\mb{M}_*}d v\|_{L_x^{\i}}\leq
C\left(\int\int\f{|\partial^{\a^\prime} (\mb{G}, \mb{P}_c F_2)|^2}{\mb{M}_*}dv
dx\right)^{\f{1}{2}}\cdot\left(\int\int\f{|\partial^{\a^\prime}
(\mb{G}, \mb{P}_c F_2)
_x|^2}{\mb{M}_*}dv dx\right)^{\f{1}{2}}\\
\di\qquad\qquad \leq C\chi_{\scriptscriptstyle T} (\chi_{\scriptscriptstyle T} +\d_0)\le C(\chi_{\scriptscriptstyle T} +\delta_0)^2.
\end{array}
\label{(4.17+)}
\end{equation}
Finally, by noticing the facts that $F_1=\mb{M}+\mb{G}$ and
$F_2=\f{n_2}{\r}\mb{M}+\mb{P}_c F_2$  and  \eqref{3.25}
with $|\a|=2$, it holds that
\begin{equation*}
\begin{array}{ll}
\di\int\int\f{|\partial^{\alpha}(\mb{G},\mb{P}_cF_2)|^2}{\mb{M}_*}dv dx\di
\le C\int\int\f{|\partial^\alpha (F_1,F_2)|^2}{\mb{M}_*}dv dx
+C\int\int\f{|\partial^\alpha\mb{M}|^2+|\partial^\a\big(\f{n_2}{\rho}\mb{M}\big)|^2}{\mb{M}_*}dv dx\\[3mm]
\di
\leq C(\chi_{\scriptscriptstyle T} +\d_0)^2,
\end{array}
\end{equation*}
where in the last inequality we have used a similar argument as used for
\eqref{(4.17+)}.
We start from the lower order estimates. First, the entropy is defined by
\begin{equation*}
-\frac{3}{2}\rho S=\int\mathbf{M}\ln \mathbf{M}d v.
\end{equation*}
Multiplying the equation \eqref{VPB1} by $\ln\mathbf{M}$ and integrating
over $v$, it holds that
$$
(-\frac{3}{2}\rho S)_t+(-\frac{3}{2}\rho u_1 S)_x+\int \Pi_x\partial_{v_1} F_2\ln \mb{M}dv
+(\int v_1\mathbf{G}\ln\mathbf{M}d v)_x=\int v_1\mathbf{G}(\ln\mathbf{M})_xd v.
$$
Direct computations yields
$$
 S=-\frac{2}{3}\ln\rho+\ln (\frac{4\pi}{3}\theta)+1,\\
\qquad  p=\frac{2}{3}\rho\theta=k\rho^{\frac53}\exp(S).
$$
Denote
$$
\begin{array}{l}
\displaystyle \mathbf{X}=(\rho,\r u_1,\r u_2,\r u_3,\r (\t+\f{|u|^2}{2}))^t,\\
\displaystyle \mathbf{Y}=(\r u_1,\r u_1^2+p,\r u_1u_2,\r u_1u_3,\r
u_1(\t+\f{|u|^2}{2})+pu_1)^t.
\end{array}
$$
Then the conservation law \eqref{F1-f} can be rewritten as
$$
\mathbf{X}_t+\mathbf{Y}_x=
\left(
\begin{array}{c}
0\\
\displaystyle \frac43\mu(\theta) u_{1x}-\int v_1^2\G d v\\
\displaystyle \mu(\theta) u_{2x}-\int v_1 v_2\G d v\\
\displaystyle \mu(\theta) u_{3x}-\int v_1 v_3\G d v\\
\displaystyle \kappa(\theta)\theta_x+\frac43\mu(\theta)u_1u_{1x}+
\sum_{i=2}^3\mu(\theta)u_iu_{ix} -\int\frac12 v_1| v|^2\G
d v
\end{array}
\right)_x.
$$
We define an entropy-entropy flux pair $(\eta,q)$ around a local
Maxwellian $\mathbf{M}_{[\bar\rho,\bar u,\bar\theta]}$ as
$$
\left\{
\begin{array}{l}
\displaystyle \eta=\bar\theta\left\{-\frac32\rho
S+\frac32\bar\rho\bar{S}+\frac32\nabla_\mathbf{X}(\rho
S)|_{\mathbf{X}=\bar{\mathbf{X}}}\cdot(\mathbf{X}-\bar{\mathbf{X}})\right\},\\
\displaystyle q=\bar\theta\left\{-\frac32\rho u_1
S+\frac32\bar\rho\bar{u}_1\bar{S}+\frac32\nabla_\mathbf{X}(\rho
S)|_{\mathbf{X}=\bar{\mathbf{X}}}\cdot(\mathbf{Y}-\bar{\mathbf{Y}})\right\}.
\end{array}
\right.
$$
Here, we can compute that
$$
\displaystyle (\rho S)_{\rho}=S+\frac{1}{\theta}+\frac{|u|^2}{2\theta}-\frac53,\qquad
\displaystyle (\rho S)_{m_i}=-\frac{u_i}{\theta},~i=1,2,3,\qquad
\displaystyle (\rho S)_{E}=\frac{1}{\theta},
$$
and
$$
\left\{\begin{array}{l} \displaystyle
\eta=\frac{3}{2}\left\{\rho\theta-\bar\theta\rho
S+\rho[(\bar{S}-\frac53)\bar\theta
+\frac{|u-\bar u|^2}{2}]+\frac23\bar\rho\bar\theta\right\} =\rho\bar\t\Psi(\f{\bar\r}{\rho})+\f32\rho\bar\t\Psi(\f{\t}{\bar\t})+\f34\rho|u-\bar u|^2,\\
\displaystyle q=u_1\eta  +(u_1-\bar
u_1)(\rho\theta-\bar\rho\bar\theta).
\end{array}
\right.
$$
Then, for $\mathbf{X}$ in any closed bounded region in
$\sum=\{\mathbf{X}:\rho>0,\theta>0\}$, there exists a positive
constant $C_0$ such that
$$
C_0^{-1}|\mathbf{X}-\bar{\mathbf{X}}|^2\leq\eta\leq
C_0|\mathbf{X}-\bar{\mathbf{X}}|^2.
$$
Direct computations yield that
\begin{equation}\label{E1}
\begin{array}{l}
\displaystyle
\quad\eta_t+q_x+\f{2\bar\t\mu(\t)}{\t}\psi_{1x}^2+\f{3\bar\t\mu(\t)}{2\t}\sum_{i=2}^3\psi_{ix}^2+\f{3\bar\r\k(\t)}{2\t^2}\omega^2_x
-\Big[\nabla_{(\bar\rho,\bar{u},\bar{S})}\eta\cdot(\bar\rho,\bar{u},\bar{S})_t
+\nabla_{(\bar\rho,\bar{u},\bar{S})}q\cdot(\bar\rho,\bar{u},\bar{S})_x\Big]\\
\di =(\cdots)_x-\f{3\bar\t\k(\t)}{2\t^2}\omega_x\bar\t_x+\f{3\k(\t)}{2\t^2}(\omega_x+\bar\t_x)\bar\t_x\omega-\f{2\bar\t\mu(\t)}{\t}\psi_{1x}\bar u_{1x}+\f{2\mu(\t)}{\t}(\bar u_{1x}+\psi_{1x})\bar u_{1x} \omega\\
\di\quad +\f32\f{\bar\t\omega_x-\omega\bar\t_x}{\t^2}\int v_1\f{|v|^2}{2}\G dv -\f32\f{\bar\t\omega_x-\omega\bar\t_x}{\t^2}\sum_{i=1}^3u_i\int v_1v_i\G dv+\f32\f{\bar\t\psi_{1x}-\bar u_{1x}\psi_1}{\t}\int v_1^2\G dv\\
\di\quad +\f32\f{\bar\t}{\t}\sum_{i=2}^3u_{ix}\int v_1v_i\G
dv+\f32\Pi_x n_2\psi_1+\f{3}{2\t}\Pi_x\omega\int v_1\mb{P}_c F_2dv.
\end{array}
\end{equation}
First, under the a-priori assumption \eqref{assumption-ap}, one has
$$
\begin{array}{ll}
\di\quad -\left[\nabla_{(\bar\rho,\bar{u},\bar{S})}\eta\cdot(\bar\rho,\bar{u},\bar{S})_t
+\nabla_{(\bar\rho,\bar{u},\bar{S})}q\cdot(\bar\rho,\bar{u},\bar{S})_x\right]\\[2mm]
\di =\f32\rho \bar u_{1x} (u_1-\bar u_1)^2+\f23\rho\bar\t \bar
u_{1x}\Psi(\f{\bar\rho}\rho)+\rho\bar\t\bar
u_{1x}\Psi(\f{\t}{\bar\t})-\f32\rho\bar\t_x(u_1-\bar
u_1)(\f23\ln\f{\bar\r}{\r}-\ln\f{\t}{\bar\t})\\
\di  \geq C^{-1} \bar u_{1x}(\phi^2+\psi_1^2+\omega^2),
\end{array}
$$
for some positive constant $C$ and the convex function
$$
\Psi(s)=s-\ln s-1.
$$
Integrating \eqref{E1} with respect to $x,t$ over
$\mb{R}^1\times[0,t]$ yields that
\begin{equation}\label{E10}
\begin{array}{ll}
\di \int\eta(x,t)dx+\int_0^t\Big[\|(\psi_x,\omega_x)\|^2+
\|\sqrt{\bar u_{1x}}(\phi,\psi_1,\omega)\|^2\Big]dt\\
\di \leq
C\int\eta(x,0)dx+C|\int_0^t\int\Big[-\f{3\bar\t\k(\t)}{2\t^2}\omega_x\bar\t_x+\f{3\k(\t)}{2\t^2}(\omega_x+\bar\t_x)\bar\t_x\omega\Big]
dxd\tau|\\
\di \quad +C|\int_0^t\int\Big[-\f{2\bar\t\mu(\t)}{\t}\psi_{1x}\bar
u_{1x}+\f{2\mu(\t)}{\t}(\bar u_{1x}+\psi_{1x})\bar u_{1x} \omega\Big]
dxd\tau|\\
\di \quad
+C|\int_0^t\int\f32\f{\bar\t\omega_x-\omega\bar\t_x}{\t^2}\Big[\int
v_1\f{|v|^2}{2}\G dv-\sum_{i=1}^3u_i\int v_1v_i\G dv\Big]dxd\tau|\\
\di \quad +C|\int_0^t\int\Big[\f32\f{\bar\t\psi_{1x}-\bar
u_{1x}\psi_1}{\t}\int v_1^2\G
dv+\f32\f{\bar\t}{\t}\sum_{i=2}^3u_{ix}\int v_1v_i\G dv\Big]
dxd\tau|\\
\di\quad +C|\int_0^t\int\Big[\f32\Pi_x
n_2\psi_1+\f{3}{2\t}\Pi_x\omega\int v_1\mb{P}_c F_2dv\Big] dxd\tau|
:=C\int\eta(x,0)dx+\sum_{i=1}^5I_i.
\end{array}
\end{equation}
First, by integration by parts and the Cauchy inequality, it holds
that
\begin{equation}\label{I1}
\begin{array}{ll}
\di
I_1=C|\int_0^t\int\Big[\f{3\bar\t\k(\t)}{2\t^2}\omega\bar\t_{xx}+\big(\f{3\bar\t\k(\t)}{2\t^2}\big)_x\omega\bar\t_{x}+\f{3\k(\t)}{2\t^2}(\omega_x+\bar\t_x)\bar\t_x\omega\Big]
dxd\tau|\\
\di\quad \leq C\int_0^t\int
|\omega|\Big[|\bar\t_{xx}|+|\bar\t_x|^2+|\omega_x||\bar\t_x|\Big]dxd\tau\\
\di\quad \leq
C\int_0^t\|\omega\|_{L^\infty}\Big[\|\bar\t_{xx}\|_{L^1}+\|\bar\t_x\|_{L^2}^2+\|\omega_x\|_{L^2}\|\bar\t_x\|_{L^2}\Big]d\tau\\
\di \quad \leq C\d^{\f18} \int_0^t\|\omega_x\|_{L^2}^2d\tau+C
\d^{\f18}\Big[\int_0^t\|\sqrt\eta\|_{L^2}^{2}(1+\tau)^{-\f76}d\tau+1\Big].
\end{array}
\end{equation}
Similarly, one has
\begin{equation*}\label{I2}
I_2 \leq C\d^{\f18} \int_0^t\|(\omega_x,\psi_{1x})\|_{L^2}^2d\tau+C
\d^{\f18}\Big[\int_0^t\|\sqrt\eta\|_{L^2}^{2}(1+\tau)^{-\f76}d\tau+1\Big],
\end{equation*}
and
\begin{equation}\label{I3}
\begin{array}{ll}
\di I_3 \leq \s \int_0^t\|(\omega_x,\omega
\bar\t_{x})\|_{L^2}^2d\tau+C_\s\ \int_0^t\int\Big[|\int
v_1\f{|v|^2}{2}\G dv|^2+\sum_{i=1}^3|\int v_1v_i\G dv|^2\Big]
dxd\tau,
\end{array}
\end{equation}
with some small positive constant $\s>0$ to be determined and the positive constant $C_\s$ depending on $\s$.
Note that by \eqref{Pi}, it holds that
\begin{equation}\label{I3+}
\begin{array}{ll}
\di |\int v_1\f{|v|^2}{2}\G dv|=|\int v_1\f{|v|^2}{2}
\mb{L}_\mb{M}^{-1}[\mb{G}_t+\mb {P}_1(v_1\mb{G}_x)+\mb
{P}_1(\Pi_x\partial_{v_1}F_2)-2Q(\mb{G}, \mb{G})] dv| :=\sum_{i=1}^4I_{3i}.
\end{array}
\end{equation}
Choose the global Maxellian $\mb{M}_*=\mb{M}_{[\rho_*,u_*,\t_*]}$ such that
\begin{equation}\label{ma1}
\rho_*>0,\qquad \f12\t(t,x)<\t_*<\t(t,x),
\end{equation}
and
\begin{equation}\label{ma2}
|\rho(x,t)-\r_*|+|u(x,t)-u_*|+|\t(x,t)-\t_*|<\eta_0
\end{equation}
with $\eta_0$ being the small positive constant in Lemma \ref{Lemma 4.2}.
Then with such chosen $\mb{M}_*$, it holds that
\begin{equation}\label{I31}
\begin{array}{ll}
\di I_{31}=|\int v_1\f{|v|^2}{2} \mb{L}_\mb{M}^{-1}\mb{G}_t dv|\leq
\Big(\int\f{\nu(|v|)|\mb{L}_\mb{M}^{-1}\mb{G}_t|^2}{\mb{M}_*}
dv\Big)^{\f12} \Big(\int\mb{M}_*\nu(|v|)^{-1}(v_1\f{|v|^2}{2})^2
dv\Big)^{\f12}\\
 \di\qquad \leq C \Big(\int\f{\nu(|v|)^{-1}|\mb{G}_t|^2}{\mb{M}_*}
dv\Big)^{\f12},
\end{array}
\end{equation}
and
\begin{equation}\label{I32}
\begin{array}{ll}
\di I_{32}=|\int v_1\f{|v|^2}{2} \mb{L}_\mb{M}^{-1}[\mb
{P}_1(v_1\mb{G}_x)] dv|\leq C
\Big(\int\f{\nu(|v|)|\mb{L}_\mb{M}^{-1}[\mb
{P}_1(v_1\mb{G}_x)]|^2}{\mb{M}_{[\r_*,u_*,2\t_*]}}
dv\Big)^{\f12} \\[3mm]
 \di\quad ~~~ \leq C \Big(\int\f{\nu(|v|)^{-1}|\mb
{P}_1(v_1\mb{G}_x)|^2}{\mb{M}_{[\r_*,u_*,2\t_*]}} dv\Big)^{\f12}\leq
C \Big(\int\f{\nu(|v|)^{-1}|\mb{G}_x|^2}{\mb{M}_*} dv\Big)^{\f12}.
\end{array}
\end{equation}
Furthermore, one has
\begin{equation}\label{I33}
\begin{array}{ll}
\di I_{33}=|\int v_1\f{|v|^2}{2} \mb{L}_\mb{M}^{-1}[\mb
{P}_1(\Pi_x\partial_{v_1}F_2)] dv|\leq C|\Pi_x|
\Big(\int\f{\nu(|v|)|\mb{L}_\mb{M}^{-1}[\mb
{P}_1(v_1\mb{G}_x)]|^2}{\mb{M}_*}
dv\Big)^{\f12} \\
 \di\quad ~~ \leq C|\Pi_x| \Big(\int\f{\nu(|v|)^{-1}|\f{n_2}{\r}\mb
{M}_{v_1}+\partial_{v_1}(\mb{P}_c F_2)|^2}{\mb{M}_*}
dv\Big)^{\f12}\\
\di \leq C|\Pi_x||n_2|+C|\Pi_x|
\Big(\int\f{\nu(|v|)^{-1}|\partial_{v_1}(\mb{P}_c F_2)|^2}{\mb{M}_*}
dv\Big)^{\f12},
\end{array}
\end{equation}
and
\begin{equation}\label{I34}
\begin{array}{ll}
\di I_{34}=|\int v_1\f{|v|^2}{2} \mb{L}_\mb{M}^{-1}[4Q(\mb{G},
\mb{G})] dv|\leq C \Big(\int\f{\nu(|v|)|\mb{L}_\mb{M}^{-1}[Q(\mb{G},
\mb{G})]|^2}{\mb{M}_*}
dv\Big)^{\f12} \\
 \di\quad ~~ \leq C\Big(\int\f{\nu(|v|)^{-1}|Q(\mb{G},
\mb{G})|^2}{\mb{M}_*} dv\Big)^{\f12}\leq
C\Big(\int\f{\nu(|v|)|\mb{G}|^2}{\mb{M}_*}
dv\Big)^{\f12}\Big(\int\f{|\mb{G}|^2}{\mb{M}_*}
dv\Big)^{\f12} \\
\di \leq C\Big(\int\f{\nu(|v|)|\widetilde{\mb{G}}|^2}{\mb{M}_*}
dv\Big)^{\f12}\Big(\int\f{|\widetilde{\mb{G}}|^2}{\mb{M}_*}
dv\Big)^{\f12}+C|(\bar\t_x,\bar
u_{1x})|\Big(\int\f{\nu(|v|)|\widetilde{\mb{G}}|^2}{\mb{M}_*}
dv\Big)^{\f12}+C|(\bar\t_x,\bar u_{1x})|^2.
\end{array}
\end{equation}
Substituting \eqref{I31}-\eqref{I34} into \eqref{I3+} and then into \eqref{I3}, one can arrive at
\begin{equation}\label{I3-e}
\begin{array}{ll}
\di I_3 \leq (\s+C\delta^{\f18}) \int_0^t\|\omega_x \|_{L^2}^2d\tau+C
\d^{\f18}\Big[\int_0^t\|\sqrt\eta\|_{L^2}^{2}(1+\tau)^{-\f76}d\tau+1\Big]\\
\di \quad +C_\s
\int_0^t\int\int\f{\nu(|v|)^{-1}|(\mb{G}_t,\mb{G}_x)|^2}{\mb{M}_*}
dv dxd\tau +C\chi_{\scriptscriptstyle T}  \int_0^t\|(\Pi_x, n_2)\|^2d\tau\\
\di \quad+C\chi_{\scriptscriptstyle T}
\int_0^t\int\int \f{\nu(|v|)|\mb{P}_c F_2|^2}{\mb{M}_*}dv dxd\tau+C(\chi_{\scriptscriptstyle T} +\delta)\int_0^t\int\int\f{\nu(|v|)^{-1}|\widetilde{\mb{G}}|^2}{\mb{M}_*}
dv dxd\tau.
\end{array}
\end{equation}
 Similar estimates holds for $I_4$. Then $I_5$ can be estimated by
 \begin{equation}\label{I5}
\begin{array}{ll}
\di I_5\leq C\chi_{\scriptscriptstyle T}  \int_0^t\|(\Pi_x, n_2)\|^2d\tau+C\chi_{\scriptscriptstyle T}
\int_0^t\int\int \f{\nu(|v|)|\mb{P}_c F_2|^2}{\mb{M}_*}dv dxd\tau.
\end{array}
\end{equation}
Substituting the estimates for $I_i(i=1,2,3,4,5)$ in \eqref{I1}-\eqref{I5} into \eqref{E10} and applying Gronwall inequality yield  the first-step
lower order estimates
\begin{equation}\label{le}
\begin{array}{ll}
\di
\|(\phi,\psi,\omega)(\cdot,t)\|^2+\int_0^t\|(\psi_x,\omega_x)\|^2d\tau+\int_0^t
\|\sqrt{\bar u_{1x}}(\phi,\psi_1,\omega)\|^2d\tau\\
\di \leq C\|(\phi,\psi,\omega)(\cdot,0)\|^2+C \d^{\f18}+C
\int_0^t\int\int\f{\nu(|v|)^{-1}|(\mb{G}_t,\mb{G}_x)|^2}{\mb{M}_*}
dv
dxd\tau\\
\di
+C(\chi_{\scriptscriptstyle T} +\delta)\int_0^t\int\int\f{\nu(|v|)|(\widetilde{\mb{G}},\mb{P}_c F_2)|^2}{\mb{M}_*}
dv dxd\tau+ C\chi_{\scriptscriptstyle T}  \int_0^t\|(\Pi_x, n_2)\|^2d\tau.
\end{array}
\end{equation}
Next we want to get the estimation of $\|\phi_x\|^2$. By
$\eqref{sys-h}_1$ and $\eqref{sys-h}_2$, it holds that
$$
\begin{array}{l}
\displaystyle
\frac{4\mu(\theta)}{3\r\bar\rho}\phi_{xt}+\psi_{1t}+\f23\omega_x
+\f{2\bar\theta}{3\bar\rho}\phi_x+\bar u_1\psi_{1x}+\bar
u_{1x}\psi_1 +\f23\r_x(\f\theta\r-\f{\bar\t}{\bar\r})
-\Pi_x\f{n_2}{\r}+\psi_1\psi_{1x}\\
\di
=-\f{4}{3\r}(\f{\mu(\t)}{\bar\r})_x\phi_t-\f{4}{3\r}\Big[\mu(\t)\f{\bar u_1\phi_x+\bar
u_{1x}\phi+\bar\r_x\psi_1+(\phi\psi_1)_x}{\bar\r}\Big]_x +\f{4}{3\r}\big(\mu(\t)\bar u_{1x}\big)_x-\f1\r\int
v_1^2\G_x dv.
\end{array}
$$
Multiplying the above equation by $\phi_x$, one has
\begin{equation}\label{ff}
\begin{array}{l}
\displaystyle
\big(\frac{2\mu(\theta)}{3\r\bar\rho}\phi_x^2+\psi_1\phi_x\big)_t+\f{2\bar\t}{3\bar\r}\phi_x^2
=(\frac{2\mu(\theta)}{3\r\bar\rho})_t\phi_x^2+\f23\omega_x\phi_x+\bar
u_{1x}\phi_x\psi_1\\
\di +\psi_{1x}\big[\bar\r\psi_{1x}+\bar u_{1x}\phi+\bar
\r_x\psi_1+(\phi\psi_1)_x\big]-\f23\r_x\phi_x(\f\theta\r-\f{\bar\t}{\bar\r})
+\Pi_x\phi_x\f{n_2}{\r}+\psi_1\psi_{1x}\phi_x-\f{4}{3\r}(\f{\mu(\t)}{\bar\r})_x\phi_t\phi_x\\
\di
-\f{4}{3\r}\Big[\mu(\t)\f{\bar
u_1\phi_x+\bar
u_{1x}\phi+\bar\r_x\psi_1\phi_x+(\phi\psi_1)_x}{\bar\r}\Big]_x\phi_x+\f{4}{3\r}\big(\mu(\t)\bar u_{1x}\big)_x\phi_x-\f1\r\int
v_1^2\G_x dv\phi_x.
\end{array}
\end{equation}
Note that
\begin{equation}\label{pix}
\begin{array}{ll}
\di \int_0^t\int|\int v_1^2\G_x d v|^2dxd\tau\leq
C\sum_{|\alpha|=2}\int_0^t\int\int\frac{\nu(|
v|)}{\mathbf{M}_*}|\partial^\alpha \mathbf{G}|^2d v
dxd\tau+C\delta^{\f12}\\
\di
+C(\chi_{\scriptscriptstyle T} +\delta)\int_0^t\Big[\|(\Pi_x,n_2,\phi_x,\psi_x,\omega_x)\|^2+\sum_{|\alpha^\prime|=
1}\int\int\frac{\nu(|
v|)|(\widetilde{\mathbf{G}},\partial^{\alpha^\prime}\mb{G})|^2}{\mathbf{M}_*}
d v dx
\\
\di \qquad\qquad\qquad\qquad +\int\int\frac{\nu(|v|)|(\mb{P}_c F_2)_{v_1}|^2}{\mathbf{M}_*}
d v dx
\Big]d\tau,
\end{array}
\end{equation}
and
\begin{equation}\label{uw}
\begin{array}{ll}
\di |\int_0^t\int -\f{4}{3\r}\Big(\mu(\t)\f{\bar
u_1\phi_x}{\bar\r}\Big)_x\phi_xdxd\tau|=|\int_0^t\int \mu(\t)\f{\bar
u_1\phi_x}{\bar\r}\Big(\f{4\phi_x}{3\r}\Big)_xdxd\tau|\\
\di =|\int_0^t\int \mu(\t)\f{\bar
u_1\phi_x}{\bar\r}\Big(\f{4\phi_{xx}}{3\r}-\f{4\phi_x\rho_x}{3\r^2}\Big)dxd\tau|\\
\di =|\int_0^t\int \Big[-\mu(\t)\f{\bar
u_1\phi_x}{\bar\r}\f{4\phi_x\rho_x}{3\r^2}-\big(\f{4\mu(\t)\bar u_1}{3\r\bar\r}\big)_x\f{\phi_x^2}{2}\Big]dxd\tau|\leq C(\chi_{\scriptscriptstyle T} +\delta)\int_0^t\|\phi_x\|^2d\tau.
\end{array}
\end{equation}
Integrating the equation \eqref{ff} with respect to $x,t$, then using Cauchy inequality and \eqref{pix}-\eqref{uw} and choosing $\chi_{\scriptscriptstyle T} ,\delta$ suitably small, one can obtain
\begin{equation}\label{phi-x-E}
\begin{array}{l}
\displaystyle
\|\phi_x(\cdot,t)\|^2+\int_0^t\|\phi_x\|^2d\tau \leq
C\Big[\|\psi_1(\cdot,t)\|^2+\|(\phi_{0x},\psi_{10})\|^2+\delta^{\f12}\Big]+C\int_0^t\|(\psi_{1x},\omega_x)\|^2d\tau
\\
\quad\displaystyle+C(\sqrt\delta+\chi_{\scriptscriptstyle T} )\int_0^t\|(\Pi_x,n_2,\psi_{1xx},\phi_{xx})\|^2d\tau+C\sum_{|\alpha|=2}\int_0^t\int\int\frac{\nu(| v|)}{\mathbf{M}_*}|\partial^\alpha
\mathbf{G}|^2d v dxd\tau \\
\di \quad +C(\chi_{\scriptscriptstyle T} +\sqrt\delta)\int_0^t\|\sqrt{\bar u_{1x}}(\psi_1,\phi,\omega)\|^2d\tau+C(\delta+\chi_{\scriptscriptstyle T} )\int_0^t \int\int\frac{\nu(|v|)|(\mb{P}_c F_2)_{v_1}|^2}{\mathbf{M}_*}
d v dx
d\tau\\
\di \quad+C(\delta+\chi_{\scriptscriptstyle T} )\sum_{|\alpha^\prime|=
1}\int_0^t\int\int\frac{\nu(| v|)}{\mathbf{M}_*}|(\widetilde{\mathbf{G}},\partial^{\alpha^\prime}\mb{G})|^2 d v dxd\tau.
\end{array}
\end{equation}
Then we estimate $\|(\phi,\psi,\omega)_t\|^2$. For this, we use the
system \eqref{sys-h-o}. By the equation $\eqref{sys-h-o}_2$,  one has
\begin{equation*}\label{psi-t-E}
\begin{array}{l}
\displaystyle \int_0^t\| \psi_{1t}\|^2 d\tau \leq C\int_0^t\|(\phi_x,\psi_x,\omega_x)\|^2d\tau+C(\chi_{\scriptscriptstyle T} +\delta)\int_0^t\Big[\|\sqrt{\bar u_{1x}}(\psi_1,\phi,\omega)\|^2+\|n_2\|^2\Big]d\tau\\
\di \qquad\qquad\qquad +C\int_0^t\int\int\frac{\nu(| v|)}{\mathbf{M}_*}|\mathbf{G}_x|^2
d v dxd\tau.
\end{array}
\end{equation*}
Similar estimates hold for $\phi_t,\psi_{2t},\psi_{3t}$ and $\omega_t$.  Therefore, one can arrive at
\begin{equation}\label{t-E}
\begin{array}{l}
\displaystyle \int_0^t\|(\phi_{t},\psi_t,\omega_t)\|^2 d\tau \leq C\int_0^t\|(\phi_x,\psi_x,\omega_x)\|^2d\tau+C(\chi_{\scriptscriptstyle T} +\delta)\int_0^t\Big[\|\sqrt{\bar u_{1x}}(\psi_1,\phi,\omega)\|^2+\|n_2\|^2\Big]d\tau\\
\di \qquad\qquad\qquad +C\int_0^t\int\int\frac{\nu(| v|)}{\mathbf{M}_*}|\mathbf{G}_x|^2
d v dxd\tau +C(\chi_{\scriptscriptstyle T} +\delta)\int_0^t\int\int\frac{\nu(| v|)}{\mathbf{M}_*}|\mathbf{P}_cF_2|^2
d v dxd\tau .
\end{array}
\end{equation}
By \eqref{le}, \eqref{phi-x-E} and \eqref{t-E}, it holds that
\begin{equation*}\label{le1}
\begin{array}{ll}
\di
\|(\phi,\psi,\omega,\phi_x)(\cdot,t)\|^2+\sum_{|\alpha|=1}\int_0^t\|\partial^\alpha(\phi,\psi,\omega)\|^2d\tau+\int_0^t
\|\sqrt{\bar u_{1x}}(\phi,\psi_1,\omega)\|^2d\tau\\
\di \leq C\|(\phi,\psi,\omega,\phi_x)(\cdot,0)\|^2+C \d^{\f18} +C(\chi_{\scriptscriptstyle T} +\delta)\int_0^t\int\int\f{\nu(|v|)|(\widetilde{\mb{G}},\mb{P}_c F_2, (\mb{P}_c F_2)_{v_1})|^2}{\mb{M}_*}
dv dxd\tau\\
\di
+C\sum_{1\leq|\alpha|\leq2}
\int_0^t\int\int\f{\nu(|v|)^{-1}|\partial^\alpha\mb{G}|^2}{\mb{M}_*}
dv
dxd\tau+ C(\chi_{\scriptscriptstyle T} +\sqrt\delta) \int_0^t\|(\Pi_x, n_2,\psi_{1xx},\phi_{xx})\|^2d\tau.
\end{array}
\end{equation*}

Next we do the microscopic estimates for the Vlasov-Poisson-Boltzmann system.
Multiplying the equation \eqref{GE} and the equation \eqref{F2-pc-1} by
$\f{\widetilde{\mb{G}}}{\mb{M_*}}$ and $\f{\mb{P}_c
F_2}{\mb{M}_*}$ , respectively, one has
\begin{equation}\label{Gle0}
\begin{array}{ll}
\di \Big(\f{|\widetilde{\mb{G}}|^2}{2\mb{M}_*}\Big)_t-\f{\widetilde{\mb{G}}}{\mb{M}_*}\mb{L}_{\mb{M}}\widetilde{\mb{G}}=\Big\{-\f{3}{2\t}\mb{L}_{\mb{M}}^{-1}\Big[\mb{P}_1\big(v_1(\f{|v-u|^2}{2\t}\omega_x+v\cdot\psi_{x})\big)\mb{M}\Big]\\[3mm]
\di \qquad\qquad\qquad\quad
-\mb{P}_1(v_1\mb{G}_x)-\mb{P}_1(\Pi_x\partial_{v_1}F_2)+2Q(\mb{G},\mb{G})-\bar{
\mb{G}}_t\Big\}\f{\widetilde{\mb{G}}}{\mb{M}_*},
\end{array}
\end{equation}
and
\begin{equation}\label{Pc-F2-le}
\begin{array}{ll}
\di \Big(\f{|\mb{P}_c F_2|^2}{2\mb{M}_*}\Big)_t-\f{\mb{P}_c
F_2}{\mb{M}_*}\mb{N}_{\mb{M}}(\mb{P_c}F_2)=\Big[-v_1\partial_x
F_2-(\f{n_2}{\rho}\mb{M})_t\\
\di\qquad\qquad\qquad\qquad\qquad -\mb{P}_c(\Pi_x\partial_{v_1}F_1)+2Q(F_2,\mb{G})\Big]\f{\mb{P}_c F_2}{\mb{M}_*},
\end{array}
\end{equation}
By using the similar  methods in obtaining \eqref{FE}, one can derive from \eqref{Gle0} and \eqref{Pc-F2-le} that
\begin{equation*}\label{FE-1}
\begin{array}{ll}
\di \int\int \f{|(\widetilde{\mb{G}},\mb{P}_c
F_2)|^2}{\mb{M}_*}(x,v,t)dxdv+\|(\Pi_x,n_2)(\cdot,t)\|^2+\int_0^t\int\int\frac{\nu(|v|)|\big(\widetilde{\mb{G}},
\mb{P}_c
F_2\big)|^2}{\mb{M}_*}dvdxd\tau\\
\di \leq C\int\int \f{|(\widetilde{\mb{G}}, \mb{P}_c
F_2)|^2}{\mb{M}_*}(x,v,0)dxdv+C\|(\Pi_{0x},n_{20})\|^2+C\|(\phi_0,\psi_0,\omega_0)\|^2+
C\delta^{\f18}\\
\di +C(\chi_{\scriptscriptstyle T} +\delta+\eta_0)\int_0^t\|(\Pi_x,n_{2x},n_{2t})\|^2
d\tau+C(\chi_{\scriptscriptstyle T} +\delta) \int_0^t
\Big[\|n_2\|^2+\sum_{|\alpha|=1}\|\partial^\alpha(\phi,\psi,\omega)\|^2\Big]d\tau\\
\di +C\sum_{|\alpha^\prime|=1}\int_0^t\int\int\frac{\nu(|v|)|\partial^{\alpha^\prime} \mb{G}|^2}{\mb{M}_*}dvdxd\tau.
\end{array}
\end{equation*}

Similar to \eqref{Phi-E-s} and \eqref{n2-E-s}, one has
\begin{equation*}\label{Phi-E}
\begin{array}{ll}
\di \|\Pi_x\|^2(t)+\int_0^t\|(\Pi_x,n_2)\|^2d\tau\leq
C\|\Pi_x\|^2(0)\\
\di \quad
+C(\chi_{\scriptscriptstyle T} +\delta)\sum_{|\alpha|=1}\int_0^t\|\partial^\alpha(\phi,\psi,\omega)\|^2dxd\tau+C~\sum_{|\alpha^\prime|=1}\int_0^t\int
\int\f{\nu(|v|)|\partial^{\alpha^\prime}(\mb{P}_cF_2))|^2}{\mb{M}_*}
dvdxd\tau\\
\di \quad +C(\chi_{\scriptscriptstyle T} +\delta)\Big[\int_0^t\int
\int\f{\nu(|v|)|(\widetilde{\mb{G}},\mb{P}_cF_2)|^2}{\mb{M}_*} dvdxd\tau+\int_0^t\int
\int\f{\nu(|v|)| \widetilde{\mb{G}}_{v_1}|^2}{\mb{M}_*}
dvdxd\tau\Big]
\end{array}
\end{equation*}
and
\begin{equation*}\label{n2-E}
\begin{array}{ll}
\di \|n_{2}\|^2(t)+\int_0^t\|(n_{2x},n_2,\sqrt{\bar
u_{1x}}n_2)\|^2d\tau\leq C\|n_{20}\|^2+ C(\chi_{\scriptscriptstyle T} +\delta)\int_0^t \|\Pi_x\|^2d\tau\\
\di
+C(\chi_{\scriptscriptstyle T} +\delta)\sum_{|\alpha|=1}\int_0^t\|\partial^\alpha(\phi,\psi,\omega)\|^2dxd\tau+C~\sum_{|\alpha^\prime|=1}\int_0^t\int
\int\f{\nu(|v|)|\partial^{\alpha^\prime}(\mb{P}_cF_2))|^2}{\mb{M}_*}
dvdxd\tau\\
\di \quad +C(\chi_{\scriptscriptstyle T} +\delta)\Big[\int_0^t\int
\int\f{\nu(|v|)|\mb{P}_cF_2|^2}{\mb{M}_*} dvdxd\tau+\int_0^t\int
\int\f{\nu(|v|)|(\widetilde{\mb{G}}, \widetilde{\mb{G}}_{v_1})|^2}{\mb{M}_*}
dvdxd\tau\Big].
\end{array}
\end{equation*}
By \eqref{phi-xtt-e}, it holds that
\begin{equation*}\label{phixt}
 \int_0^t\|\Pi_{xt}\|^2d\tau \leq
C\int_0^t\Big[\|n_{2}\|^2+\int\int\f{\nu(|v|)|\mb{P}_c
F_2|^2}{\mb{M}_*}dxdv\Big]d\tau.
\end{equation*}
By the equation \eqref{n21}, one has
\begin{equation*}\label{n2t-E}
\begin{array}{ll}
\di \int_0^t\|n_{2t}\|^2d\tau=\int_0^t\|\int v_1F_{2x}dv\|^2d\tau=\int_0^t\|\int v_1\big(\f{\mb{M}}{\rho}n_2+\mb{P}_c F_2\big)_xdv\|^2d\tau\\
\di\quad \leq C\int_0^t\Big[\|n_{2x}\|^2+\int\int\f{\nu(|v|)|(\mb{P}_cF_2)_{x}|^2}{\mb{M}_*}dxdv\Big]d\tau+C(\chi_{\scriptscriptstyle T} +\delta)\int_0^t\|(n_2,\phi_x,\psi_x,\omega_x)\|^2d\tau.
\end{array}
\end{equation*}
In summary, collecting all the above lower order estimates and choosing suitably small $\chi_{\scriptscriptstyle T} $, $\delta$ and $\eta_0$, we arrive that
\begin{equation*}\label{LE-A}
\begin{array}{ll}
\di \|(\phi,\psi,\omega,\phi_x,\Pi_x,n_2)(\cdot,t)\|^2+\int\int\f{|(\widetilde{\mb{G}},\mb{P}_c F_2)|^2}{\mb{M}_*}(x,v,t) dxdv+\int_0^t\|\sqrt{\bar u_{1x}}(\phi,\psi_1,\omega)\|^2d\tau\\
\di +\sum_{|\alpha^\prime|=1}\int_0^t\|\partial^{\alpha^\prime}(\phi,\psi,\omega,n_2)\|^2d\tau+\int_0^t\|(\Pi_x, \Pi_{xt}, n_2)\|^2d\tau+\int_0^t\int\int\f{\nu(|v|)|(\widetilde{\mb{G}},\mb{P}_c F_2)|^2}{\mb{M}_*}dxdvd\tau\\
\di \leq C\|(\phi_0,\psi_0,\omega_0,\phi_{0x},\Pi_{0x},n_{20})\|^2+C\delta^{\f18}+C\int\int\f{|(\widetilde{\mb{G}}_0,\mb{P}_c F_{20})|^2}{\mb{M}_*}(x,v) dxdv\\[3mm]
\di +C(\chi_{\scriptscriptstyle T} +\delta)\int_0^t\Big[\|(\psi_{1xx},\phi_{xx})\|^2+\int\int\f{\nu(|v|)|(\widetilde{\mb{G}},  \mb{P}_cF_2)_{v_1}|^2}{\mb{M}_*} dxdv\Big]d\tau\\
\di +C\sum_{|\alpha|=2}\int_0^t\int\int \f{\nu(|v|)|\partial^\alpha\mb{G}|^2}{\mb{M}_*}dxdvd\tau+C\sum_{|\alpha^\prime|=1}\int_0^t\int\int \f{\nu(|v|)|\partial^{\alpha^\prime}(\mb{G}, \mb{P}_cF_2)|^2}{\mb{M}_*}dxdvd\tau.
\end{array}
\end{equation*}
The higher order estimates can be done similarly as in Appendix A and will be skipped for brevity.

\bigskip

\noindent\textbf{Acknowledgements:} The authors are grateful with the referee's helpful suggestions on the revision of the manuscript. \par
The research of H.-L. Li is partially supported by NNSFC grants No.11231006, 11671384, and 11225102, NNFC-RGC grant No. 11461161007 and Beijing New Century Baiqianwan Talent Project.
The research of Y. Wang is partially supported by the National Natural Sciences Foundation
of China No. 11671385, Youth Innovation Promotion Association of CAS and Young top-notch talent Program. The  research of T. Yang is supported by the General Research Fund of Hong Kong, CityU 11302215.  The  research of M. Y. Zhong is partially supported by the NNSFC grant No. 11301094 and the Beijing Postdoctoral Research Foundation No. 2014ZZ-96.

\noindent\textbf{Conflict of Interest:} The authors declare that they have no conflict of interest.

\end{CJK*}
\end{document}